\documentclass[11pt]{amsart} 
\usepackage{amsmath} 
\usepackage{amsxtra}
\usepackage{amscd} 
\usepackage{amsthm} 
\usepackage{amsfonts}
\usepackage{amssymb} 
\usepackage{eucal}

\newtheorem{cor}[subsection]{Corollary}
\newtheorem{lem}[subsection]{Lemma}
\newtheorem{prop}[subsection]{Proposition}
\newtheorem{propconstr}[subsection]{Proposition-Construction}

\newtheorem{conj}[subsection]{Conjecture}
\newtheorem{thm}[subsection]{Theorem}
\newtheorem{defn}[subsection]{Definition}

\theoremstyle{definition}

\theoremstyle{remark}

\newcommand{\propconstrref}[1]{Proposition-Construction~\ref{#1}}
\newcommand{\thmref}[1]{Theorem~\ref{#1}}
\newcommand{\secref}[1]{Sect.~\ref{#1}}
\newcommand{\lemref}[1]{Lemma~\ref{#1}}
\newcommand{\propref}[1]{Proposition~\ref{#1}}
\newcommand{\corref}[1]{Corollary~\ref{#1}}
\newcommand{\conjref}[1]{Conjecture~\ref{#1}}
\newcommand{\defnref}[1]{Definition~\ref{#1}}

\newcommand{\nc}{\newcommand}
\nc{\renc}{\renewcommand}
\nc{\ssec}{\subsection}
\nc{\sssec}{\subsubsection}
\nc{\on}{\operatorname}

\nc\ol{\overline}
\nc\wt{\widetilde}
\nc\tboxtimes{\wt{\boxtimes}}
\nc{\alp}{\alpha}
\nc{\hl}{\overset{\leftarrow}h}
\nc{\hr}{\overset{\rightarrow}h}

\nc{\ZZ}{{\mathbb Z}}
\nc{\NN}{{\mathbb N}}
\nc{\PP}{{\mathbb P}}
\nc{\FF}{{\mathbb F}}
\nc{\CC}{{\mathbb C}}
\nc{\OO}{{\mathbb O}}
\renc{\SS}{{\mathbb S}}
\nc{\DD}{{\mathbb D}}
\nc{\GG}{{\mathbb G}}
\renewcommand{\AA}{{\mathbb A}}
\nc{\Fq}{{\mathbb F}_q}
\nc{\Fqb}{\ol{{\mathbb F}_q}}
\nc{\Ql}{\ol{{\mathbb Q}_\ell}}
\nc{\id}{\text{id}}
\nc\X{\mathcal X}

\nc{\Hom}{\on{Hom}}
\nc{\Lie}{\on{Lie}}
\nc{\Loc}{\on{Loc}}
\nc{\Pic}{\on{Pic}}
\nc{\Bun}{\on{Bun}}
\nc{\IC}{\on{IC}}
\nc{\Aut}{\on{Aut}}
\nc{\rk}{\on{rk}}
\nc{\Sh}{\on{Sh}}
\nc{\Perv}{\on{Perv}}
\nc{\pos}{{\on{pos}}}
\nc{\Conv}{\on{Conv}}
\nc{\Sph}{\on{Sph}}
\nc{\Sym}{\on{Sym}}
\nc{\BunBb}{\overline{\Bun}_B}
\nc{\Buno}{\overset{o}{\Bun}}
\nc{\BunPb}{{\overline{\Bun}_P}}
\nc{\BunBM}{\overline{\Bun}_{B(M)}}
\nc{\BunPbw}{{\widetilde{\Bun}_P}}
\nc{\BunBP}{\widetilde{\Bun}_{B,P}}
\nc{\GUb}{\overline{G/U}}
\nc{\GUPb}{\overline{G/U(P)}}
\nc{\Dt}{\widetilde{\on{D}}}
\nc{\Pt}{\widetilde{\on{P}}}

\nc{\Hhom}{\underline{\on{Hom}}}
\nc\syminfty{\on{Sym}^{\infty}}
\nc\lal{\ol{\lambda}}
\nc\xl{\ol{x}}
\nc\thl{\ol{\theta}}
\nc\nul{\ol{\nu}}
\nc\mul{\ol{\mu}}
\nc\Sum\Sigma
\nc{\oX}{\overset{o}{X}{}}

\nc{\M}{{\mathcal M}}
\nc{\N}{{\mathcal N}}
\nc{\F}{{\mathcal F}}
\nc{\D}{{\mathcal D}}
\nc{\Q}{{\mathcal Q}}
\nc{\Y}{{\mathcal Y}}
\nc{\G}{{\mathcal G}}
\nc{\E}{{\mathcal E}}
\nc{\CalC}{{\mathcal C}}
\nc\Dh{\widehat{\D}}

\nc{\C}{{\mathcal C}}
\nc{\K}{{\mathcal K}}
\renewcommand{\H}{{\mathcal H}}
\renewcommand{\S}{{\mathcal S}}
\nc{\T}{{\mathcal T}}
\nc{\V}{{\mathcal V}}
\renc{\P}{{\mathcal P}}
\nc{\A}{{\mathcal A}}
\nc{\B}{{\mathcal B}}
\nc{\U}{{\mathcal U}}
\renewcommand{\L}{{\mathcal L}}
\nc{\Gr}{\on{Gr}}

\nc{\fA}{{\mathfrak A}}
\nc{\fP}{{\mathfrak P}}
\nc{\frn}{{\check{\mathfrak u}(P)}}
\nc{\p}{\mathfrak p}
\nc{\q}{\mathfrak q}
\nc\f{{\mathfrak f}}

\nc{\s}{{\mathfrak s}}
\nc\w{\text{w}}
\renewcommand{\r}{{\mathfrak r}}

\nc\mathi\iota
\nc\Spec{\on{Spec}}
\nc\Mod{\on{Mod}}
\nc{\tw}{\widetilde{\mathfrak t}}
\nc{\pw}{\widetilde{\mathfrak p}}
\nc{\qw}{\widetilde{\mathfrak q}}
\nc{\jw}{\widetilde j}

\nc{\grb}{\overline{\Gr}}
\nc{\I}{\mathcal I}

\nc{\lambdach}{\check\lambda}
\nc{\Lambdach}{\check\Lambda}
\nc{\much}{\check\mu}
\nc{\omegach}{\check\omega}
\nc{\nuch}{\check\nu}
\nc{\etach}{\check\eta}
\nc{\alphach}{\check\alpha}
\nc{\betach}{\check\beta}
\nc{\rhoch}{\check\rho}

\nc{\Hb}{\overline{\H}}

\nc{\Qb}{\ol{\Q}}
\nc{\db}{\ol{d}}
\nc{\yb}{\ol{y}}

\nc{\Dwk}{\on{D}^W(\Qb_k)}
\nc{\Dwkex}{\on{D}^W(\Qb_{k,ex})}
\nc{\Dwkexp}{\on{D}^W(\Qb_{k+1,ex})}
\nc{\Pwk}{\on{P}^W(\Qb_k)}
\nc{\Pwkex}{\on{P}^W(\Qb_{k,ex})}
\nc{\Pwkexp}{\on{P}^W(\Qb_{k+1,ex})}

\title[On a vanishing conjecture]
{On a vanishing conjecture appearing in the \\
geometric Langlands correspondence}

\author{D. Gaitsgory }

\address{Department of Mathematics, the University of Chicago,
5734 University Ave., Chicago IL 60637}

\email{gaitsgde@math.uchicago.edu}

\thanks{The author is a prize fellow at the Clay Mathematics
Institute}

\begin{document}

\maketitle

\section*{Introduction}

\ssec{}

This paper should be regarded as a sequel to \cite{FGV1}. In {\it loc.cit.}
it was shown that the geometric Langlands conjecture for $GL_n$ follows
from a certain vanishing conjecture. The goal of the present paper is to prove
this vanishing conjecture.

Let $X$ be a smooth projective curve over a ground field $k$.
Let $E$ be an $m$-dimensional local system on $X$, and let $\Bun_m$ be the moduli
stack of rank $m$ vector bundles on $X$.

The geometric Langlands conjecture says that to $E$ we can associate a
perverse sheaf $\F_E$ on $\Bun_m$, which is a Hecke eigensheaf with respect
to $E$. 

The vanishing conjecture of \cite{FGV1} says that for all integers $n<m$,
a certain functor $\on{Av}^d_E$, depending on $E$ and a parameter $d\in \ZZ^+$,
which maps the category $\on{D}(\Bun_n)$ to itself,
vanishes identically, when $d$ is large enough.

The fact that the vanishing conjecture implies the geometric Langlands conjecture
may be regarded as a geometric version of the converse theorem. Moreover, as will
be explained in the sequel, the vanishing of the functor $\on{Av}^d_E$ is 
analogous to the condition that the Rankin-Selberg convolution of $E$, viewed
as an $m$-dimensional Galois representation, and an automorphic form on $GL_n$ with
$n<m$ is well-behaved.

\medskip

Both the geometric Langlands conjecture and the vanishing conjecture can be
formulated in any of the sheaf-theoretic situations, e.g., $\Ql$-adic sheaves (when
char$(k)\neq \ell$), D-modules (when  char$(k)=0$),
and sheaves with coefficients in a finite field $\FF_\ell$ (again, when 
char$(k)\neq \ell$).

When the ground field is the finite field $\FF_q$ and we are working with $\ell$-adic
coefficients, it was shown in \cite{FGV1} that the vanishing conjecture can 
be deduced from Lafforgue's theorem that establishes the full Langlands 
correspondence for global fields of positive characteristic, cf. \cite{Laf}.

The proof that will be given in this paper treats the cases of various ground fields
and coefficients uniformly, and in particular, it will be independent of Lafforgue's
results. 

However, we will be able to treat only the case of characteristic $0$ coefficients,
or, more generally, the case of $\FF_\ell$-coefficients when $\ell$ is $>d$, where $d$
is the parameter appearing in the formulation of the vanishing conjecture. 

\ssec{}

Let us briefly indicate the main steps of the proof.

\medskip

First, we show that instead of proving that the functor $\on{Av}^d_E$ vanishes,
it is sufficient to prove that it is exact, i.e., that it maps perverse sheaves
to perverse sheaves. The $\{$ exactness $\}\to\{$ vanishing $\}$ implication is 
achieved by an argument involving the comparison
of Euler-Poincar\'e characteristics of complexes obtained by applying the functor
$\on{Av}^d_E$ for various local systems $E$ of the same rank.

\medskip

Secondly, we show that the functor $\on{Av}^d_E$ can be expressed in terms 
of the ``elementary'' functor $\on{Av}^1_E$ using the action of the symmetric
group $\Sigma_d$. (It is this step that does not allow to treat the case of 
$\FF_\ell$-coefficients if $\ell\leq d$.)

\medskip

Thirdly, we define a certain quotient triangulated category $\Dt(\Bun_n)$
of $\on{D}(\Bun_n)$ by ``killing'' objects that one can call degenerate.
(This notion of degeneracy is spelled out using what we call Whittaker 
functors.) 

The main properties of the quotient $\Dt(\Bun_n)$ are as follows:
(0) $\Dt(\Bun_n)$ inherits the perverse $t$-structure from $\on{D}(\Bun_n)$,
(1) the Hecke functors defined on $\on{D}(\Bun_n)$ descend to $\Dt(\Bun_n)$
and are exact, and (2) the subcategory of objects of $\on{D}(\Bun_n)$ that map 
to $0$ in $\Dt(\Bun_n)$ is orthogonal to cuspidal complexes.

\medskip

Next we show that properties (0) and (1) above and the irreducibility assumption on $E$
formally imply that the elementary functor $\on{Av}^1_E$ is exact on the quotient category.
From that, we deduce that the functor $\on{Av}^d_E$ is also exact modulo the
subcategory of degenerate sheaves.

\medskip

Finally, by induction on $n$ we show that $\on{Av}^d_E$ maps $\on{D}(\Bun_n)$
to the subcategory of cuspidal sheaves, and, using property (2) above, we deduce that 
once $\on{Av}^d_E$ is exact modulo degenerate sheaves, it must be exact.

\ssec{}

Let us now explain how the the paper is organized.

\medskip

In Sect. 1 we recall the formulation of the vanishing conjecture.
In addition, we discuss some properties of the Hecke functors.

In Sect. 2 we outline the proof of the vanishing conjecture,
parallel to what we did above. We reduce the proof to two statements:
one is \thmref{exactness of Av^1} which says that the functor
$\on{Av}^1_E$ is exact on the quotient category,
and the other is the existence
of the quotient category $\Dt(\Bun_n)$ with the desired properties.

In Sect. 3 we prove \thmref{exactness of Av^1}.

Sect. 4-8 are devoted to the construction of the quotient category
and verification of the required properties. Let us describe the main ideas
involved in the construction.

\medskip

We start with some motivation from the theory of automorphic 
functions, following \cite{PS} and \cite{Sha}.

Let $\K$ be a global field, and $\AA$ the ring of adeles.
Let $P$ be the mirabolic subgroup of $GL_n$.
It is well-known that there is an isomorphism between the space of
cuspidal functions on $P(\K)\backslash GL_n(\AA)$ and the space
of Whittaker functions on $N(\K)\backslash GL_n(\AA)$, where $N\subset GL_n$
is the maximal unipotent subgroup. Moreover, this isomorphism can be written
as a series of $n-1$ Fourier transforms along the topological group $\K\backslash\AA$.

\medskip

In \secref{Whittaker categories} and \secref{Whittaker functs} we develop
the corresponding notions in the geometric context.
For us, the space of functions on $P(\K)\backslash GL_n(\AA)$ is replaced
by the category $\on{D}(\Bun_n')$, and the space of Whittaker functions 
is replaced by a certain subcategory in $\on{D}(\Qb)$ 
(cf. \secref{Whittaker categories}, where the notation is introduced).

The main result of these two sections is that there exists an exact ``Whittaker''
functor $W:\on{D}(\Bun_n')\to \on{D}^W(\Qb)$. The exactness is guaranteed by an
interpretation of $W$ as a series of Fourier-Deligne transform functors.

\medskip

In \secref{cuspidality} we show that the kernel $ker(W)\subset \on{D}(\Bun_n')$
is orthogonal to the subcategory $\on{D}_{cusp}(\Bun_n')$ of cuspidal sheaves.

\medskip

In \secref{Hecke functors} we define the action of the Hecke functors
on $\on{D}(\Bun_n')$ and $\on{D}^W(\Qb)$, and show that the Whittaker functor $W$ 
commutes with the Hecke functors.
The key result of this section is \thmref{Hecke is right-exact}, which
says that the Hecke functor acting on $\on{D}^W(\Qb)$ is right-exact.
This fact ultimately leads to the desired property (1) above,
that the Hecke functor is exact on the quotient category.

\medskip

Finally, in \secref{construction of quotients} we define our quotient 
category $\Dt(\Bun_n)$. 

\ssec{Conventions}

In the main body of the paper we will be working over a ground field $k$
of positive characteristic $p$ (which can be assumed algebraically closed) 
and with $\ell$-adic sheaves.
All the results carry over automatically to the D-module context for
schemes over a ground field of characteristic $0$, where instead of
the Artin-Schreier sheaf we use the corresponding D-module "$e^x$"
on the affine line. This paper allows to treat the case of $\FF_\ell$ 
coefficients, when $\ell>d$ (cf. below) in exactly the same manner.

\medskip

We follow the conventions of \cite{FGV1} in everything related to
stacks and derived categories on them. In particular, for a stack $\Y$
{\it of finite type}, we will denote by $\on{D}(\Y)$ the corresponding 
bounded derived category of sheaves on $\Y$. If $\Y$ is of infinite type,
but has the form $\Y=\underset{i}\cup \Y_i$, where $\Y_i$ is an increasing family 
of open substacks of finite type (the basic example being $\Bun_n$),
$\on{D}(\Y)$ is by definition the inverse limit of $\on{D}(\Y_i)$.

Throughout the paper we will be working with the perverse 
t-structure on $\on{D}(\Y)$, and will denote by $\on{P}(\Y)\subset \on{D}(\Y)$ the abelian
category of perverse sheaves. 
For $\F\in \on{D}(\Y)$, we will denote by $h^i(\F)$ its perverse cohomology sheaves.

For a map $\Y_1\to \Y_2$ and
$\F\in \on{D}(\Y_2)$ we will sometimes write $\F|_{\Y_1}$ for the
$*$ pull-back of $\F$ on $\Y_1$.

\medskip

For a group $\Sigma$ acting on $\Y$ we will denote by $\on{D}^\Sigma(\Y)$ the
corresponding equivariant derived category. In most applications, the
group $\Sigma$ will be finite, which from now on we will assume.

If the action of $\Sigma$ on
$\Y$ is trivial, we have the natural functor of invariants
$\F\mapsto (\F)^\Sigma:\on{D}^\Sigma(\Y)\to \on{D}(\Y)$.
This functor is {\it exact} when we work with coefficients of
characteristic zero, or when the order of $\Sigma$ is co-prime
with the characteristic.

The exactness of this functor is crucial for this paper, and it is
the reason why we have to assume that $\ell>d$, since the finite groups in
question will be the symmetric groups $\Sigma_{d'}$, $d'\leq d$.

\ssec{Acknowledgments} 

I would like to express my deep gratitude to V.~Drinfeld for 
his attention and many helpful discussions. His ideas are
present in numerous places in this paper. In particular, the
definition of Whittaker functors, which is one of the main 
technical tools, follows a suggestion of his. 

I would also like to thank D.~Arinkin, A.~Beilinson, A.~Braverman, E.Frenkel, 
D.~Kazhdan, I.~Mirkovi\'c, V.~Ostrik, K.~Vilonen and V.~Vologodsky
for moral support and stimulating discussions, and especially my thesis 
adviser J.~Bernstein, who has long 
ago indicated the ideas that are used in the argument proving 
\thmref{exactness of Av^1}.

\section{The conjecture}    \label{the conjecture}

\ssec{}

We will first recall the formulation of the Vanishing Conjecture,
as it was stated in \cite{FGV1}.

Let $\Bun_n$ be the moduli stack of rank $n$ vector bundles on our
curve $X$.
Let $\Mod_n^d$ denote the stack classifying the data of
$(\M,\M',\beta)$, where $\M,\M'\in\Bun_n$, and $\beta$ is an embedding
$\M\hookrightarrow \M'$ as coherent sheaves, and the quotient $\M'/\M$
(which is automatically a torsion sheaf) has length $d$.

We have the two natural projections
$$\Bun_n  \overset{\hl} \longleftarrow \Mod_n^d \overset{\hr}\longrightarrow \Bun_n,$$
which remember the data of $\M$ and $\M'$, respectively.

Let $X^{(d)}$ denote the $d$-th symmetric power of $X$. We have a natural map 
$s:\Mod_n^d\to X^{(d)}$, which sends a triple
$(\M,\M',\beta)$ to the divisor of the map $\Lambda^n(\M)\to \Lambda^n(\M')$.
In addition, we have a smooth map $\s:\Mod_n^d\to\on{Coh}^d_0$, where $\on{Coh}^d_0$
is the stack classifying torsion
coherent sheaves of length $d$. The map $\s$ sends a triple as above to $\M'/\M$.

\medskip

Recall that to a local system $E$ on $X$, Laumon associated a perverse sheaf
$\L^d_E\in \on{P}(\on{Coh}^d_0)$. The pull-back $\s^*(\L^d_E)$ (which is
perverse up to a cohomological shift) can be described as follows:

Let $\overset{\circ}X{}^d$ denote the complement to the diagonal divisor in $X^{(d)}$.
Let $\overset{\circ}{\on{Mod}}{}^d_n$ denote the preimage of $\overset{\circ}X{}^d$ under
$s$, and let $\overset{\circ}s:\overset{\circ}{\on{Mod}}{}^d_n\to \overset{\circ}X{}^d$ be the
corresponding map. Unlike $s$, the map $\overset{\circ}s$ is smooth.
Finally, let $j$ denote the open embedding of $\overset{\circ}{\on{Mod}}{}^d_n$ into $\Mod_n^d$.

Consider the symmetric power of $E$ as a sheaf $E^{(d)}\in \on{D}(X^{(d)})$,
and let $\overset{\circ}E{}^{(d)}$ denote its restriction to $\overset{\circ}X{}^d$.
It is easy to see that $\overset{\circ}E{}^{(d)}$ is lisse.

We have:
\begin{equation} 
\s^*(\L^d_E)\simeq j_{!*}\left(\overset{\circ}s{}^*(\overset{\circ}E{}^{(d)})\right).
\end{equation}

\ssec{}

We introduce the averaging functor
$\on{Av}^d_E: \on{D}(\Bun_n)\to \on{D}(\Bun_n)$ as follows:
$$\F\in \on{D}(\Bun_n)\mapsto \hl_!\left(\hr{}^*(\F)\otimes \s^*(\L^d_E)\right)[nd].$$

Let us note immediately, that this functor is essentially Verdier self-dual, 
in the sense that
$$\DD(\on{Av}^d_E(\F))\simeq \on{Av}^d_{E^*}(\DD(\F)),$$
where $E^*$ is the dual local system. This follows from the fact that
the map $\s\times\hr:\Mod_n^d\to \on{Coh}_0^d\times \Bun_n$
is smooth of relative dimension $nd$, and the map $\hl$ is proper.

\medskip

The following conjecture was proposed in \cite{FGV1}:

\begin{conj} \label{vanishing conjecture}
Assume that $E$ is irreducible, of rank $>$ n. Then for
$d$, which is greater than $(2g-2)\cdot n\cdot \on{rk}(E)$, the functor
$\on{Av}^d_E$ is identically equal to zero.
\end{conj}

\ssec{}

Let us discuss some rather tautological reformulations of \conjref{vanishing conjecture}.
Consider the map $\hl\times \hr:\Mod^d_n\to \Bun_n\times \Bun_n$;
it is representable, but not proper, and set
$\K^d_E:=(\hl\times \hr)_!(\s^*(\L^d_E))\in \on{D}(\Bun_n\times \Bun_n)$.

Let $\M\in \Bun_n$ be a geometric point (corresponding to a morphism denoted
$\iota_\M:\on{Spec}(k)\to \Bun_n$), and let $\delta_\M\in \on{D}(\Bun_n)$
be $(\iota_\M)_!(\ol{\mathbb Q}_l)$. Note that since the stack $\Bun_n$ is not
separated, $\iota_\M$ need not be a closed embedding, therefore,
$\delta_\M$ is a priori a complex of sheaves.

\begin{lem} \label{reformulation as kernel}
The vanishing of the functor $\on{Av}^d_E$ is equivalent to each of the following
statements:

\smallskip

\noindent {\em (1)} 
For every $\M\in \Bun_n$, the object $\on{Av}^d_E(\delta_\M)\in \on{D}(\Bun_n)$ vanishes.

\smallskip

\noindent {\em (2)} 
The object $\K^d_E\in \on{D}(\Bun_n\times \Bun_n)$ 
vanishes.
\end{lem}

\begin{proof}

First, statements (1) and (2) above are equivalent:
For $\M$, the stalk of $\on{Av}^d_E(\delta_\M)$ at $\M'\in \Bun_n$
is isomorphic to the stalk of $\K^d_E$ at $(\M\times \M')\in 
\Bun_n\times \Bun_n$.

Obviously, \conjref{vanishing conjecture} implies statement (1).
Conversely, assume that statement (1) above holds. Let
$\on{Av}^{-d}_{E^*}$ be the (both left and right) adjoint functor of
$\on{Av}^{d}_{E}$; explicitly,
$$\on{Av}^{-d}_{E^*}(\F)=\hr_!\left(\hl{}^*(\F)\otimes \s^*(\L^d_{E^*})\right)[nd].$$

It is enough to show that $\on{Av}^{-d}_{E^*}$ identically vanishes. However,
by adjointness, for an object $\F\in \on{D}(\Bun_n)$, the stalk of $\on{Av}^{-d}_{E^*}(\F)$
at $\M\in \Bun_n$ is isomorphic to $\on{RHom}_{\on{D}(\Bun_n)}(\on{Av}^d_E(\delta_\M),\F)$.

\end{proof}

\ssec{}

The assertion of the above conjecture is a geometric analog of the 
statement that the Rankin-Selberg convolution $L(\pi,\sigma)$, where
$\pi$ is an automorphic representation of $GL_{n}$ and $\sigma$
is an irreducible $m$-dimensional Galois representation with 
$m>n$ has an analytic continuation and satisfies a functional equation.

More precisely, let $X$ be a curve over a finite field, and
$\K$ the corresponding global field. 
Then it is known that the double quotient $GL_n(\K)\backslash GL_n(\AA)/GL_n(\OO)$
can be identified with the set (of isomorphism classes) of points of the stack
$\Bun_n$.

By passing to the traces of the Frobenius, we have a function-theoretic version
of the averaging functor; let us denote it by $\on{Funct(\on{Av}^d_E)}$,
which is now an operator from the space of functions on 
$GL_n(\K)\backslash GL_n(\AA)/GL_n(\OO)$ to itself.
 
Let now $f_\pi$ be a spherical vector in some unramified automorphic
representation $\pi$ of $GL_n(\AA)$. One can show that
\begin{equation} \label{L function}
\underset{d\geq 0}\Sigma \on{Funct(\on{Av}^d_E)}(f_\pi)=
L(\pi,E)\cdot f_\pi,
\end{equation}
where the L-function $L(\pi,E)$ is regarded as a formal series in $d$.

The assertion of \conjref{vanishing conjecture} implies that the above
series is a polynomial of degree $\leq m\cdot n\cdot (2g-2)$. And this is
the same estimate as the one following from the functional equation,
which $L(\pi,E)$ is supposed to satisfy.
 
\ssec{}

In the rest of this section we will make several preparatory
steps towards the proof of \conjref{vanishing conjecture}.

Recall that the Hecke functor $\on{H}:\on{D}(\Bun_n)\to \on{D}(X\times \Bun_n)$
is defined using the stack $\H=\Mod_n^1$, as 
$$\F\mapsto (s\times \hl)_!(\hr{}^*(\F))[n].$$

In the sequel it will be important to introduce parameters in all
our constructions. Thus, for a scheme $S$, we have a similarly
defined functor
$$\on{H}_S:\on{D}(S\times\Bun_n)\to \on{D}(S\times X\times \Bun_n).$$

\medskip

For an integer $d$ let us consider the d-fold iteration 
$\on{H}_{S\times X^{d-1}}\circ...\circ\on{H}_{S\times X}\circ \on{H}_S$, denoted
$$\on{H}_S^{\boxtimes d}:\on{D}(S\times\Bun_n)\to \on{D}(S\times X^d\times \Bun_n).$$

\begin{prop} \label{symmetric group acts}
The functor $\on{H}_S^{\boxtimes d}$ maps $\on{D}(S\times\Bun_n)$ to
the equivariant derived category $\on{D}^{\Sigma_d}(S\times X^d\times\Bun_n)$,
where $\Sigma_d$ is the symmetric group acting naturally on $X^d$.
\end{prop}

\begin{proof}

In the proof we will suppress $S$ to simplify the notation.

Let $\on{ItMod}^d_n$ denote the stack of iterated modifications, i.e.,
it classifies the data of a pair of vector bundles $\M,\M'\in \Bun_n$
together with a flag
$$\M=\M_0\subset \M_1\subset...\subset \M_d=\M',$$
where each $\M_i/\M_{i-1}$ is a torsion sheaf of length $1$.

Let $\r$ denote the natural map $\on{ItMod}^d_n\to \Mod^d_n$, and
let $\wt{\hl}$ and $\wt{\hr}$ be the two maps from
$\on{ItMod}^d_n$ to $\Bun_n$ equal to $\hl\circ \r$ and $\hr\circ \r$,
respectively. We will denote by $\wt{s}$ the map 
$\on{ItMod}^d_n\to X^d$, which remembers the supports of the successive quotients
$\M_i/\M_{i-1}$.

\medskip

It is easy to see that
the functor $\F\mapsto \on{H}^{\boxtimes d}(\F)$ can be rewritten
as
\begin{equation} \label{hecke formula}
\F\mapsto (\wt{s}\times \wt{\hl})_!(\wt{\hr}{}^*(\F))[nd].
\end{equation}

\medskip

We will now introduce a stack intermediate between $\Mod^d_n$ and $\on{ItMod}^d_n$.
Consider the Cartesian product 
$$\on{IntMod}^d_n:=\Mod^d_n\underset{X^{(d)}}\times X^d.$$
Note that $\on{IntMod}^d_n$ carries a natural action of the symmetric group
$\Sigma_d$ via its action on $X^d$.
Let $\wt{\wt{\hl}}, \wt{\wt{\hr}}$ be the corresponding projections from
$\on{IntMod}^d_n$ to $\Bun_n$, and $\wt{\wt{s}}$ the map
$\on{IntMod}^d_n\to X^d$. All these maps are $\Sigma_d$-invariant.

We have a natural map $\r_{\on{Int}}:\on{ItMod}^d_n\to \on{IntMod}^d_n$.

\begin{lem}
The map $\r_{\on{Int}}$ is a small resolution of singularities.
\end{lem}

The proof of this lemma follows from the fact $\on{IntMod}^d_n$
is squeezed between $\on{ItMod}^d_n$ and $\Mod^d_n$, and the fact that
the map $\r:\on{ItMod}^d_n\to \Mod^d_n$ is known to be small from 
the Springer theory, cf. \cite{BM}.

Hence, the direct image of the constant sheaf on $\on{ItMod}^d_n$ under
$\r_{\on{Int}}$ is isomorphic to the intersection cohomology sheaf
$\on{IC}_{\on{IntMod}^d_n}$, up to a cohomological shift. 

Therefore, by the projection formula, the expression in \eqref{hecke formula} can 
be rewritten as
\begin{equation} \label{another hecke formula}
(\wt{\wt{s}}\times \wt{\wt{\hl}})_!\left(\wt{\wt{\hr}}{}^*(\F)\otimes
\on{IC}_{\on{IntMod}^d_n}\right) [-\on{dim}(\Bun_n)]
\end{equation}

However, since the map $\wt{\wt{\hr}}$ is $\Sigma_d$-invariant, and
$\on{IC}_{\on{IntMod}^d_n}$ is a $\Sigma_d$-equivariant object of
$\on{D}(\on{IntMod}^d_n)$, we obtain that 
$\wt{\wt{\hr}}{}^*(\F)\otimes \on{IC}_{\on{IntMod}^d_n}$ is naturally an object of
$\on{D}^{\Sigma_d}(\on{IntMod}^d_n)$. Similarly, since the map $\wt{\wt{\hl}}$
is $\Sigma_d$-invariant, the expression in \eqref{another hecke formula}
is naturally an object of $\on{D}^{\Sigma_d}(\Bun_n)$.

\end{proof}

\ssec{}

Let $\Delta(X)\subset X^i$ be the main diagonal. Obviously, the symmetric
group $\Sigma_i$ acting on $X^i$ stabilizes $\Delta(X)$. Hence, for
an object $\F\in \on{D}^{\Sigma_i}(S\times X^i\times \Bun_n)$, it makes sense
to consider 
$$\Hom_{\Sigma_i}(\rho,\F|_{S\times \Delta(X)\times \Bun_n})\in \on{D}(S\times X\times \Bun_n)$$
for various representations $\rho$ of $\Sigma_i$. In particular,
let us consider the following functor $\on{D}(S\times \Bun_n)\to 
\on{D}(S\times X\times \Bun_n)$ that sends $\F$ to
$$\Hom_{\Sigma_i}
(sign,\on{H}_S^{\boxtimes i}(\F)|_{S\times \Delta(X)\times \Bun_n}),$$
where $sign$ is the sign representation of $\Sigma_i$. 

The following has been established in \cite{FGV1}:

\begin{prop} \label{hecke}
The functor
$$\F\mapsto \Hom_{\Sigma_i}
(sign,\on{H}_S^{\boxtimes i}(\F)|_{S\times \Delta(X)\times \Bun_n})$$
is zero if $i>n$ and for $i=n$ it is canonically isomorphic to
$$\F\mapsto (\on{id}_S\times m)^*(\F)[n],$$
where $m:X\times \Bun_n\to \Bun_n$ is the multiplication map,
i.e., $m(x,\M)=\M(x)$. 
\end{prop}

\begin{proof}

Again, to simplify the notation we will suppress the scheme $S$.

Let $\Mod^{i,\Delta}_n$ denote the preimage of $\Delta(X)\subset X^i$
inside $\on{IntMod}^i_n$. Note that the symmetric group $\Sigma_i$
acts trivially on $\Mod^{i,\Delta}_n$, and the $*$-restriction 
$\on{IC}_{\on{IntMod}^i_n}|_{\Mod^{i,\Delta}_n}$
is a $\Sigma_i$-equivariant object of $\on{D}(\Mod^{i,\Delta}_n)$.

Note also that for $i=n$, $\Mod^{i,\Delta}_n$ contains $X\times \Bun_n$
as a closed subset via
$$(x,\M)\mapsto (\M,\M(x),x^i)\in \Mod^i_n\underset{X^{(i)}}\times X^i.$$

The following is also a part of the Springer correspondence, 
cf. \cite{BM}, Sect. 3:

\begin{lem}
The object 
$$\Hom_{\Sigma_i}
(sign,\on{IC}_{\on{IntMod}^i_n}|_{\Mod^{i,\Delta}_n})$$ is
zero if $i>n$, and for $i=n$ it is isomorphic to the constant sheaf on
$X\times \Bun_n\subset \Mod^{i,\Delta}_n$ cohomologically shifted by
$[\on{dim}(\Bun_n)+n]$.
\end{lem}

This lemma and the projection formula imply the proposition.

\end{proof}

\ssec{}   \label{Averaging as an iteration}

We will now perform manipulations analogous to the ones of 
\propref{symmetric group acts} and \propref{hecke} with the
averaging functor $\on{Av}^d_E$.

Let us observe that for $d=1$, the averaging functor can be described
as follows:
$$\on{Av}^1_E(\F)\simeq p_!(\on{H}(\F)\otimes q^*(E)),$$
where $p$ and $q$ are the projections $X\times \Bun_n\to \Bun_n$
and $X\times \Bun_n\to X$, respectively.

We introduce the functor $\on{ItAv}^d_E:\on{D}(\Bun_n)\to \on{D}(\Bun_n)$
as a $d$-fold iteration of $\on{Av}^1_E$.

\begin{prop} \label{iteratet functor}
The functor $\on{ItAv}^d_E$ maps $\on{D}(\Bun_n)$ to the equivariant
derived category $\on{D}^{\Sigma_d}(\Bun_n)$. 
\end{prop}

\begin{proof}

First, it is easy to see that $\on{ItAv}^d_E(\F)$ can be rewritten as
$$p_!(\on{H}^{\boxtimes d}(\F)\otimes q^*(E^{\boxtimes d})),$$
where $p,q$ are the two projections from $X^d\times \Bun_n$ to $\Bun_n$
and $X^d$, respectively.

Hence, the assertion that $\on{ItAv}^d_E(\F)$ naturally lifts to an object
of the equivariant derived category $\on{D}^{\Sigma_d}(\Bun_n)$ follows
from \propref{symmetric group acts}.

\end{proof}

The next assertion allows to express the functor $\on{Av}^d_E$ via
$\on{Av}^1_E$. This is the only essential place in the paper where we use
the assumption that we are working with characteristic zero coefficients.

\begin{prop}   \label{iteration}
There is a canonical isomorphism of functors 
$$\on{Av}^d_E(\F)\simeq (\on{ItAv}^d_E(\F))^{\Sigma_d}.$$
\end{prop}

\begin{proof}

The following lemma has been proved in the original paper
of Laumon (cf. \cite{La}):

\begin{lem}  \label{Laumon's sheaf}
The direct image ${\mathcal Spr}_E^d:=
\r_!(\wt{s}^*(E^{\boxtimes d}))\in \on{D}(\Mod^d_n)$ 
is naturally $\Sigma_d$-equivariant. Moreover,
$$\s^*(\L_E^d)\simeq ({\mathcal Spr}_E^d)^{\Sigma_d}.$$
\end{lem}

Using the projection formula and the lemma, we can rewrite
$\on{ItAv}^d_E(\F)$ as
$$\hl_!(\hr{}^*(\F)\otimes {\mathcal Spr}_E^d)[nd].$$
(It is easy to see that the $\Sigma_d$-equivariant structure
on $\on{ItAv}^d_E(\F)$, which arises from the last expression
is the same as the one constructed before.)

Using \lemref{Laumon's sheaf} we conclude the proof.

\end{proof}

\section{Strategy of the proof}  \label{strategy}

In this section we will reduce the assertion of \conjref{vanishing conjecture}
to a series of theorems, which will be proved in the subsequent sections.

\ssec{}

By induction we will assume that \conjref{vanishing conjecture}
holds for all $n'$ with $n'<n$. 
We will deduce \conjref{vanishing conjecture} for $n$
from the following weaker statement:

\begin{thm}  \label{exactness of Av}
Let $E$, $n$ and $d$ be as in \conjref{vanishing conjecture}. Then
the functor $\on{Av}_E^d:\on{D}(\Bun_n)\to \on{D}(\Bun_n)$
is exact in the sense of the perverse t-structure.
\end{thm}

First we will prove that \thmref{exactness of Av} implies
\conjref{vanishing conjecture}. In fact, we will give two proofs:
the one discussed below is somewhat simpler, but at some point it 
resorts to some non-trivial results from the classical
theory of automorphic functions. The second proof,
which is due to A.~Braverman, will be given in the Appendix.

\medskip

Thus, let us assume that \thmref{exactness of Av} holds. Using 
\lemref{reformulation as kernel}(1), to prove \conjref{vanishing conjecture}, 
it suffices to show that $\on{Av}_E^d(\F)=0$, whenever $\F$ is a perverse sheaf, 
which appears as a constituent in some $\delta_\M$ for $\M\in \Bun_n$.
Set $\F'=\on{Av}_E^d(\F)$. By \thmref{exactness of Av}, we know that $\F'$ is perverse.

\begin{lem}  \label{van Euler char}
To show that a perverse sheaf $\F'$ on a stack $\Y$ vanishes, 
it is sufficient to show that the Euler-Poincar\'e characteristics of 
its stalks $\F'_y$ at all $y\in \Y$ are zero.
\end{lem}

\begin{proof}
If $\F'\neq 0$, there exists a locally closed substack
$\Y^0\subset \Y$, such that $\F'|_{\Y^0}$ is a lisse sheaf,
up to a cohomological shift. But then the Euler-Poincar\'e characteristics of
$\F'$ on $\Y^0$ are obviously non-zero.
\end{proof}

Now we have the following assertion, which states that the 
Euler-Poincar\'e characteristics of $\on{Av}_E^d(\F)$ do not
depend on the local system.

\begin{lem}  \label{independence}
Let $E'$ be any other local system on $X$ (irreducible or not)
with $\on{rk}(E')=\on{rk}(E)$.
Then the pointwise Euler-Poincar\'e characteristics of
$\on{Av}_E^d(\F)$ and $\on{Av}_{E'}^d(\F)$ are the same for any
$\F\in \on{D}(\Bun_n)$.
\end{lem}

\begin{proof}

We will deduce the lemma from the following theorem of Deligne, cf. \cite{Ill}:

{\it Let $f:\Y_1\to \Y_2$ be a proper map of schemes, and let $\S$ and $\S'$ be two
objects of $\on{D}(\Y_1)$, which are \'etale-locally isomorphic.
Then the Euler-Poincar\'e characteristics of $f_!(\S)$ and $f_!(\S')$ at
all points of $\Y_2$ coincide.}

\medskip

We apply this theorem in the following situation:

$\Y_2=\Bun_n$, $\Y_1=\on{Mod}^d_n$, $\S:=\hr{}^*(\F)\otimes \s^*(\L^d_E)$,
$\S':=\hr{}^*(\F)\otimes \s^*(\L^d_{E'})$, and $f=\hl$.

The assertion of the lemma follows from the fact that $s^*(\L^d_E)$
and $s^*(\L^d_{E'})$ are \'etale-locally isomorphic, because $E$ and $E'$ are.

\end{proof}

Thus, it suffices to show that for our $\F\in \on{P}(\Bun_n)$ and {\it some}
local system $E'$ of rank equal to that of $E$, the Euler-Poincar\'e characteristics of
the stalks of $\on{Av}_{E'}^d(\F)$ vanish.

When we are working in the $\ell$-adic situation over a finite field, the 
required fact was established in \cite{FGV1}. In fact, in {\it loc.cit.}
we exhibited a local system $E'$, for which the functor $\on{Av}_{E'}^d$ was zero.
\footnote{This part of the argument will be replaced by a different one in
the Appendix.}

In particular, we obtain that in the $\ell$-adic situation over a finite field
the vanishing of the Euler-Poincar\'e characteristics takes place when $E'$ is the
trivial local system. 

\medskip

Using the fact that our initial perverse sheaf $\F$ was of geometric
origin, the standard reduction argument (cf. \cite{BBD}, Sect. 6.1.7)
implies the vanishing of the Euler-Poincar\'e characteristics for
the trivial local system in the setting of $\ell$-adic sheaves over any 
ground field, and, when the field equals ${\mathbb C}$, also 
for constructible sheaves with complex coefficients. 

By the Riemann-Hilbert 
correspondence, this translates to the required vanishing statement
in the setting of D-modules over ${\mathbb C}$, and, hence, over
any field of characteristic zero.

\ssec{}

Thus, from now on, our goal will be to prove \thmref{exactness of Av}.
In view of \propref{iteration}, a natural idea would be to show that the ``elementary''
functor $\on{Av}^1_E$ is exact. The latter, however, is false. 

Recall that $\on{Av}_E^1$ is a composition of $\on{H}:\on{D}(\Bun_n)\to \on{D}(X\times \Bun_n)$
followed by the functor $\F\mapsto p_!(q^*(E)\otimes \F):
\on{D}(X\times \Bun_n)\to \on{D}(\Bun_n)$.

As it turns out, the source of the non-exactness of $\on{Av}^1_E$ is the fact that
the Hecke functor $\on{H}$ is not exact, except when $n=1$. Therefore, we will 
first consider the latter case, which would be the prototype of the argument
in general.

\ssec{The case n=1}

Of course, the assertion of \conjref{vanishing conjecture} in this case 
is known, cf. \cite{Dr}. However, the proof we are giving below is completely different.

\medskip

First, let us note that it is indeed sufficient to show that the functor
$\on{Av}^1_E$ is exact: 

The exactness of $\on{Av}^1_E$ implies that the functor $\on{ItAv}^d_E$
is exact for any $d$. Since the coefficients of our sheaves are of characteristic $0$,
from \propref{iteration} we obtain that $\on{Av}^d_E$ is a direct summand of 
$\on{ItAv}^d_E$, and, therefore, is exact as well.

\medskip

To prove that $\on{Av}^1_E$ exact, it is enough to show that
for an irreducible perverse sheaf $\F$, $\on{Av}^1_E(\F)$
has no cohomologies above $0$ (because $\on{Av}^1_E$ is essentially
Verdier self-dual).

For $n=1$, $\Bun_n$ is the Picard stack $\on{Pic}$, and the Hecke functor
can be identified with the pull-back $\F\mapsto m^*(\F)[1]$, where
$m:X\times \on{Pic}\to \on{Pic}$ is the multiplication map.
We have:
$$\on{Av}^1_E(\F)\simeq p_!(m^*(\F)[1]\otimes q^*(E)),$$
where $p$ and $q$ are the two projections $X\times \on{Pic}$ to
$\on{Pic}$ an $X$, respectively.

Since the map $m$ is smooth, the sheaf $m^*(\F)[1]$ is also perverse and irreducible, 
\footnote{The fact
that we can control irreducibility under the Hecke functors is another
simplification of the $n=1$ case.}
and $m^*(\F)[1]\otimes q^*(E)$ is perverse. Since $p$ is a projection with 
$1$-dimensional fibers, it is enough to show that
$$h^1\left(p_!\left(m^*(\F)[1]\otimes q^*(E)\right)\right)=0.$$

We will argue by contradiction. If $\F^1=h^1\left(p_!\left(m^*(\S)[1]\otimes q^*(E)\right) 
\right)\neq 0$, by adjunction we have a surjective map
\begin{equation}  \label{surjection in abelian case}
m^*(\F)[1]\otimes q^*(E)\to p^*(\F^1)[1],
\end{equation}
which gives rise to a map
\begin{equation}  \label{adjunction in the abelian case}
m^*(\F)[1]\to E^*[1]\boxtimes \F^1.
\end{equation}

Since $E$ was assumed irreducible, sub-objects of the right-hand side of
\eqref{adjunction in the abelian case} are in bijection with 
sub-objects of $\F^1$. Therefore, since the map of \eqref{surjection in abelian case}
is surjective, so is the map in \eqref{adjunction in the abelian case}.
By the irreducibility of $\F$, it must, therefore, be an isomorphism. 

\medskip

We claim that this cannot happen if the rank of $E$ is greater than $1$.

Indeed, let us consider the pull-back
$$(\on{id}\times m)^*\left(m^*(\F)\right)[2]\in \on{P}(X\times X\times \on{Pic}).$$

On the one hand, we know that it is isomorphic to
$E^*[1]\boxtimes m^*(\F^1)[1]$. On the other hand, it is equivariant with
respect to the permutation group $\Sigma_2$ acting on $X\times X$.

\begin{lem}
Let $\S$ be an irreducible perverse sheaf on a variety of the form
$X\times X\times \Y$, which, on the one hand, is $\Sigma_2$-equivariant, and 
on the other hand, has a form $E[1]\boxtimes \S'$, where $E$ is an irreducible
local system, and $\S'\in \on{P}(X\times \Y)$.
Then $\S$ must be of the form $\S\simeq E[1]\boxtimes E[1]\boxtimes \S''$,
moreover the $\Sigma_2$-equivariant structure on $\S$ is the standard one on
$E[1]\boxtimes E[1]$ times some $\Sigma_2$-action on $\S''$.
\end{lem}

\begin{proof}

Let $q$ and $p$ be the projections from $X\times \Y$ to $X$ and $\Y$,
respectively. It is enough to show that
$$h^1\left(p_!(\S'\otimes q^*(E^*))\right)\neq 0.$$
For $i=1,2$ let $q_i$ be the projection $X\times X\times \Y\to X$
on the $i$-th factor, and let $p_i$ be the complementary projection
on $X\times \Y$. 
We have
$$E\boxtimes p_!(\S'\otimes q^*(E^*))\simeq
p_2{}_!(\S\otimes q_2^*(E^*)),$$
which, due to the $\Sigma_2$-equivariance assumption,
is isomorphic to $p_1{}_!(\S\otimes q_1^*(E^*))$, and the latter
has non-trivial cohomology in dimension $1$.

\end{proof}

Thus, from the lemma, we obtain that
$(\on{id}\times m)^*\left(m^*(\F)\right)$ has the from
$E^*\boxtimes E^*\boxtimes \F''$. Let us restrict
$(\on{id}\times m)^*\left(m^*(\F)\right)$ to the diagonal
$(\Delta\times\on{id}):X\times \on{Pic}\subset X\times X\times \on{Pic}$,
and take $\Sigma_2$ anti-invariants.

On the one hand, from \propref{hecke} (which is especially easy in the $n=1$ case) 
we know that for any $\F\in \on{D}(\on{Pic})$,
$$\Hom_{\Sigma_2}
\biggl(sign, (\on{id}\times m)^*\left(m^*(\F)\right)|_{X\times \on{Pic}}\biggr)=0.$$

But on the other hand,
$(\on{id}\times m)^*\left(m^*(\F)\right)|_{X\times \on{Pic}}\simeq
(E^*)^{\otimes 2}\boxtimes \F''$, and when we take $\Sigma_2$ anti-invariants
we obtain
$$\biggl(\on{Sym}^2(E^*)\boxtimes \Hom_{\Sigma_2}\left(sign,\F''\right)\biggr)
\oplus
\biggl(\Lambda^2(E^*)\boxtimes (\F'')^{\Sigma_2}\biggr).$$
Now, since $\on{rk}(E)>1$, neither $\Lambda^2(E^*)$ nor $\on{Sym}^2(E^*)$ are
$0$, therefore, the entire expression cannot vanish.

\ssec{}   \label{subcategory compatible with t-structure}

The key fact used in the above argument was that the Hecke 
functor, which in this case acts as $\F\mapsto m^*(\F)[1]$, is exact.

For $n\geq 1$, our approach will consist of making the Hecke functors exact by
passing to a quotient triangulated category.

Recall that if $\C$ is a triangulated category, and $\C'\subset \C$ is
a full triangulated subcategory, one can form a quotient $\C/\C'$.
This quotient is a triangulated category endowed with a projection
functor
$$\C\to \C/\C',$$
which is universal with respect to the property that it makes any arrow
$\S_1\to \S_2$ in $\C$, whose cone belongs to $\C'$, into an isomorphism.

Note that the inclusion of $\C'$ into $ker(\C\to \C')$ is not
necessarily an equivalence. Rather, $ker(\C\to \C')$ is the
full subcategory, consisting of objects, which appear as direct
summands of objects of $\C'$. 

\medskip

Suppose now that $\C$ is endowed with a t-structure. Let $\on{P}(\C)$ be
the corresponding abelian subcategory, and let $\C'\subset \C$ be as above.

\begin{defn}   \label{defn of compatibility}
We say that $\C'$ is compatible with the t-structure if

\smallskip

\noindent{\em (1)}
$\on{P}(\C'):=\on{P}(\C)\cap \C'$ is a Serre subcategory of $\on{P}(\C)$.
\footnote{We remind that a Serre subcategory of an abelian category is
a full subcategory stable under taking sub-objects and extensions.}

\smallskip

\noindent{\em (2)}
If an object $\S\in \C$ belongs to $\C'$, then so do its
cohomological truncations $\tau^{\leq 0}(\S)$ and $\tau^{>0}(\S)$.
\end{defn}

A typical way of producing categories $\C'$ satisfying this definition
is given by the following lemma:

\begin{lem}  \label{kernel is compatible}
Let $\C_1$, and $\C_2$ be two triangulated categories
endowed with t-structures. Let $F:\C_1\to \C_2$ be a functor,
which is t-exact. 
Then $\C'_1:=ker(F)\subset \C_1$ is compatible with the t-structure.
\end{lem}

The following proposition is in some sense a converse to the above lemma:

\begin{prop} \label{t-structure on quotient}
Let $\C$ be as above, and $\C'\subset \C$ be compatible with the t-structure.
Then the quotient category $\wt{\C}:=\C/\C'$ carries a canonical t-structure, such that

\smallskip

\noindent{\em (1)}
The projection functor $\C\to \wt{\C}$ is exact.

\smallskip

\noindent{\em (2)}
The abelian category $\on{P}(\wt{\C})$ identifies with
the Serre quotient $\on{P}(\C)/\on{P}(\C')$.

\end{prop}

\begin{proof}

Let $\S$ be an object of $\C$, viewed as an object of the quotient category
$\C/\C'$.
We say that it belongs to $\wt{\C}^{\leq 0}$ (resp., $\wt{\C}^{>0}$)
if $\tau^{>0}(\S)$ (resp., $\tau^{\leq 0}(\S)$) belongs
to $\C'$.

If $\S_1\to \S_2$ is a morphism, whose cone belongs to $\C'$, it is easy 
to see that $\S_1$ belongs to $\wt{\C}^{\leq 0}$ (resp., $\wt{\C}^{>0}$)
if and only if $\S_2$ does.

We have to check now that if $\S_1\in \wt{\C}^{\leq 0}$, and 
$\S_2\in \wt{\C}^{>0}$,
then $\Hom_{\wt{\C}}(\S_1,\S_2)=0$.

Indeed, with no restriction of generality, by applying the cohomological 
truncation functor, we can assume that
$\S_1$ is represented by an object of $\C$, which lies in $\C^{\leq 0}$,
and $\S_2$ is represented by an object, which belongs to $\C^{>0}$.

Each element of the $\Hom$ group can be represented by a diagram
$$\S_1 \leftarrow \S_3 \to \S_2,$$
where the cone of $\S_3\to \S_1$ belongs to $\C'$. Hence, this diagram
can be replaced by an equivalent one
$$\S_1 \leftarrow \tau^{\leq 0}(\S_3)\to \S_2,$$
where $\tau$ is the cohomological truncation.

But now, any map $\tau^{\leq 0}(\S_3)\to\S_2$ is zero already in $\C$,
since $\S_2\in \C^{>0}$.

\medskip

The projection $\C\to \wt{\C}$ is exact by construction. By the universal property
of the Serre quotient, we have a functor $\on{P}(\C)/\on{P}(\C')\to \on{P}(\wt{\C})$.
Again, by construction, this functor is surjective on objects, and to prove that it is
fully-faithful it is sufficient to show that for $\S_1,\S_2\in \on{P}(\C)$ a map
$\S_1\to \S_2$ is an isomorphism in $\on{P}(\wt{\C})$ if and only if its kernel and
cokernel belong to $\on{P}(\C')$. 

Let $\S$ denote the cone of this map, regarded as an object of $\C$. By assumption,
it belongs to $\C'$; therefore $h^0(\S)$ and $h^1(\S)$ both belong to $\C'$,
by \defnref{defn of compatibility}. But the above $h^0(\S)$ and $h^1(\S)$,
both of which are objects of $\C'\cap \on{P}(\C)=\on{P}(\C')$, are the kernel
and cokernel, respectively, of $\S_1\to \S_2$.

\end{proof}

%

\ssec{}   \label{properties}

Thus, our strategy will be to find an appropriate quotient category
of $\on{D}(\Bun_d)$.

More precisely, we will construct for every base $S$ a category
$\Dt(S\times \Bun_n)$, which is the quotient of $\on{D}(S\times \Bun_n)$
by a triangulated subcategory $\on{D}_{degen}(S\times \Bun_n)$, such that
$\on{D}_{degen}(S\times \Bun_n)$ is compatible with the perverse t-structure,
and such that the following properties will hold:

\medskip

\noindent Property 0:
The categories $\on{D}(S\times \Bun_n)$ inherit the standard 4 functors.
In other words, for a map of schemes $f:S_1\to S_2$ the four direct and inverse 
image functors $\on{D}(S_1\times \Bun_n)\rightleftarrows
\on{D}(S_2\times \Bun_n)$ preserve the corresponding subcategories,
and thus define the functors 
$\Dt(S_1\times \Bun_n)\rightleftarrows \Dt(S_2\times \Bun_n)$.
Moreover, the same is true for the Verdier duality functor on 
$\on{D}(S\times \Bun_n)$, and for the functor
$\on{D}(S)\times \on{D}(S\times \Bun_n)\to \on{D}(S\times \Bun_n)$,
given by the tensor product along $S$.

\medskip

\noindent Property 1:
The Hecke functor
$\on{H}_S:\on{D}(S\times \Bun_n)\to \on{D}(S\times X\times \Bun_n)$
preserves the corresponding triangulated subcategories, and the 
resulting functor 
$$\wt{\on{H}}_S:\Dt(S\times \Bun_n)\to \Dt(S\times X\times \Bun_n)$$
is exact.

\medskip

\noindent Property 2:
There exists an integer $d_0$ large enough such that the following holds:
if $\F_1\in \on{D}(\Bun_n)$ is supported on the connected component $\Bun^d_n$
with $d\geq d_0$
(cf. \secref{sect Hecke right exact} for our conventions regarding the connected
components of the stack $\Bun_n$), and is cuspidal (cf. \cite{FGV1} or
\secref{cuspidality} for the notion of cuspidality), and $\F_2\in \on{D}_{degen}(\Bun_n)$, 
then the $\Hom$ group $\Hom_{\on{D}(\Bun_n)}(\F_1,\F_2)$ vanishes.

\ssec{}   \label{functors on the quotient}

Assuming the existence of such a family of quotient categories,
we will now derive \thmref{exactness of Av}.

First, let us observe that the functor
$\on{Av}^1_E:\on{D}(\Bun_n)\to \on{D}(\Bun_n)$
descends to a functor $\wt{\on{Av}}^1_E:\Dt(\Bun_n)\to \Dt(\Bun_n)$.

Indeed, according to \secref{Averaging as an iteration}, the functor
$\on{Av}^1_E$ is the composition of 
$$\wt{\on{H}}:\Dt(\Bun_n)\to \Dt(X\times \Bun_n),$$
well-defined according to Property 1 above, followed by a functor
$\Dt(X\times \Bun_n)\to \Dt(\Bun_n)$ that sends
$\S\in \on{D}(X\times \Bun_n)$ to
$p_!(q^*(E)\otimes \S)$, which is well-defined due to Property 0
(the maps $p$ and $q$ here are as in \secref{Averaging as an iteration}.)

\medskip

The first step in the proof of \thmref{exactness of Av} is the following 
theorem, which is the key result of this paper. 
The proof will be given in the next section, and it mimics
the argument for the $n=1$ case, discussed above.

\begin{thm}  \label{exactness of Av^1}
The functor $\wt{\on{Av}}^1_E:\Dt(\Bun_n)\to \Dt(\Bun_n)$
is exact.
\end{thm}

\medskip

To state a corollary of \thmref{exactness of Av^1}, which we will actually
use in the proof of \thmref{exactness of Av}, we need to make some 
preparations.

Let $S$ be a base and $\Sigma$ a finite group acting on $S$. 
(Here it becomes important that the characteristic of the
coefficients of our sheaves is either $0$ or coprime with $|\Sigma|$.)
We define the category $\Dt^{\Sigma}(S\times \Bun_n)$ as the quotient of
the equivariant derived category $\on{D}^{\Sigma}(S\times \Bun_n)$ by
the triangulated subcategory $\on{D}^{\Sigma}_{degen}(S\times \Bun_n)$ equal to
the preimage of $\on{D}_{degen}(S\times \Bun_n)$ under the forgetful functor
$$\on{D}^{\Sigma}(S\times \Bun_n)\to \on{D}(S\times \Bun_n).$$
The quotient acquires a t-structure, according to \lemref{kernel is compatible}
and \propref{t-structure on quotient}.

Let $\Pt^\Sigma(S\times \Bun_n)$ be the corresponding abelian subcategory
in $\Dt^{\Sigma}(S\times \Bun_n)$. By construction, this is the quotient
of $\on{P}^\Sigma(S\times \Bun_n)$ by a Serre subcategory consisting of
objects, whose image in $\Pt(S\times \Bun_n)$ is zero.

In the applications we will take $S=X^d$ and $\Sigma$ to be the symmetric group
$\Sigma_d$.

\medskip

Assume now that the action of $\Sigma$ on $S$ is actually trivial.
Then we have the functor of invariants 
denoted $\F\mapsto (\F)^\Sigma$ from $\on{D}^{\Sigma}(S\times \Bun_n)$ to
$\on{D}(S\times \Bun_n)$. We claim that it descends to a well-defined
functor $\Dt^{\Sigma}(S\times \Bun_n)\to \Dt(S\times \Bun_n)$. Indeed,
for an object $\F\in\on{D}^{\Sigma}(S\times \Bun_n)$, the image
of $(\F)^\Sigma$ under the forgetful functor
$\on{D}^{\Sigma}(S\times \Bun_n)\to \on{D}(S\times \Bun_n)$ is a 
direct summand of the image of $\F$. (In fact, in the case of a trivial action,
every object of $\Dt^\Sigma(S\times
\Bun_n)$ can be canonically written as $\underset{\rho}\oplus\, \S_\rho\otimes \rho$, 
where $\rho$ runs over the set of irreducible representations of $\Sigma$,
and $\S_\rho$ is an object of $\Dt(S\times \Bun_n)$.)

\medskip

That said, first, from \propref{symmetric group acts} we obtain that 
the functor $$\on{H}_S^{\boxtimes d}:\on{D}(S\times \Bun_n)\to
\on{D}^{\Sigma_d}(S\times X^d\times \Bun_n)$$ gives rise to a functor
$\wt{\on{H}}_S^{\boxtimes d}:\Dt(S\times \Bun_n)\to
\Dt^{\Sigma_d}(S\times X^d\times \Bun_n)$.

Secondly, from \propref{hecke} we obtain that the functor
$$(\on{id}_S\times m)^*:\on{D}(S\times \Bun_n)\to 
\on{D}(S\times X\times \Bun_n)$$ descends to a well-defined
functor $\Dt(S\times \Bun_n)\to \Dt(S\times X\times \Bun_n)$,
which is isomorphic to
\begin{equation} \label{hecke on the quotient}
\S\mapsto 
\Hom_{\Sigma_i}\biggl(sign,
\wt{\on{H}}_S^{\boxtimes i}(\S)|_{S\times \Delta(X)\times \Bun_n}\biggr)[-n],
\end{equation}
for $i=n$, and for $i>n$ the latter functor is zero.

Thirdly, from \propref{iteration}, we obtain that 
the functor $\on{ItAv}^d_E:\on{D}(\Bun_n)\to \on{D}^{\Sigma_d}(\Bun_n)$
gives rise to a functor $\wt{\on{ItAv}}^d_E:\Dt(\Bun_n)\to \Dt^{\Sigma_d}(\Bun_n)$

And finally, we obtain that the functor $\on{Av}^d_E:\on{D}(\Bun_n)\to \on{D}(\Bun_n)$
gives rise to a well-defined functor $\wt{\on{Av}}^d_E:\Dt(\Bun_n)\to \Dt(\Bun_n)$
with
\begin{equation} \label{iteration on the quotient}
\wt{\on{Av}}^d_E(\S)\simeq (\wt{\on{ItAv}}^d_E(\S))^{\Sigma_d}.
\end{equation}

Now \thmref{exactness of Av^1} implies the following:

\begin{cor}  \label{iteration is exact}
The functor $\wt{\on{Av}}_E^d:\Dt(\Bun_n)\to \Dt(\Bun_n)$ is exact. 
\end{cor}

\begin{proof}   

\thmref{exactness of Av^1} readily implies that the functor
$\wt{\on{ItAv}}^d_E$ is exact.

Since the functor $\F\mapsto (\F)^{\Sigma_d}:\on{D}^{\Sigma_d}(\Bun_n)\to 
\on{D}(\Bun_n)$ is exact (which follows from our assumption on the 
characteristic of the coefficients),
we obtain that the same is true for the 
corresponding functor $\Dt^{\Sigma_d}(\Bun_n)\to \Dt(\Bun_n)$.

Hence, the assertion follows from \eqref{iteration on the quotient}.

\end{proof}

\ssec{}

We proceed with the proof of \thmref{exactness of Av} modulo the existence of
the categories $\Dt(S\times \Bun_n)$ and \thmref{exactness of Av^1}, and the induction
hypothesis that \conjref{vanishing conjecture} holds for all $n'<n$. 
The following assertion is essentially borrowed from \cite{FGV1}:

\begin{lem}  \label{Av is cuspidal}
For any $\F\in \on{D}(\Bun_n)$, the object
$\on{Av}^d_E(\F)$ is cuspidal, provided that $d> (2g-2)\cdot n\cdot \on{rk}(E)$. 
\end{lem}

\begin{proof}

We have to show that the constant term functors
$\on{CT}^n_{n_1,n_2}(\on{Av}^d_E(\F))$ all vanish.

However, it was shown in \cite{FGV1}, Lemma 9.8, that
$\on{CT}^n_{n_1,n_2}(\on{Av}^d_E(\F))$ is an extension 
of objects of the form
$$(\on{Av}^{d_1}_E\boxtimes \on{Av}^{d_2}_E)(\on{CT}^n_{n_1,n_2}(\F)),$$
for all possible $d_1,d_2\geq 0$, $d_1+d_2=d$,
where $\on{Av}^{d_1}_E\boxtimes \on{Av}^{d_2}_E$ denotes the corresponding
functor $\on{D}(\Bun_{n_1}\times \Bun_{n_2})\to \on{D}(\Bun_{n_1}\times \Bun_{n_2})$.

\medskip

However, for every pair $d_1,d_2$ as above, at least one of the parameters
satisfies $d_i>(2g-2)\cdot n_i\cdot \on{rk}(E)$. Hence, the corresponding
functor $\on{Av}^{d_i}_E:\on{D}(\Bun_{n_i})\to \on{D}(\Bun_{n_i})$ vanishes
by the induction hypothesis.

\end{proof}

Now we are ready to finish the proof of \thmref{exactness of Av}.
Since $\on{Av}^d_E$ is essentially Verdier self-dual, it is
enough to show that $\on{Av}^d_E$ is right-exact.

We can assume that we start with a perverse sheaf $\F$ supported on
$\Bun_n^{d'}$ with $d'\geq d+d_0$, and we have to show that
$\on{Av}^d_E(\F)\in \on{P}(\Bun_n^{d'-d})$ has no cohomologies in degrees $> 0$.

Suppose not, and let
$$\on{Av}^d_E(\F)\to \tau^{>0}(\on{Av}^d_E(\F))$$
be the truncation map. This map cannot be zero, unless
$\tau^{>0}(\on{Av}^d_E(\F))$ vanishes.

By \lemref{Av is cuspidal}, we know that $\on{Av}^d_E(\F)$ is cuspidal.
On the other hand, by \corref{iteration is exact}, we know that
$\tau^{>0}(\on{Av}^d_E(\F))$ projects to zero in $\Dt(\Bun_n)$. This
is a contradiction in view of Property 2 of $\Dt(\Bun_n)$.

\section{The symmetric group argument}

The goal of this section is to prove \thmref{exactness of Av^1},
assuming the existence of the quotient categories $\Dt(S\times \Bun_n)$
which satisfy Properties 0 and 1 of \secref{properties}.

\ssec{} Since the situation is essentially Verdier self-dual, it would 
be sufficient to prove that the functor
$\wt{\on{Av}}_E^1:\Dt(\Bun_n)\to \Dt(\Bun_n)$ is right-exact.

Let us suppose that it is not, and we will arrive to a contradiction.

\medskip

By definition, the functor $\wt{\on{Av}}_E^1$ is a composition of
an exact functor $$\wt{\on{H}}:\Dt(\Bun_n)\to \Dt(X\times \Bun_n)$$
followed by the functor 
$$\S\mapsto p_!(q^*(E)\otimes \S):
\Dt(X\times \Bun_n)\to \Dt(\Bun_n)$$
of cohomological amplitude $[-1,1]$. Thus, $\S\mapsto h^1(\wt{\on{Av}}_E^1(\S))$ 
is a right-exact functor $\Pt(\Bun_n)\to \Pt(\Bun_n)$.

Similarly, the amplitude of 
$\wt{\on{ItAv}}_E^i:\Dt(\Bun_n)\to \Dt(\Bun_n)$ is at most
$[-i,i]$, and $\S\mapsto h^i(\wt{\on{ItAv}}_E^i(\S))$ 
is a right-exact functor $\Pt(\Bun_n)\to \Pt^{\Sigma_i}(\Bun_n)$.

\begin{propconstr}   \label{canonical quotient}
For $\S\in \Pt(\Bun_n)$, 
there is a natural map in \newline
$\Pt^{\Sigma_i}(X^i\times \Bun_n)$:
$$\wt{\on{H}}^{\boxtimes i}(\S)\to (E^*)^{\boxtimes i}[i]\boxtimes 
h^i(\wt{\on{ItAv}}_E^i(\S)).$$
When $E$ is irreducible, the above map is surjective.
\end{propconstr}

\begin{proof}

The adjointness of the functors 
$p^!,p_!:\on{D}(X^i\times \Bun_n)\rightleftarrows \on{D}(\Bun_n)$
gives rise to a pair of mutually adjoint functors
$\on{P}(X^i\times \Bun_n)\rightleftarrows \on{P}(\Bun_n)$
given by
$$\F'\mapsto h^i(p_!(\F')) \text{ and } \F\mapsto p^*(\F)[i],$$
(the former being the left adjoint of the latter).
Since $X^i$ is smooth and connected, for $\F\in \on{P}(\Bun_n)$,
every sub-quotient of $p^*(\F)[-i]$ is of the form 
$p^*({}'\F)[-i]$, where $'\F$ is a sub-quotient of $\F$. This
implies that for any $\F'\in \on{P}(X^i\times \Bun_n)$, the adjunction
morphism $\F'\to p^*\left(h^i(p_!(\F'))\right)[i]$ is surjective.

\medskip

We have another pair of mutually adjoint functors between
the same categories:
$$\F'\mapsto h^i\left(p_!\left(q^*(E^{\boxtimes i})\otimes \F'\right)\right)
\text{ and } \F\mapsto (E^*)^{\boxtimes i}[i]\boxtimes \F,$$
and, when $E$ is irreducible, the adjunction map
$$\F'\mapsto (E^*)^{\boxtimes i}[i]\boxtimes 
h^i\left(p_!\left(q^*(E^{\boxtimes i})\otimes \F'\right)\right)$$
is also surjective as follows from the next lemma:

\begin{lem}
If for $\F'\in \on{P}(X^i\times \Bun_n)$ and $\F''\in \on{P}(\Bun_n)$ there
is a surjective map $q^*(E^{\boxtimes i})\otimes(\F')\to p^*(\F'')[i]$,
and $E$ is irreducible, then the adjoint map
$$\F'\to (E^*)^{\boxtimes i}[i]\boxtimes \F''$$
is also surjective.
\end{lem}

Moreover, the same assertions remain true for the corresponding functors
that act on the level of equivariant categories:
$\on{P}^{\Sigma_i}(X^i\times \Bun_n)\rightleftarrows \on{P}^{\Sigma_i}(\Bun_n)$.

\medskip

By passing to the quotient $\Dt^{\Sigma_i}(X^i\times \Bun_n)$, and using
Property 0 of the quotient categories, every $\S'\in \Dt^{\Sigma_i}(X^i\times \Bun_n)$
we obtain a functorial map
$$\S'\to (E^*)^{\boxtimes i}[i]\boxtimes 
h^i\left(p_!\left(q^*(E^{\boxtimes i})\otimes \S'\right))\right),$$
which is surjective if $E$ is irreducible.

By setting now $\S'=\wt{\on{H}}^{\boxtimes i}(\S)$ we arrive to the 
assertion of the proposition.

\end{proof}

\ssec{The case $i=n+1$}   \label{i=n+1}

Note that for $\S\in \Pt(\Bun_n)$ we have 
$$h^i(\wt{\on{ItAv}}_E^i(\S))\simeq 
h^1(\wt{\on{Av}}_E^1)\circ...\circ h^1(\wt{\on{Av}}_E^1)(\S),$$
as functors $\Pt(\Bun_n)\to \Pt(\Bun_n)$.

Therefore,
if for some $\S\in \Pt(\Bun_n)$, $h^i(\wt{\on{ItAv}}_E^i(\S))\neq 0$,
then $h^j(\wt{\on{ItAv}}_E^j(\S))\neq 0$ for all $j\leq i$.

Our first step will be to show that for all $\S\in \Pt(\Bun_n)$, 
$h^i(\wt{\on{ItAv}}_E^i(\S))=0$ for $i=n+1$, which would imply that the
same is true for all $i\geq n+1$.

\medskip

Consider the restriction of the surjection of \propconstrref{canonical quotient}
to the diagonal
$X\times \Bun_n\overset{\Delta\times \on{id}}
\longrightarrow X^i\times \Bun_n$, i.e., we have the following map
in $\Dt^{\Sigma_i}(X\times \Bun_n)$:
\begin{equation} \label{map on the diagonal}
(\Delta\times \on{id})^*\left(\wt{\on{H}}^{\boxtimes i}(\S)\right)[1-i]
\to 
(\Delta\times\on{id})^*
\Bigl((E^*)^{\boxtimes i}[i]\boxtimes h^i\left(\wt{\on{ItAv}}_E^i(\S)\right)\Bigr)[1-i].
\end{equation}

A key technical result, that we will need, states that both sides
of \eqref{map on the diagonal} belong in fact to $\Pt^{\Sigma_i}(X\times \Bun_n)$
and that the map in \eqref{map on the diagonal} is still surjective. 
In fact, we will prove the following:

\begin{prop}  \label{surjectivity on the diagonal}
Let $\K\in \on{P}(X^i\times\Bun_n)$ be a perverse sheaf, which appears
as a sub-quotient of some $h^k(\on{H}^{\boxtimes i}(\S))$ for some
object $\S\in \on{D}(\Bun_n)$. Then for any smooth sub-variety $X'\subset X^i$,
the *-restriction $\K|_{X'\times \Bun_n}$ lives in the cohomological
dimension $-\on{codim}(X',X^i)$.
\end{prop}

This proposition will be proved in \secref{proof of surj on diag}.
Let us now explain how it
implies what we need about \eqref{map on the diagonal}.

Indeed, $\wt{\on{H}}^{\boxtimes i}(\S)$ can be represented
by a sub-quotient $\K$ of $h^0(\on{H}^{\boxtimes i}(\F))$ for some
$\F\in \on{P}(\Bun_n)$. Hence, the left-hand side of \eqref{map on the diagonal}
can be represented by $(\Delta\times \on{id})^*(\K)[1-i]$, which belongs
to $\on{P}(X\times \Bun_n)$ according to \propref{surjectivity on the diagonal}.

The fact that the right-hand side of \eqref{map on the diagonal} belongs to
$\Pt(X\times \Bun_n)$ is obvious. In fact, it is isomorphic to
$(E^*)^{\otimes i}[1]\boxtimes h^i(\wt{\on{ItAv}}_E^i(\S))$.

Finally, the map of \propconstrref{canonical quotient} can be represented by a 
surjective map of perverse sheaves
$\K\to \K'$, where $\K$ is as above. 
By the long exact sequence, the cokernel of \eqref{map on the diagonal}
injects into 
$$h^1\Bigl((\Delta\times\on{id})^*(ker(\K\to \K'))[1-i]\Bigr),$$ 
which vanishes according to \propref{surjectivity on the diagonal}.

\medskip

Now we are ready to prove that $h^i(\wt{\on{ItAv}}_E^i(\S))=0$ for $i=n+1$.

Since the functor of taking $\Sigma_{n+1}$-invariants is exact,
the map in \eqref{map on the diagonal} will continue to be surjective
when we pass to the $sign$-isotypic components on both sides, i.e., we have:
\begin{align*}
&\Hom_{\Sigma_{n+1}}
\biggl(sign,
(\Delta\times \on{id})^*\left(\wt{\on{H}}^{\boxtimes n+1}(\S)\right)[-n]\biggr)
\twoheadrightarrow \\
&\Hom_{\Sigma_{n+1}}\biggl(sign,(E^*)^{\otimes n+1}[1]
\boxtimes h^{n+1}(\wt{\on{ItAv}}_E^{n+1}(\S))\biggr).
\end{align*}

Now, by \eqref{hecke on the quotient},
the left-hand side in the above formula is zero. By surjectivity,
the right-hand side must also be zero. But we claim that this can only happen if 
$h^{n+1}(\wt{\on{ItAv}}_E^{n+1}(\S))=0$.

Indeed, let $\rho$ be an irreducible $\Sigma_{n+1}$-representation,
which has a non-trivial isotypic component in $h^{n+1}(\wt{\on{ItAv}}_E^{n+1}(\S))$.
However, since $\on{rk}(E)\geq n+1$, by the 
Schur-Weyl theory, $\rho^*\otimes sign$ appears with a non-zero
multiplicity in $(E^*)^{\otimes n+1}$. Hence, $sign$ appears with a non-zero
multiplicity in $(E^*)^{\otimes n+1}[1]\boxtimes h^{n+1}(\wt{\on{ItAv}}_E^{n+1}(\S))$.

\ssec{Proof of \propref{surjectivity on the diagonal}}   \label{proof of surj on diag}

Recall the notion of universal local acyclicity
in the situation when we have an object $\F\in \on{D}(Z)$ on a scheme
(or stack) $Z$ over a smooth base $\Y$ (cf. \cite{SGA41/2} or \cite{BG}).
In our case $Z=X^i\times \Bun_n$, $\Y=X^i$.

The first observation is:

\begin{lem}
For any $\F\in \on{D}(\Bun_n)$, the object $\F'=\on{H}^{\boxtimes i}(\F)$
is ULA with respect to the projection $q:X^i\times \Bun_n\to X^i$.
\end{lem}

\begin{proof}

The lemma is proved by induction. Supposing the validity
for an integer $i=j$, let us deduce the corresponding 
assertion for $i=j+1$. In other words, it suffices to show that
if $\F''\in \on{D}(X^i\times \Bun_n)$ is ULA with respect to 
$X^i\times \Bun_n\to X^i$, then $\on{H}_{X^i}(\F'')\in 
\on{D}(X^{i+1}\times \Bun_n)$ is ULA with respect to 
$X^{i+1}\times \Bun_n\to X^{i+1}$.

Consider the diagram
$$X^i\times X\times \Bun_n 
\overset{\on{id}_{X^i}\times s\times \hl}\longleftarrow X^i\times \on{Mod}_n^1 
\overset{\on{id}_{X^i}\times s\times \hr}\longrightarrow X^i\times X \times \Bun_n.$$

By definition, 
$\on{H}_{X^i}(\F'')=(\on{id}_{X^i}\times s\times \hl)_!\circ 
(\on{id}_{X^i}\times \hr)^*(\F'')[i\cdot n]$. 

The ULA property is stable under direct images under proper morphisms.
Since the map
$\on{id}_{X^i}\times s\times \hl$ is proper, it is enough to show that
$(\on{id}_{X^i}\times \hr)^*(\F'')\in 
\on{D}(X^i\times \on{Mod}_n^1)$ is ULA with respect
to $X^i\times \on{Mod}_n^1\overset{\on{id}\times s}\to X^{i+1}$.
However, this follows from the assumption on $\F''$, since the map
$\on{id}_{X^i}\times s\times \hl: 
X^i\times \on{Mod}_n^1\to X^i\times X\times \Bun_n$
is smooth.

\end{proof}

The proposition will now follow from the next general observation:

Let $\F\in \on{D}(Z)$ be a complex, which is ULA with respect
to a projection $Z\to \Y$, where $\Y$ is smooth. Let $\K$ be a
sub-quotient of $h^k(\F)$ for some $k$, and let $\Y'\subset \Y$
be a smooth sub-variety. Denote by $Z'$ the preimage of $\Y'$ in $Z$.
In the above circumstances we have:

\begin{lem}
The *-restriction $\K|_{Z'}$ lives in the cohomological dimension
$-d$, where $d:=\on{codim}(\Y',\Y)$.
\end{lem}

\begin{proof}

Note that the assertion of the lemma implies that
$\K|_{Z'}[-d]$ is a sub-quotient of $h^{k-d}(\F|_{Z'})$.

Therefore, to prove the lemma, we can assume by induction
that $d=1$, and that, moreover, $\Y'$
is cut by the equation of a function with a non-vanishing differential.

Let $\Psi$, $\Phi$ be the corresponding near-by and vanishing cycles 
functors: $\on{D}(Z)\to \on{D}(Z')$.
By assumption, we have $\Phi(\F)=0$. The exactness of
$\Phi$ implies that $\Phi(\K)=0$ as well. Therefore,
$\K|_{Z'}\simeq \Psi(\K)[1]$, which is what we had to prove.

\end{proof}

\ssec{}

Let now $i$ be the maximal integer, for which the functor
$\S\mapsto h^i(\wt{\on{ItAv}}_E^i(\S)): \Dt(\Bun_n)\to \Dt(\Bun_n)$
is non-identically zero. We know already that $i\leq n$.
We are assuming that $i\geq 1$ and we want to arrive to a contradiction.

For $\S\in \Pt(\Bun_n)$, let us denote by $\S^i$ the object
$h^i(\wt{\on{ItAv}}_E^i(\S))\in \Pt^{\Sigma_i}(\Bun_n)$ and 
consider the canonical surjection of \propconstrref{canonical quotient}
$$\wt{\on{H}}^{\boxtimes i}(\S)\to (E^*)^{\boxtimes i}[i]\boxtimes \S^i.$$
Let us now apply the functor 
$\wt{\on{H}}_{X^i}^{\boxtimes n}: 
\Pt(X^i\times \Bun_n)\to \Pt(X^{i+n}\times \Bun_n)$ to both sides.
This functor maps $\Pt^{\Sigma_i}(X^i\times \Bun_n)$ to
$\Pt^{\Sigma_i\times \Sigma_n}(X^{i+n}\times \Bun_n)$

We obtain a morphism
\begin{equation} \label{map after another hecke}
\wt{\on{H}}^{\boxtimes i+n}(\S)\to (E^*)^{\boxtimes i}[i]\boxtimes 
\wt{\on{H}}^{\boxtimes n}(\S^i),
\end{equation}
which is still surjective, by the right exactness
of $\wt{\on{H}}_{X^i}^{\boxtimes n}$.

Note that the left-hand side of \eqref{map after another hecke}
is in fact an object of $\Pt^{\Sigma_{i+n}}(X^{i+n}\times \Bun_n)$.
We have a natural induction functor 
$$\on{Ind}_{\Sigma_i\times \Sigma_n}^{\Sigma_{i+n}}:
\on{P}^{\Sigma_i\times \Sigma_n}(X^{i+n}\times \Bun_n)\to 
\on{P}^{\Sigma_{i+n}}(X^{i+n}\times \Bun_n),$$ which is the left
(and right) adjoint to the forgetful functor.
By passing to the quotient we obtain the corresponding
induction functor from
$\Pt^{\Sigma_i\times \Sigma_n}(X^{i+n}\times \Bun_n)$ to 
$\Pt^{\Sigma_{i+n}}(X^{i+n}\times \Bun_n)$.

Thus, we obtain a map in $\Pt^{\Sigma_{i+n}}(X^{i+n}\times \Bun_n)$:
$$\wt{\on{H}}^{\boxtimes i+n}(\S)\to 
\on{Ind}_{\Sigma_i\times \Sigma_n}^{\Sigma_{i+n}}
\left((E^*)^{\boxtimes i}[i]\boxtimes \wt{\on{H}}^{\boxtimes n}(\S^i)\right).$$

The assumption that $i$ was maximal will yield the following:

\begin{prop}  \label{surjectivity of induction}
The above map
$$\wt{\on{H}}^{\boxtimes i+n}(\S)\to 
\on{Ind}_{\Sigma_i\times \Sigma_n}^{\Sigma_{i+n}}
\left((E^*)^{\boxtimes i}[i]\boxtimes \wt{\on{H}}^{\boxtimes n}(\S^i)\right)$$
is surjective.
\end{prop}

Let us conclude the proof of \thmref{exactness of Av^1} using this
proposition.

Let $\Delta_i:X\to X^i$, $\Delta_n:X\to X^n$, $\Delta_2:X\to X\times X$,
and $\Delta_{i+n}:X\to X^{i+n}$ be the corresponding diagonal
embeddings.
According to \propref{surjectivity on the diagonal}, 
in the formula
\begin{align*}
&(\Delta_{i+n}\times \on{id})^*
\left(\wt{\on{H}}^{\boxtimes i+n}(\S)\right)[1-i-n]\to  \\
&(\Delta_{i+n}\times \on{id})^*
\Bigl(\on{Ind}_{\Sigma_i\times \Sigma_n}^{\Sigma_{i+n}}
\left((E^*)^{\boxtimes i}[i]\boxtimes \wt{\on{H}}^{\boxtimes n}(\S^i)\right)\Bigr)
[1-i-n]
\end{align*}
both sides belong to $\Pt^{\Sigma_{i+n}}(X\times \Bun_n)$, and the map
is surjective. Therefore, the above map will still be surjective when we 
pass to the $sign$-isotypic
component on both sides with respect to $\Sigma_{i+n}$.

By \eqref{hecke on the quotient}, the left-hand side, i.e.,
$$\Hom_{\Sigma_{i+n}}\biggl(sign, 
(\Delta_{i+n}\times\on{id})^*
\left(\on{H}^{\boxtimes i+n}(\S)\right)[1-i-n]\biggr)$$ vanishes. Therefore,
so must the right-hand side.

Since the induction functors commute with the restriction functor
$(\Delta_{i+n}\times\on{id})^*$, by adjunction we obtain that
\begin{align*}
&\Hom_{\Sigma_{i+n}}\biggl(sign,
(\Delta_{i+n}\times\on{id})^*
\Bigl(\on{Ind}_{\Sigma_i\times \Sigma_n}^{\Sigma_{i+n}}
\left((E^*)^{\boxtimes i}[i]\boxtimes \wt{\on{H}}^{\boxtimes n}(\S^i)\right)\Bigr)
[1-i-n]\biggr)\simeq \\
&\Hom_{\Sigma_i\times \Sigma_n}
\biggl(\on{Res}^{\Sigma_{i+n}}_{\Sigma_i\times \Sigma_n}(sign),
(\Delta_{i+n}\times\on{id})^*\left((E^*)^{\boxtimes i}[i]\boxtimes 
\wt{\on{H}}^{\boxtimes n}(\S^i))[1-i-n]\right)\biggr).
\end{align*}

We have: 
$\on{Res}^{\Sigma_{i+n}}_{\Sigma_i\times \Sigma_n}(sign)\simeq
sign\times sign$, and 
$\Delta_{i+n}=\Delta_2\circ (\Delta_i\times\Delta_n)$. Let us,
therefore, rewrite the last expression as

\begin{equation} \label{expression for sign}
(\Delta_2\times \on{id})^*
\Biggl(\Hom_{\Sigma_i}\biggl(sign,(E^*)^{\otimes i}\boxtimes
\Hom_{\Sigma_n}
\left(sign,(\Delta_n\times\on{id})^*(\wt{\on{H}}^{\boxtimes n}(\S^i))[1-n]\right)
\biggr)\Biggr).
\end{equation}

Recall the multiplication map $m:X\times \Bun_n\to \Bun_n$, and recall
also from \eqref{hecke on the quotient}, that for $\S\in \Dt(\Bun_n)$
$$\Hom_{\Sigma_n}\left(sign,
(\Delta_n\times\on{id})^*(\wt{\on{H}}^{\boxtimes n}(\S))[1-n]\right)\simeq
m^*(\S)[1].$$

Therefore, \eqref{expression for sign} can be rewritten as
\begin{align*}
&(\Delta_2\times \on{id})^*
\biggl( \Hom_{\Sigma_i}\left(sign, (E^*)^{\otimes i}\boxtimes
m^*(\S^i)[1]\right)\biggr)\simeq \\
&\Hom_{\Sigma_i}\biggl(sign,q^*\left((E^*)^{\otimes i}\right)\otimes 
m^*(\S^i)[1]\biggr).
\end{align*}
As in \secref{i=n+1}, since $i\leq \on{rk}(E)$,
we obtain that the vanishing of the latter expression implies that
$m^*(\S^i)=0$. Therefore, the functor
$$\S\mapsto m^*\bigl(h^i\left(\wt{\on{ItAv}}_E^i(\S)\right)\bigr)$$ vanishes identically.

We claim that this implies that the functor
$\S\to h^i\left(\wt{\on{ItAv}}_E^i(\S)\right)$ vanishes. Indeed, for any fixed
$x\in X$, consider the pull-back map 
$m_x^*:\on{D}(\Bun_n)\to \on{D}(\Bun_n)$,
which is the composition of $m^*$ and the restriction
to $x\times \Bun_n\subset X\times \Bun_n$.

Obviously, $$m^*_x\circ h^i\left(\on{ItAv}_E^i(\S)\right)\simeq
h^i\left(\on{ItAv}_E^i(m^*_x(\S))\right).$$ Hence, the corresponding
functors on the level of $\Dt(\Bun_n)$ are also isomorphic.

Thus, we obtain that the functor
$\S\mapsto h^i\left(\wt{\on{ItAv}}_E^i(\S)\right)$ ``kills'' the image of
$m_x^*:\Pt(\Bun_n)\to \Pt(\Bun_n)$.

However, since $m_x^*:\on{D}(\Bun_n)\to \on{D}(\Bun_n)$ is
essentially surjective (i.e., surjective on objects), the same is true for
$m_x^*:\Pt(\Bun_n)\to \Pt(\Bun_n)$, in other words, 
$h^i\left(\wt{\on{ItAv}}_E^i(\S)\right)$ vanishes on the entire 
$\Pt(\Bun_n)$.

\ssec{Proof of \propref{surjectivity of induction}}

Observe that as an object of $\Pt(X^{i+n}\times \Bun_n)$,
$\on{Ind}_{\Sigma_i\times \Sigma_n}^{\Sigma_{i+n}}
((E^*)^{\boxtimes i}[i]\boxtimes \wt{\on{H}}^{\boxtimes n}(\S^i))$ can be
written as 
\begin{equation} \label{induction as a sum}
\underset{\sigma \in \Sigma_{i+n}}\oplus 
\sigma^*((E^*)^{\boxtimes i}[i]\boxtimes \wt{\on{H}}^{\boxtimes n}(\S^i)),
\end{equation}
where the sum is taken over the coset representatives of 
$\Sigma_{i+n}/\Sigma_i\times \Sigma_n$.

The proof of the proposition is based on the following observation:

\begin{lem}
Let $\K\to \underset{i}\oplus \K_i$ be a map of objects of an
Artinian abelian category, such that each of the maps
$\K\to \K_i$ is surjective. Assume that for $i\neq j$,
$\K_i$ and $\K_j$ have no isomorphic quotients. Then
the map $\K\to \underset{i}\oplus \K_i$ is surjective as well.
\end{lem}

We know that the map
$\wt{\on{H}}^{\boxtimes i+n}(\S)\to 
(E^*)^{\boxtimes i}[i]\boxtimes \wt{\on{H}}^{\boxtimes n}(\S^i)$ is surjective.
By the $\Sigma_{i+n}$-equivariance of $\wt{\on{H}}^{\boxtimes i+n}(\S)$, 
we obtain that each 
$\wt{\on{H}}^{\boxtimes i+n}(\S)\to 
\sigma^*\left((E^*)^{\boxtimes i}[i]\boxtimes \on{H}^{\boxtimes n}(\S^i)\right)$
is surjective as well.

To apply this lemma we need to verify that for
$\sigma_1,\sigma_2\in \Sigma_{i+n}$, which belong to different
cosets, the objects 
$\sigma_1^*\left((E^*)^{\boxtimes i}[i]\boxtimes \on{H}^{\boxtimes n}(\S^i)\right)$
and $\sigma_2^*\left((E^*)^{\boxtimes i}[i]\boxtimes \on{H}^{\boxtimes n}(\S^i)\right)$
of $\Pt(X^i\times X^n\times \Bun_n)$ have no isomorphic quotients.

Again, by $\Sigma_{i+n}$-equivariance, we can assume that $\sigma_1$ is 
the unit element, and $\sigma=\sigma_2$ is such that the permutation
that it defines on the set $\{1,...,i+n\}$ satisfies $\sigma(1)=i+1$.

For any $j\in \{1,...,k\}$, let $q_{k,j}$ denote the projection on the
$j$-th factor $X^k\times \Bun_n\to X$, and $p_{k,j}$ the complementary projection
on $X^{k-1}\times \Bun_n$.

Let $\S'$ be a quotient common to 
$(E^*)^{\boxtimes i}[i]\boxtimes\wt{\on{H}}^{\boxtimes n}(\S^i)$ and
$\sigma^*\left((E^*)^{\boxtimes i}[i]\boxtimes \wt{\on{H}}^{\boxtimes n}(\S^i)\right)$.
Since $E$ is irreducible, every (sub)-quotient of $(E^*)^{\boxtimes i}[i]
\boxtimes\wt{\on{H}}^{\boxtimes n}(\S^i)$ is of the form 
$(E^*)^{\boxtimes i}[i]\boxtimes \S''$, where $\S''\in \Pt(X^n\times \Bun_n)$ is a (sub)-quotient of 
$\wt{\on{H}}^{\boxtimes n}(\S^i)$. Therefore,

\begin{equation}
h^1\left(p_{i+n,1}{}_!(q_{i+n,1}^*(E)\otimes \S')\right)\neq 0.
\end{equation}

As in \propref{canonical quotient}, this implies:
$$h^1\Biggl(p_{i+n,1}{}_!\biggl(q_{i+n,1}^*(E)\otimes 
\sigma^*\left((E^*)^{\boxtimes i}[i]\boxtimes 
\wt{\on{H}}^{\boxtimes n}(\S^i)\right)\biggr)\Biggr)\neq 0,$$
which is equivalent to
$h^1\Biggl
(p_{i+n,i+1}{}_!\biggl(q_{i+n,i+1}^*(E)\otimes 
\left((E^*)^{\boxtimes i}[i]\boxtimes 
\wt{\on{H}}^{\boxtimes n}(\S^i)\right)\biggr)\Biggr)\neq 0$, and hence
$h^1\biggl(p_{n,1}{}_!\left(q_{n,1}^*(E)\otimes 
\wt{\on{H}}^{\boxtimes n}(\S^i)\right)\biggr)\neq 0$.

A simple diagram chase shows:
\begin{lem}
For any $\F\in \on{D}(\Bun_n)$,
$$p_{k,1}{}_!\left(q_{k,1}^*(E)\otimes \on{H}^{\boxtimes k}(\F)\right)\simeq
\on{H}^{\boxtimes k-1}(\on{Av}^1_E(\F)).$$
\end{lem}

The lemma implies that we have also an isomorphism
$$p_{k,1}{}_!\left(q_{k,1}^*(E)\otimes \wt{\on{H}}^{\boxtimes k}(\S)\right)\simeq
\wt{\on{H}}^{\boxtimes k-1}\left(\wt{\on{Av}}^1_E(\S)\right)$$
as functors $\Dt(\Bun_n)\to \Dt(X^{k-1}\times \Bun_n)$.

Therefore, the fact that
$h^1\biggl(p_{n,1}{}_!\left(q_{n,1}^*(E)\otimes \wt{\on{H}}^{\boxtimes n}(\S^i)
\right)\biggr)$
is non-zero as an  object of $\Pt(X^{n-1}\times \Bun_n)$
implies that 
$h^1\biggl(\wt{\on{H}}^{\boxtimes n-1}\left(\wt{\on{Av}}^1_E(\S^i)\right)\bigg)\neq 0$.
The exactness  of 
$\wt{\on{H}}^{\boxtimes n-1}:\Pt(\Bun_n)\to \Pt(X^{n-1}\times \Bun_n)$ forces
$h^1\left(\wt{\on{Av}}^1_E(\S^i)\right)
\simeq h^{i+1}(\wt{\on{ItAv}}_E^{i+1}(\S))\neq 0$.

However, this contradicts the assumption that $i$ was the
maximal integer for which
$h^i(\wt{\on{ItAv}}_E^i(\S))\neq 0$.
 
\section{Whittaker categories}     \label{Whittaker categories}

From this moment on we will be occupied with construction of the quotient categories
$\Dt(S\times \Bun_n)$. This will be done using the formalism of Whittaker
categories and functors between then. The first step, i.e., the definition of
the appropriate categories, is the goal of the present section, which we carry out
in a way similar to the definition of Whittaker categories in Sect. 6 of \cite{FGV}.

\ssec{Drinfeld's compactifications}  \label{compact}

Let $\Bun_n'$ be the stack classifying pairs $(\M,\kappa_1)$, where
$\kappa_1$ is a non-zero map $\Omega^{n-1}\to \M$. Let $\pi$ denote
the natural projection $\Bun_n'\to \Bun_n$.

Recall also the stack $\Qb$ introduced in \cite{FGV1}. We will now
introduce a series of stacks $\Qb_1,...\Qb_n$ with $\Qb_1=\Bun'_n$,
$\Qb_n=\Qb$, which interpolate between the two.

Namely, $\Qb_k$ classifies the data of a rank $n$ bundle $\M$
and a collection of non-zero maps
$$\kappa_i: \Omega^{n-1+...+n-i}\to \Lambda^i(\M),\,\,i=1,...,k,$$
which satisfy the Pl\"ucker relations in the sense of \cite{FGV1}.
  
Let, in addition, $\Qb_{k,ex}$ be the stack classifying the same data
$(\M,\kappa_1,...,\kappa_k)$ as above, but where we allow the last map,
i.e., $\kappa_k$, to vanish. In particular, $\Qb_k$ is an open
substack in $\Qb_{k,ex}$.

We have the natural forgetful maps $\pi_{k+1,k}:\Qb_{k+1}\to \Qb_k$,
and $\pi_{k+1,ex,k}:\Qb_{k+1,ex}\to \Qb_k$.

\medskip

We will introduce certain triangulated categories $\Dwk$ (resp., $\Dwkex$) 
of sheaves on $\Qb_k$ (resp., $\Dwkex$), that we will call the 
Whittaker categories.

Each $\Dwk$ will be a full triangulated subcategory of $\on{D}(\Qb_k)$ 
defined by the condition that
$\F\in \Dwk$ if its perverse cohomologies belong to a certain 
Serre subcategory $\Pwk\subset \on{P}(\Qb_k)$, singled out by 
some equivariance condition; and similarly for $\Dwkex$.

By definition, for $k=1$, $\Dwk$ is the entire $\on{D}(\Qb_1)=
\on{D}(\Bun_n')$, i.e., the equivariance condition in this case is
vacuous.

\ssec{}

For a fixed point $y\in X$, let $\Qb^y_k$
be an open substack that corresponds to the condition that neither
of the maps $\kappa_1,...,\kappa_k$ has a zero at $y$.

If $(\M,\kappa_1,...,\kappa_k)$ is a point of
$\Qb^y_k$, on the formal disk $\D_y$ around $y$ we obtain
a flag
$$0=\M_0\subset \M_1\subset...\subset \M_k\subset \M|_{\D_y}$$
with $\M_j/\M_{j-1}\simeq \Omega^{n-j}|_{\D_y}$. 

Let $N_{k,\D_y}$ be the group-scheme (of infinite type) over
$\Qb^y_k$ defined as follows: its fiber over a point
$(\M,\kappa_1,...,\kappa_k)$ as above consists of all
automorphisms of $\M|_{\D_y}$, which are strictly upper-triangular
with respect to the flag of the $\M_i$'s.

In addition, we have a group-indscheme $N_{k,\D^*_y}$,
which contains $N_{k,\D_y}$ as a group-subscheme, and whose fiber
over $(\M,\kappa_1,...,\kappa_k)$ consists of all automorphisms, strictly
upper-triangular with respect to the given flag, of 
$\M$ over the formal punctured disk $\D^*_y$.

As in \cite{FGV}, one can show that $N_{k,\D^*_y}$ is in fact
an ind-groupscheme. More precisely, $N_{k,\D^*_y}$ can be represented
as a union of group-schemes $N^i_{k,\D^*_y}, \, i\in\NN$ with 
$N^i_{k,\D^*_y}\supset N_{k,\D_y}$ and $N^i_{k,\D^*_y}/N_{k,\D_y}$
finite-dimensional.

The quotient $\H_{N^y_k}:=N_{k,\D^*_y}/N_{k,\D_y}$ is an
ind-scheme over the stack
$\Qb^y_k$ and is a version of the Hecke stack for the corresponding 
unipotent group. We have: $\H_{N^y_k}=\underset{i\in \NN}\cup \H_{N^y_k}^i$,
where $\H_{N^y_k}^i:=N^i_{k,\D^*_y}/N_{k,\D_y}$; the latter is
isomorphic to a tower of fibrations into affine spaces over $\Qb^y_k$.

We let $\on{pr}_k$ (resp., $\on{pr}^i_k$) denote the natural
projection $\H_{N^y_k}\to \Qb^y_k$ (resp., $\H_{N^y_k}^i\to \Qb^y_k$).

\ssec{The groupoids}

We claim that $\H_{N^y_k}$ carries a natural structure of groupoid over
$\Qb^y_k$. This is the standard procedure that makes the Hecke stack
a groupoid over the moduli space of bundles, cf. \cite{FGV}.

Namely, we define the second projection $\on{act}_k:\H_{N^y_k}\to \Qb^y_k$
as follows:

Given a point $(\M,\kappa_1,...,\kappa_k)\in \Qb^y_k$ and
an automorphism $g:\M|_{\D^*_y}\to \M|_{\D^*_y}$, the new bundle
$\M'$ is defined to be equal to $\M$ on $X-y$, and a meromorphic 
section $m'\in \Gamma(X-y,\M')$ is regular if $g(m')$, viewed as
an element of $\Gamma(\D^*_y,\M)$, belongs to $\Gamma(\D_y,\M)$.

The condition that $g$ is strictly upper-triangular means that
$\M'|_{\D_y}$ is still endowed with a filtration
$$0=\M'_0\subset \M'_1\subset...\subset \M'_k\subset \M'_n:=\M'|_{\D_y}$$
with $\M'_j/\M'_{j-1}\simeq \Omega^{n-j}|_{\D_y}$. 

Again, from the construction, the ``old'' maps 
$\kappa_i:\Omega^{n-1+...+n-i}\to \Lambda^i(\M)$, which are 
a priori meromorphic as maps 
$\Omega^{n-1+...+n-i}\to \Lambda^i(\M')$, are in fact regular, and thus define
the data $\kappa_i'$ for $\M'$.

\medskip

Let $\Qb^y_{k+1,ex}:=\Qb^y_k\underset{\Qb_k}\times \Qb_{k+1,ex}$
be the preimage of $\Qb^y_k$ in $\Qb_{k+1,ex}$. (Note that
$\Qb^y_{k+1}$ denotes a completely different stack; we have inclusions
$\Qb^y_{k+1}\subset \Qb_{k+1}\subset  \Qb_{k+1,ex} \supset \Qb^y_{k+1,ex}$.)

Consider the pull-back $\H_{N^y_k}\underset{\Qb^y_k}\times \Qb^y_{k+1,ex}$
as an ind-scheme over $\Qb^y_{k+1,ex}$. The next assertion follows
from the construction:

\begin{lem}   \label{groupoid lifts}
The fiber product $\H_{N^y_k}\underset{\Qb^y_k}\times \Qb^y_{k+1,ex}$
has a natural structure of groupoid over $\Qb^y_{k+1,ex}$, i.e.,
here exists a naturally defined map 
$\on{act}_{k,ex}:\H_{N^y_k}\underset{\Qb^y_k}\times \Qb^y_{k+1,ex}\to 
\Qb^y_{k+1,ex}$, which makes the following diagram Cartesian
$$
\CD
\Qb^y_{k+1,ex}  @<{\on{act}_{k,ex}}<< 
\H_{N^y_k}\underset{\Qb^y_k}\times \Qb^y_{k+1,ex}  \\
@V{\pi_{k+1,ex,k}}VV               @V{\on{id}\times \pi_{k+1,ex,k}}VV   \\
\Qb^y_k  @<{\on{act}_k}<< \H_{N^y_k}.
\endCD
$$
\end{lem}

We will denote by $\on{act}^i_{k,ex}$ the restriction of
$\on{act}_{k,ex}$ to the sub-groupoid
$\H^i_{N^y_k}\underset{\Qb^y_k}\times \Qb^y_{k+1,ex}$, and by
$\on{pr}_{k,ex}$ (resp., $\on{pr}^i_{k,ex}$)
the natural projection from
$\H_{N^y_k}\underset{\Qb^y_k}\times \Qb^y_{k+1,ex}$ (resp., from
$\H^i_{N^y_k}\underset{\Qb^y_k}\times \Qb^y_{k+1,ex}$) to 
$\Qb^y_{k+1,ex}$.

\ssec{The characters}

One more observation we need to make before introducing
the categories of interest is the following:

We claim that there exists a natural morphism
$$\chi_k:\H_{N^y_k}\underset{\Qb^y_k}\times \Qb^y_{k+1,ex}\to \AA^1.$$

Indeed, a point of $\H_{N^y_k}\underset{\Qb^y_k}\times \Qb^y_{k+1,ex}$ 
is the data of
$$(\M,\kappa_1,...,\kappa_k)\in \Qb^y_k,\,\,
\kappa_{k+1}:\Omega^{n-1+...+n-k-1}\to \Lambda^{k+1}(\M),\,\, 
g\in \on{Aut}(\M|_{\D^*_y}).$$
The endomorphism $(g-\on{Id})$ defines for every $i=1,...,k$ a map
$(\M/\M_i)|_{\D^*_y}\to (\M_i)_{\D^*_y}$, which we compose with
$$\Omega^{n-i-1}|_{\D^*_y}\to 
(\M/\M_i)|_{\D^*_y} \text{ and } (\M_i)_{\D^*_y}\to (\M_i/\M_{i-1})_{\D^*_y}\simeq
\Omega^{n-i}|_{\D^*_y}.$$
As a result, for every $i=1,...,k$ we obtain a map
$\Omega^{n-i-1}|_{\D^*_y}\to \Omega^{n-i}|_{\D^*_y}$,
well-defined up to a map regular on $\D_y$, due to the
corresponding ambiguity in $g$. By taking residues at $y$ 
we obtain $k$ points of $\AA^1$, i.e., we obtain well-defined
maps $^i\chi_k: \H_{N^y_k}\underset{\Qb^y_k}\times \Qb^y_{k+1,ex}\to\AA^1$
for $i=1,...,k$.

The map $\chi_k$ is defined as a composition
$$\H_{N^y_k}\underset{\Qb^y_k}\times \Qb^y_{k+1,ex}\overset{^i\chi_k}\to
\underset{i}\Pi (\AA^1)\overset{\on{sum}}\to \AA^1.$$
We will denote by $\chi^i_k$ the restriction of $\chi_k$ to
$\H^i_{N^y_k}\underset{\Qb^y_k}\times \Qb^y_{k+1,ex}\subset
\H_{N^y_k}\underset{\Qb^y_k}\times \Qb^y_{k+1,ex}$.

In what follows $\on{A-Sch}$ will denote the Artin-Schreier sheaf on $\AA^1$.

\ssec{}

Everything said above can be generalized in a straightforward
way when one point $y$ is replaced by a finite collection
of pairwise distinct points $\yb:=y_1,...,y_m$.

Namely, we have the open substack 
$$\Qb^{\yb}_k:=\underset{j}\cap\, \Qb^{y_j}_k\subset \Qb_k,$$
and the groupoid $\H_{N^{\yb}_k}$ over it.
In fact,
$$\H_{N^{\yb}_k}\simeq 
\H_{N^{y_1}_k}|_{\Qb^{\yb}_k}\underset{\Qb^{\yb}_k}\times...
\underset{\Qb^{\yb}_k}\times \H_{N^{y_m}_k}|_{\Qb^{\yb}_k}.$$

In other words, the groupoids $\H_{N^{y_j}_k}|_{\Qb^{\yb}_k}$,
$j=1,...,m$ acting on $\Qb^{\yb}_k$ pairwise commute in the 
natural sense, hence we can form the product groupoid, which can
be identified with $\H_{N^{\yb}_k}$.

\ssec{The categories on $\Qb^y_{k+1,ex}$}   \label{def of cat with y}

We define the category $\on{P}^W(\Qb^y_{k+1,ex})
\subset \on{P}(\Qb^y_{k+1,ex})$
to consist of all perverse sheaves $\F\in \on{P}(\Qb^y_{k+1,ex})$,
for which the following holds:

For each $i\in \NN$, there exists an isomorphism between the following
two sheaves on $\H^i_{N_k}\underset{\Qb^y_k}\times \Qb^y_{k+1,ex}$:
\begin{equation}  \label{N_k equivariance}
\on{act}^i_{k,ex}{}^*(\F)
\text{ and }
\on{pr}^i_{k,ex}{}^*(\F)\otimes \chi^i_k{}^*(\on{A-Sch})
\end{equation}
such that the restriction of this isomorphism to the unit section
$\Qb^y_{k+1,ex}\subset \H^i_{N_k}\underset{\Qb^y_k}\times \Qb^y_{k+1,ex}$
is the identity map.

Note that both sides of \eqref{N_k equivariance}
are objects of $\on{D}(\H^i_{N_k}\underset{\Qb^y_k}\times \Qb^y_{k+1,ex})$,
which become perverse after the cohomological shift by
$\on{dim.rel.}(\H^i_{N_k},\Qb^y_k)$, since both maps
$\on{pr}^i_{k,ex}$ and $\on{act}^i_{k,ex}$ are smooth of that
relative dimension.

Since $\on{A-Sch}$ is a $1$-dimensional lisse sheaf
and the fibers of $\on{pr}^i_{k,ex}$ are connected, if an isomorphism
of \eqref{N_k equivariance} exists, it is unique. Moreover
a family of such isomorphisms for $i\in\NN$ is necessarily compatible.
All this follows from the following general lemma:

\begin{lem}  \label{general equivariance}
Let $p:\Y_1\to \Y_2$ be a smooth surjective map between schemes (or stacks)
of relative dimension $d$, which has connected fibers. Then

\noindent{\em (1)}
$\F\mapsto p^*(\F)[d]$ is a full embedding of 
$\on{P}(\Y_2)$ into $\on{P}(\Y_1)$; its image is stable under
sub-quotients.

\noindent{\em (2)}
If, moreover, the fibers of $p$ are contractible, then the same is true
when we replace $\on{P}(\Y_i)$ by $\on{D}(\Y_i)$, i.e.,
$\on{D}(\Y_2)$ is a full triangulated subcategory of $\on{D}(\Y_1)$.
In particular, $\on{P}(\Y_2)\subset \on{P}(\Y_1)$ is stable 
under extensions, and is therefore a Serre subcategory.
\end{lem}

\medskip

Since $\on{pr}^i_{k,ex}:\H^i_{N_k}\underset{\Qb^y_k}\times \Qb^y_{k+1,ex}\to
\Qb^y_{k+1,ex}$ is a tower of affine fibrations, from 
\lemref{general equivariance} above, we obtain that 
$\on{P}^W(\Qb^y_{k+1,ex})$
is indeed a Serre subcategory of $\on{P}(\Qb^y_{k+1,ex})$.

\medskip

We define $\on{D}^W(\Qb^y_{k+1,ex})$  
as the full triangulated subcategory of $\on{D}(\Qb^y_{k+1,ex})$, 
consisting of objects whose perverse cohomologies belong to
$\on{P}^W(\Qb^y_{k+1,ex})$.

\medskip

From \lemref{general equivariance} it follows that for any
$\F\in \on{D}^W(\Qb^y_{k+1,ex})$ there exists a unique
isomorphism
$$\on{act}^i_{k,ex}{}^*(\F)\simeq \on{pr}^i_{k,ex}{}^*(\F)\otimes 
\chi^i_k{}^*(\on{A-Sch}),$$ compatible
with the restrictions of both sides to the unit section.

\medskip

In the same way, for a collection of 
pairwise distinct points $\yb=y_1,...,y_m$, one defines the categories 
$\on{P}^W(\Qb^{\yb}_{k+1,ex})$ and $\on{D}^W(\Qb^{\yb}_{k+1,ex})$,
the former being a Serre subcategory of
$\on{P}(\Qb^{\yb}_{k+1,ex})$, and the latter a full triangulated
subcategory of $\on{D}(\Qb^{\yb}_{k+1,ex})$.

\ssec{}   \label{stratifications}

To proceed, we need to recall a natural stratification defined 
on the stacks $\Qb_k$.

For a string of non-negative integers $\db:=(d_1,...,d_k)$
let $X^{(\db)}$ be the corresponding partially symmetrized power
of the curve $X$:
$$X^{(\db)}=\underset{j=1,...,k}\Pi \, X^{(d_j)}.$$

Let $^{\db}\Qb_k$ be the stack
that classifies the data of $(\M,\kappa_1,...,\kappa_k,D_1,...,D_k)$,
where $\M$ is as before, $D_i\in X^{(d_j)}$, and each $\kappa_i$
is an injective bundle map
$$\Omega^{n-1+...+n-i}(D_i)\to \Lambda^i(\M),$$
such that $\{\kappa_1,...,\kappa_k\}$ satisfy the Pl\"ucker relations.

We have a natural map $^{\db}\Qb_k\to \Qb_k$. It was shown in
\cite{BG} that each $^{\db}\Qb_k$ becomes a locally closed substack of
$\Qb_k$, and that, moreover, these substacks for various $\db$ define a locally finite 
decomposition of $\Qb_k$ into locally closed pieces.

Observe that each $^{\db}\Qb_k$ can be alternatively viewed as
a stack classifying the data of a vector bundle $\M$ endowed
with a filtration
$$0=\M_0\subset \M_1\subset...\subset\M_k\subset \M,$$
and identifications
$\M_i/\M_{i-1}\simeq \Omega^{n-i}(D_i-D_{i-1})$ for
$(D_1,...,D_k) \in X^{(\db)}$.

\ssec{}

For a string of integers $\db$ as above, let
$^{\db}\Qb^y_k$ denote the intersection $^{\db}\Qb_k\cap \Qb^y_k$.
Let $^{\db}\Qb^y_{k+1,ex}\subset {}^{\db}\Qb_{k+1,ex}$ be the 
preimages of $^{\db}\Qb^y_k$ and $^{\db}\Qb_k$, respectively, in
$\Qb_{k+1,ex}$.

Note that $^{\db}\Qb^y_{k+1,ex}$ is the stack that classifies the data
of $(D_1,...,D_k)\in (X-y)^{(\db)}$, a vector bundle $\M$, a filtration
$$0=\M_0\subset \M_1\subset...\subset\M_k\subset \M$$
with $\M_i/\M_{i-1}\simeq \Omega^{n-i}(D_i-D_{i-1})$ and, finally, a map
$\wt{\kappa}_{k+1}:\Omega^{n-k-1}(-D_k)\to \M/\M_k$.

Note that each $^{\db}\Qb^y_{k+1,ex}$ is stable under the action
of the groupoid $\H^i_{N_k}\underset{\Qb^y_k}\times \Qb^y_{k+1,ex}$.

Therefore, in the same way as above, we can introduce the category
$\on{P}^W({}^{\db}\Qb^y_{k+1,ex})$ 
(resp., $\on{D}^W({}^{\db}\Qb^y_{k+1,ex})$), which is a
Serre subcategory (resp., a full triangulated subcategory) of
$\on{P}({}^{\db}\Qb^y_{k+1,ex})$ (resp., $\on{D}({}^{\db}\Qb^y_{k+1,ex})$).

\medskip

From \lemref{general equivariance} we obtain:

\begin{lem} \label{restrictions determine}

\noindent{\em (1)}
The $*$ and $!$ restrictions $\on{D}(\Qb^y_{k+1,ex})\to 
\on{D}({}^{\db}\Qb^y_{k+1,ex})$ map the category $\on{D}^W(\Qb^y_{k+1,ex})$
into $\on{D}^W({}^{\db}\Qb^y_{k+1,ex})$.

\noindent{\em (2)}
The $*$ and $!$ direct image functors map
$\on{D}^W({}^{\db}\Qb^y_{k+1,ex})$ to
$\on{D}^W(\Qb^y_{k+1,ex})\subset \on{D}(\Qb^y_{k+1,ex})$.

\noindent{\em (3)}
An object $\F\in \on{D}(\Qb^y_{k+1,ex})$ belongs to
$\on{D}^W(\Qb^y_{k+1,ex})$ if and only if its $*$-restrictions 
(or, equivalently, $!$-restrictions)
to $^{\db}\Qb^y_{k+1,ex}$ belong to
$\on{D}^W({}^{\db}\Qb^y_{k+1,ex})$ for all $\db$.

\end{lem}

Of course, the same assertion holds when we replace
one point $y$ by a finite collection of points $\yb$.

\begin{proof}

The fact that each $^{\db}\Qb^y_{k+1,ex}$ is stable under the action
of $\H^i_{N_k}\underset{\Qb^y_k}\times \Qb^y_{k+1,ex}$
means that we have a commutative diagram
$$
\CD
^{\db}\Qb^y_{k+1,ex} @<<< \H^i_{N_k}\underset{\Qb^y_k}\times{}^{\db}\Qb^y_{k+1,ex}
@>>> ^{\db}\Qb^y_{k+1,ex} \\
@VVV   @VVV  @VVV  \\
\Qb^y_{k+1,ex} @<{\on{act}^i_{k,ex}}<< \H^i_{N_k}\underset{\Qb^y_k}\times\Qb^y_{k+1,ex}
@>{\on{pr}^i_{k,ex}}>> \Qb^y_{k+1,ex},
\endCD
$$
in which both squares are Cartesian. As was remarked above, an object 
$\F\in \on{D}(\Qb^y_{k+1,ex})$ belongs to $\on{D}^W(\Qb^y_{k+1,ex})$
if and only if we have a compatible system of isomorphisms
$\on{act}^i_{k,ex}{}^*(\F)\simeq \on{pr}^i_{k,ex}{}^*(\F)\otimes \chi_k^i{}^*(\on{A-Sch})$,
and similarly for $^{\db}\Qb^y_{k+1,ex}$. This implies assertions (1) and (2) of the
lemma.

Assertion (3) follows from (1) and (2), since the decomposition of $\Qb^y_{k+1,ex}$
into the strata $^{\db}\Qb^y_{k+1,ex}$ is locally finite.

\end{proof}

\ssec{}    \label{description of cat on strata}
We will now analyze how objects of $\on{D}^W(\Qb^y_{k+1,ex})$ can look like
when restricted to $^{\db}\Qb^y_{k+1,ex}$.

For $\db$ as above, let $^{\db}\Qb_{k+1,ex}{}'\subset {}^{\db}\Qb_{k+1,ex}$
denote the closed substack that corresponds to the condition that
for $i=1,...,k$, each $D_i':=D_i-D_{i-1}$ is an effective divisor
(by definition, $D'_1=D_1$) and that moreover for $i=1,...,k-1$ each
$D'_{i+1}-D'_i$ is effective, and 
$\wt{\kappa}_{k+1}:\Omega^{n-k-1}(-D_k)\to \M/\M_k$ factors as
$\Omega^{n-k-1}(-D_k)\to \Omega^{n-k-1}(D'_k)
\overset{\wt{\kappa}_{k+1}'}\to \M/\M_k$.

Note that we have a natural map 
$^{\db}\chi_k':{}^{\db}\Qb_{k+1,ex}{}'\to \AA^1$,
defined in a way similar to how the map $\chi_k$ was defined. Namely,
we have to sum up the classes of the successive extensions
$$0\to \Omega^{n-i}(D'_i)\to \M_{i+1}/\M_{i-1}\to 
\Omega^{n-i-1}(D'_{i+1})\to 0$$
in $$Ext^1\left(\Omega^{n-i-1}(D'_{i+1}),\Omega^{n-i}(D'_i)\right)\simeq
H^1(X,\Omega(D'_i-D'_{i+1}))\to H^1(X,\Omega)\simeq \AA^1$$ for
$i=1,...,k-1$ and
the class of the induced extension of $\Omega^{n-k-1}(D'_k)$ by
$\Omega^{n-k}(D'_k)$ by means of $\wt{\kappa}_{k+1}'$.

\medskip

Let $^{\db}\P_k$ denote the stack classifying the data of 
$(D_1,...,D_k)\in X^{(\db)}$ such that each $D'_i:=D_i-D_{i-1}$ and $D'_{i+1}-D'_i$ is
effective, a vector bundle $\M'$ of rank $n-k$ with a map 
$\Omega^{n-k-1}(D'_k)\to \M'$.

Note that we have a natural projection 
$\phi_k:{}^{\db}\Qb_{k+1,ex}{}'\to {}^{\db}\P_k$, with 
$\M':=\M/\M_k$ in the above notation. The map
$\phi_k$ is smooth and has contractible fibers.

Let $^{\db}\P^y_k\subset {}^{\db}\P_k$ be the open substack that 
corresponds to the condition that the divisors $D_i$ avoid $y$.
Let $^{\db}\Qb^y_{k+1,ex}{}':={}^{\db}\Qb_{k+1,ex}{}'\cap \Qb^y_{k+1,ex}$.

\medskip

The following proposition is a version of Lemma 6.2.8 of \cite{FGV}.

\begin{prop}  \label{description of category}

\noindent{\em (1)}
Every object of $\on{D}^W({}^{\db}\Qb^y_{k+1,ex})$ is supported on
$\on{D}^W({}^{\db}\Qb^y_{k+1,ex}{}')$. 

\noindent{\em (2)}
The functor 
$\F\mapsto \phi_k{}^*(\F)\otimes \chi_k'{}^*(\on{A-Sch})$
defines an equivalence of categories 
$\on{D}({}^{\db}\P^y_k)\to \on{D}^W({}^{\db}\Qb^y_{k+1,ex})$.
\end{prop}

A similar assertion holds when a single point $y$ is replaced 
by a finite collection $\yb$.

\begin{cor}  \label{restrictions well-behaved}
Let $\yb=y_1,...,y_m$ be a collection of points with $y$
being one of them. Then the restriction functor
$\on{D}(\Qb^{y}_{k+1,ex})\to \on{D}(\Qb^{\yb}_{k+1,ex})$
maps the category $\on{D}^W(\Qb^{y}_{k+1,ex})$ to 
$\on{D}^W(\Qb^{\yb}_{k+1,ex})$.
\end{cor}

\begin{proof} 

Let $\F$ be an object of $\on{D}^W(\Qb^{y}_{k+1,ex})$.
According to \lemref{restrictions determine}, it suffices to check 
that the $*$-restrictions $\F|_{{}^{\db}\Qb^{\yb}_{k+1,ex}}$ all 
belong to $\on{D}^W({}^{\db}\Qb^{\yb}_{k+1,ex})$.

We know that $\F|_{{}^{\db}\Qb^y_{k+1,ex}}$ can be described 
as a pull-back from $^{\db}\P^y_k$ as in 
\propref{description of category}, tensored by the pull-back of 
the Artin-Schreier sheaf, and this description obviously survives
the further restriction to $^{\db}\Qb^{\yb}_{k+1,ex}$.

\end{proof}

\ssec{}  \label{def of whit cat}

Finally, we define the category $\on{P}^W(\Qb_{k+1,ex})$
to consist of all perverse sheaves  $\F\in \on{P}(\Qb_{k+1,ex})$
for which $\F|_{\Qb^{\yb}_{k+1,ex}}$ belongs to
$\on{P}^W(\Qb^{\yb}_{k+1,ex})$ for all finite collections $\yb$.

According to \lemref{general equivariance}, $\on{P}^W(\Qb_{k+1,ex})$
is a Serre subcategory of $\on{P}(\Qb_{k+1,ex})$. We set
$\on{D}^W(\Qb_{k+1,ex})$ to be the full triangulated subcategory
of $\on{D}(\Qb_{k+1,ex})$ generated by $\on{P}^W(\Qb_{k+1,ex})$.
In other words, $\F\in \on{D}^W(\Qb_{k+1,ex})$ if and only if
all of its perverse cohomologies belong to $\on{P}^W(\Qb_{k+1,ex})$.

According to \corref{restrictions well-behaved}, in order to check
that $\F\in \on{D}(\Qb_{k+1,ex})$ belongs to $\on{D}^W(\Qb_{k+1,ex})$,
it is sufficient to check that 
$\F|_{\Qb^y_{k+1,ex}}\in \on{D}^W(\Qb^y_{k+1,ex})$ for
all points $y\in X$, i.e., it is enough to consider $1$-element sets.

\medskip

Consider now $\Qb_{k+1}\subset \Qb_{k+1,ex}$. Since this open substack
is stable under the action of the groupoids used in the definition
of $\on{P}^W(\Qb_{k+1,ex})$, the categories
$\on{P}^W(\Qb_{k+1})\subset \on{P}(\Qb_{k+1})$ and
$\on{D}^W(\Qb_{k+1})\subset \on{D}(\Qb_{k+1})$ are well-defined.

By \lemref{general equivariance}, the direct and inverse image functors
$\on{D}(\Qb_{k+1})\rightleftarrows \on{D}(\Qb_{k+1,ex})$
map the subcategories $\on{P}^W(\Qb_{k+1})$ and $\on{D}^W(\Qb^y_{k+1,ex})$
to one-another.

\medskip

We emphasize that by definition, for $k=0$, 
$\on{D}^W(\Qb_{k+1})=\on{D}(\Qb_{k+1})$.

\ssec{}   \label{base S}

Let now $S$ be an arbitrary ``base'' scheme. All the constructions of this
section go through when we replace $\Qb_k$ by the product
$S\times \Qb_k$.

In other words, we have well-defined categories $\on{P}^W(S\times \Qb_{k+1,ex})$,
$\on{D}^W(S\times \Qb_{k+1,ex})$, $\on{P}^W(S\times \Qb_k)$, and
$\on{D}^W(S\times \Qb_k)$.

Moreover, for a morphism $S_1\to S_2$, the two pairs of direct and inverse
image functors $\on{D}(S_1\times \Qb_k)\rightleftarrows 
\on{D}(S_2\times \Qb_k)$ map the categories $\on{D}^W(S_1\times \Qb_k)$
and $\on{D}^W(S_2\times \Qb_k)$ to one-another; and similarly for
the ``$ex$''-version.

\section{Whittaker functors}    \label{Whittaker functs}

The goal of this section is to prove the following theorem:

\begin{thm}   \label{Whittaker functors}
For each $k=1,...,n-1$ we have an equivalence of categories
$W_{k,k+1,ex}:\on{D}^W(\Qb_k)\to \on{D}^W(\Qb_{k+1,ex})$,
which maps $\on{P}^W(\Qb_k)$ to $\on{P}^W(\Qb_{k+1,ex})$.
The quasi-inverse functor is given by $\F'\mapsto 
\pi_{k+1,ex,k}{}_!(\F')$, which in this case is isomorphic to 
$\pi_{k+1,ex,k}{}_*(\F')$.
\end{thm}

\ssec{}   \label{explicit on strata}

As a first step, we will describe the functor 
$W_{k,k+1,ex}$ on the strata $^{\db}\Qb_{k+1,ex}$ for
$\db=d_1,...,d_k$. Namely, let $\on{D}^W({}^{\db}\Qb_{k+1,ex})$
(resp., $\on{D}^W({}^{\db}\Qb_k)$) be the corresponding subcategory of
$\on{D}({}^{\db}\Qb_{k+1,ex})$ (resp., $\on{D}({}^{\db}\Qb_k)$), and let
us describe the functor
$$^{\db}W_{k,k+1,ex}:\on{D}^W({}^{\db}\Qb_k)\to \on{D}^W({}^{\db}\Qb_{k+1,ex}).$$

Recall the substack $^{\db}\Qb_{k+1,ex}{}'\subset {}^{\db}\Qb_{k+1,ex}$,
cf. \secref{description of cat on strata},
and define the corresponding substack $^{\db}\Qb_k{}'\subset {}^{\db}\Qb_k$ as 
the intersection $^{\db}\Qb_{k,ex}{}'\cap \Qb_{k}$. 
Let $\on{D}^W({}^{\db}\Qb_{k+1,ex}{}')$
(resp., $\on{D}^W({}^{\db}\Qb_k{}')$) be the corresponding subcategory of
$\on{D}({}^{\db}\Qb_{k+1,ex}{}')$ (resp., $\on{D}({}^{\db}\Qb_k{}')$).

As in \propref{description of category}(1), every object of 
$\on{D}({}^{\db}\Qb_{k+1,ex})$ is supported on the closed substack
$^{\db}\Qb_{k+1,ex}{}'$, and similarly for $^{\db}\Qb_k{}'$.
Therefore, the required functor $^{\db}W_{k,k+1,ex}$  amounts to
a functor $\on{D}^W({}^{\db}\Qb_k{}')\to \on{D}^W({}^{\db}\Qb_{k+1,ex}{}')$.

As in \propref{description of category}(2), the category
$\on{D}^W({}^{\db}\Qb_{k+1,ex}{}')$ 
is equivalent to the category $\on{D}({}^{\db}\P_k)$
by means of 
$$\F\mapsto \phi_k{}^*(\F)\otimes \chi_k'{}^*(\on{A-Sch}).$$

We will now give a similar explicit description of
$\on{D}^W({}^{\db}\Qb_k{}')$.

\medskip

For $\db=d_1,...,d_k$ let $^{\db}\P_{k-1}'$
denote the stack classifying the data of
$(D_1,...,D_k)\in X^{(\db)}$ with $D'_i=D_i-D_{i-1}$ effective for
$i=1,...,k$, and $D'_{i+1}-D'_i$ effective for $i=1,...,k-1$, a
vector bundle $\M''$ of rank $n-k+1$ with an {\it injective bundle map} 
$\wt{\kappa}_k:\Omega^{n-k}(D'_k)\to \M''$.

We have a natural projection
$\phi_{k-1}':{}^{\db}\Qb_k{}'\to {}^{\db}\P_{k-1}'$ that sends a point
$$\biggl((D_1,...,D_k) \in X^{(\db)};\,\,\,
0=\M_0\subset \M_1\subset...\subset\M_k\subset \M;\,\,\,
\M_i/\M_{i-1}\simeq \Omega^{n-i}(D'_i)\biggr)\in {}^{\db}\Qb_k{}'$$
to $\M'':=\M/\M_{k-1}$ and
$$\Omega^{n-k}(D'_k)\simeq \M_k/\M_{k-1}\hookrightarrow \M/\M_{k-1}=\M''.$$

Again, as in \propref{description of category}, the category
$\on{D}^W({}^{\db}\Qb_k{}')$ is equivalent to $\on{D}({}^{\db}\P_{k-1}')$,
by means of $\F\mapsto \phi_{k-1}'{}^*(\F)\otimes \chi_{k-1}'{}^*(\on{A-Sch})$.

\medskip

Note that the last two pieces of data in the definition of
$^{\db}\P_{k-1}'$, i.e., $(\M'',\wt{\kappa}_k)$,
can be rewritten as a short exact sequence
$$0\to \Omega^{n-k}(D'_k)\to \M''\to \M'\to 0,$$
where $\M'$ is a vector bundle of rank $n-k$. From this
it is easy to see that the stacks $^{\db}\P_{k-1}'$
and $^{\db}\P_k$ form a pair of mutually dual vector bundles
over the base classifying $(\M',D_1,...,D_k)$, which is isomorphic
to the product $\Bun_{n-k}\times \underset{i=1,...,k}\Pi\, X^{(d'_i-d'_{i-1})}$,
where $\db=d_1,...,d_k$ and $d'_i=d_i-d_{i-1}$.

We define the functor 
$^{\db}W_{k,k+1,ex}: \on{D}^W({}^{\db}\Qb_k{}')\to 
\on{D}^W({}^{\db}\Qb_{k+1,ex}{}')$
as a composition
$$\on{D}^W({}^{\db}\Qb_k{}')\simeq \on{D}({}^{\db}\P_{k-1}')\overset{\on{Four}'}
\longrightarrow \on{D}({}^{\db}\P_k)\simeq \on{D}^W({}^{\db}\Qb_{k+1,ex}{}'),$$
where $\on{Four}'$ is the Fourier transform functor
$\on{D}({}^{\db}\P_{k-1}')\to \on{D}({}^{\db}\P_k)$ followed
by the cohomological shift by 
$\on{dim.rel.}({}^{\db}\Qb_{k+1,ex}{}',{}^{\db}\P_k)-
\on{dim.rel.}({}^{\db}\Qb_k{}',{}^{\db}\P_{k-1}')$.

The functor $^{\db}W_{k,k+1,ex}$ is an equivalence of categories
mapping perverse sheaves to perverse sheaves, because the same it true for
the Fourier transform functor.

\medskip

Let $^{\db}\pi_{k+1,ex,k}:{}^{\db}\Qb_{k+1,ex}\to {}^{\db}\Qb_k$
be the restriction of $\pi_{k+1,ex,k}$ to the corresponding stratum.

\begin{lem}   \label{Whittaker functor on strata}
The functor $\F'\mapsto {}^{\db}\pi_{k+1,ex,k}{}_!(\F')$ maps
$\on{D}^W({}^{\db}\Qb_{k+1,ex})$ to $\on{D}^W({}^{\db}\Qb_k)$
and induces a functor quasi-inverse to $^{\db}W_{k,k+1,ex}$.
Moreover, in the above formula the $!$-direct image coincides
with the $*$ one.
\end{lem}

\begin{proof}

We have the following Cartesian diagram:
$$
\CD
^{\db}\Qb_k{}'  @<{^{\db}\pi_{k+1,ex,k}}<<  ^{\db}\Qb_{k+1,ex}{}'   \\
@V{\phi_{k-1}'}VV                    @V{\phi_k}VV   \\
^{\db}\P_{k-1}'  @<<<   
^{\db}\P_{k-1}'\underset{\Bun_{n-k}\times \underset{i=1,...,k}\Pi\, X^{(d'_i-d'_{i-1})}}\times
{}^{\db}\P_k.
\endCD
$$

Therefore, the assertion of the lemma can be translated to the following
general situation:

Let $\varphi:\E\to \Y$ be a vector bundle, and $\check \varphi:\check \E\to \Y$
its dual.
Consider the functor 
$W_\E:\on{D}(\E)\to \on{D}(\E\underset{\Y}\times \check\E)$
given by
$$\F\mapsto (\varphi\times \on{id})^*(\on{Four}(\F))\otimes 
\on{ev}^*(\on{A-Sch})[\on{dim.rel.}(\E,\Y)],$$
where $\varphi\times \on{id}$ is the natural projection
$\E\underset{\Y}\times \check\E\to \check \E$, and
$\on{ev}:\E\underset{\Y}\times \check\E\to\AA^1$ is the evaluation map.

Then for $\on{id}\times \check\varphi:\E\underset{\Y}\times \check\E\to \E$
we have:
$$(\on{id}\times \check\varphi)_!(W_\E(\F))\simeq
(\on{id}\times \check\varphi)_*(W_\E(\F))\simeq \F,$$
and this follows from the standard properties of the Fourier
transform functor.

\end{proof}

\ssec{}   \label{two groupoids}

We are now going to extend the above stratum-by-stratum definition of
$^{\db}W_{k,k+1,ex}$ to a globally defined functor
$W_{k,k+1,ex}:\on{D}^W(\Qb_k)\to \on{D}^W(\Qb_{k+1,ex})$.

We will first construct the functor
$W^y_{k,k+1,ex}:\on{D}^W(\Qb^y_k)\to \on{D}^W(\Qb^y_{k+1,ex})$. (The same
definition works for $y$ replaced by a finite collection of points
$\yb$.)

For that we will single out two sub-groupoids in the groupoid
$\H_{N^y_k}$ over $\Qb^y_k$, denoted $'\H_{N^y_k}$ and $''\H_{N^y_k}$,
respectively. Both these groupoids correspond to certain group
sub-schemes $'N_{k,\D_y},{}''N_{k,\D_y}$ of $N_{k,\D_y}$
(resp., $'N_{k,\D^*_y},{}''N_{k,\D^*_y}\subset N_{k,\D^*_y}$).

Recall that $N_{k,\D_y}$ (resp., $N_{k,\D^*_y}$) consists
of automorphisms of $\M|_{\D_y}$ (resp., $\M|_{\D^*_y}$), which
are strictly upper-triangular with respect to the filtration
$0=\M_0\subset \M_1\subset...\subset \M_k\subset \M|_{\D_y}$,
defined by our point of $\Qb^y_k$.

The group $'N_{k,\D_y}$ (resp., $'N_{k,\D^*_y}$)
consists of those automorphisms, which induce the identity
map on $\M_k$.
The group $''N_{k,\D_y}$ (resp., $''N_{k,\D^*_y}$)
consists of those automorphisms, 
which induce the identity map on $\M/\M_{k-1}$.

Note that the fiber of $'\H_{N^y_k}$ over a point
$(\M,\kappa_1,...,\kappa_k)\in \Qb^y_k$ is the vector space
$\Hom_{\D_y^*}(\M/\M_k,\M_k)/\Hom_{\D_y}(\M/\M_k,\M_k)$.
Without restricting the generality we can assume that the filtration
$\H_{N^y_k}=\underset{i\in \NN}\cup \H_{N^y_k}^i$ induces on $'N_{k,\D_y}$
the standard filtration: 
$$'\H_{N^y_k}^i=\Hom_{\D_y}\left(\M/\M_k,\M_k(i\cdot y)\right)/
\Hom_{\D_y}\left(\M/\M_k,\M_k\right).$$

Let $'\on{pr}_k$ and $'\on{act}_k$ (resp., $'\on{pr}^i_k$, $'\on{act}^i_k$)
be the restrictions to $'\H_{N^y_k}$ (resp., $'\H_{N^y_k}^i$) of the maps
$\on{pr}_k,\on{act}_k:\H_{N^y_k}\to \Qb^y_k$, respectively.
We will denote $'\H_{N^y_k}^i$ also by $\E_k^i$ and think of it as a
vector bundle over $\Qb^y_k$. Let $\check\E_k^i$ denote the dual
vector bundle, and ${'\check{\on{pr}}^i_k}$ its projection to $\Qb^y_k$.

By Serre's duality, the fiber of $\check\E_k^i$ over
$(\M,\kappa_1,...,\kappa_k)\in \Qb^y_k$ can be identified with the
vector space $$\Hom_{\D_y}
\biggl(\M_k,\left((\M/\M_k)/(\M/\M_k)(-iy)\right)\otimes\Omega\biggr).$$
For $i'\geq i$ we have a natural map $\on{pr}_{i',i}:\check\E_k^{i'}\to \check\E_k^i$.

\begin{prop}   \label{closed embeddings}
We have a natural map $f_i:\Qb^y_{k+1,ex}\to \check\E_k^i$ for any $i$.
Moreover,

\noindent{\em(1)}
For $i'\geq i$ the composition 
$\Qb^y_{k+1,ex}\overset{f_{i'}}\to 
\check\E_k^{i'}\overset{\on{pr}_{i',i}}\to \check\E_k^i$ equals $f_i$.

\noindent{\em(2)}
For each open substack $U\subset \Qb^y_k$ of finite type, there
exists an integer $i(U)$ large enough such that over $U$, the 
map $f_i:\Qb^y_{k+1,ex}\to \check\E_k^i$ is a closed embedding for every
$i\geq i(U)$.
\end{prop}

\begin{proof}

Let $\check\E_k^i{}'\subset \check\E_k^i$ be a vector sub-bundle, whose fiber
over a point $(\M,\kappa_1,...,\kappa_k)\in \Qb^y_k$ is the vector
space
\begin{equation} \label{fiber of vector bundle}
\Hom_{\D_y}\biggl(\Omega^{n-k},
\left((\M/\M_k)/(\M/\M_k)(-iy)\right)\otimes\Omega\biggr),
\end{equation}
which maps to 
$\Hom_{\D_y}
\biggl(\M_k,\left((\M/\M_k)/(\M/\M_k)(-iy)\right)\otimes\Omega\biggr)$
by means of the projection $\M_k\to \Omega^{n-k}|_{\D^y}$.

Note that given a filtration 
$$0=\M_0\subset \M_1\subset...\subset \M_k\subset \M$$
with $\M_j/\M_{j-1}\simeq \Omega^{n-j}$,
specifying a map $\Omega^{n-k}\to 
\left((\M/\M_k)/(\M/\M_k)(-iy)\right)\otimes\Omega$ is 
the same as to specify a map
$\Omega^{n-1+...+n-k-1}\to \Lambda^{k+1}(\M)/\Lambda^{k+1}(\M)(-iy)$,
which satisfies the Pl\"ucker relations with all the maps
$\kappa_j:\Omega^{n-1+...+n-j}\to \Lambda^j(\M)$, $j=1,...,k$.

\medskip

The map $\Qb^y_{k+1,ex}\to \check\E^i$ is defined now as follows:
having $(\M,\kappa_1,...,\kappa_k)\in \Qb^y_k$, to the data of
$\kappa_{k+1}:\Omega^{n-1+...+n-k-1}\to \Lambda^{k+1}(\M)$
we attach the corresponding map
$\Omega^{n-1+...+n-k-1}\to \Lambda^{k+1}(\M)/\Lambda^{k+1}(\M)(-iy)$
over $\D_y$. According to the above discussion, this defines a point
of $\check\E_k^i{}'$, and hence of $\check\E_k^i$.

\medskip

Thus, the map $f_i$ has been constructed.
Point (1) of the proposition is straightforward from the
construction.

For an open substack $U$ of finite type, let $i(U)$ 
be such that the vector space 
$\Hom(\Omega^{n-1+...+n-k-1},\Lambda^{k+1}(\M)(-iy))$ is zero
for $(\M,\kappa_1,...,\kappa_k)\in U$.

Let $\check\E_k^i{}''$ be the vector bundle over $\Qb^y_k$, whose fiber
over a point as above is the vector space
$\Hom\left(\Omega^{n-1+...+n-k-1},
\Lambda^{k+1}(\M)/\Lambda^{k+1}(\M)(-iy)\right)$. For $i\geq U(i)$
the natural map $\Qb^y_{k+1,ex}\to \check\E_k^i{}''$ is a closed embedding
over $U$.

Then for $i\geq i(U)$, we have a sequence of maps
$$\Qb^y_{k+1,ex}\to \check\E_k^i{}'\to \check\E_k^i{}''.$$
We know that the second arrow is a closed embedding, being
the set of those sections that satisfy the Pl\"ucker relations.
We also know that the composed map is a closed embedding.
Hence, so is the first map.

\end{proof}

\ssec{}

Consider now the action map
${}'\on{act}^i_k:{}'\H_{N^y_k}^i\to \Qb^y_k$. Since the projection
${}'\on{pr}^i_k:{}'\H_{N^y_k}^i\to \Qb^y_k$ is a smooth map
and $'\H_{N^y_k}^i$ is a groupoid, the map $'\on{act}^i_k$ is smooth 
as well; let $\dim(i,k)$ denote its relative dimension.

We define the functor 
$W^{y,i}_{k,k+1,ex}:\on{D}^W(\Qb^y_k)\to \on{D}(\check\E_k^i)$
as 
$$\F\mapsto \on{Four}\left({}'\on{act}^i_k{}^*(\F)[\dim(i,k)]\right),$$
where $\on{Four}:\on{D}(\E_k^i)\to \on{D}(\check\E_k^i)$ is
the Fourier transform functor. Evidently, this functor is exact.
For $i'\geq i$ we have:
$$\on{pr}_{i',i}{}_!(W^{y,i'}_{k,k+1,ex}(\F))\simeq W^{y,i}_{k,k+1,ex}(\F).$$

\begin{prop}   \label{support of sheaves}
For an open substack $U\subset \Qb^y_k$ of finite type and any 
integer $i$ which is large enough (and in particular $i\geq i(U)$ of 
\propref{closed embeddings}), over the preimage of $U$, any object 
of the form $W^{y,i}_{k,k+1,ex}(\F)$ for $\F\in \on{D}^W(\Qb^y_k)$
is supported on $\Qb^y_{k+1,ex}\subset \check\E_k^i$.
\end{prop}

\begin{proof}

Recall that to a string of non-negative integers
$\db=d_1,...,d_k$ we attached a locally closed substack
$^{\db}\Qb^y_k\subset \Qb^y_k$. Let $|\db|$ be $\underset{i}\Sum\, d_i$.

It is easy to see that for every open substack $U\subset \Qb_k$ of finite type
there exists an integer $d$ such that
\begin{equation} \label{stacks of finite type}
U\subset \underset{|\db|\leq d}\cup\, {}^{\db}\Qb_k.
\end{equation}

Thus, for a given $U$ there exists an integer $i'(U)$ such that
for any $\db$ and $\M\in U\cap {}^{\db}\Qb^y_k$
we have: $\Hom\left(\M_k,\M/\M_k\otimes\Omega(-iy)\right)=0$,
for $i\geq i'(U)$ where
$$0=\M_0\subset \M_1\subset...\subset \M_k\subset \M$$
is the filtration with $\M_i/\M_{i-1}\simeq \Omega^{n-i}(D'_i)$,
$(D_1,...,D_k)\in (X-y)^{(\db)}$.

\medskip

Let now $\F$ be an object of $\on{D}^W({}^{\db}\Qb^y_k)$ with
$|\db|\leq d$. Then $W^{y,i}_{k,k+1,ex}(\F)$ yields an object of
$\on{D}\left(({'\check{\on{pr}}^i_k})^{-1}({}^{\db}\Qb^y_k)\right)$,
where $({'\check{\on{pr}}^i_k})^{-1}({}^{\db}\Qb^y_k)$ is the preimage 
of $^{\db}\Qb^y_k$ in $\check\E_k^i$.

\begin{lem}  \label{compatible on strata}
For $i\geq i'(U),i(U)$, over the preimage of $U$,
the object $W^{y,i}_{k,k+1,ex}(\F)$ is supported
on $^{\db}\Qb^y_{k+1,ex}$ and is isomorphic to  
$^{\db}W_{k,k+1,ex}(\F)$ of \secref{explicit on strata}.
\end{lem}

\begin{proof}

Recall the stacks $^{\db}\P_k$ and $^{\db}\P'_{k-1}$, which form a pair
of mutually dual (generalized) vector bundles over the base
$\Bun_{n-k}\times \underset{i=1,...,k}\Pi\, X^{(d'_i-d'_{i-1})}$, where the
latter classifies the data of a rank $n-k$ vector bundle $\M'$ and a collection
of divisors $D_1,...,D_k$ with $D'_i=D_i-D_{i-1}$ effective and
$D'_i-D'_{i-1}$ effective as well. To simplify the notation, let
us temporarily denote this base by $\Y$, $^{\db}\P'_{k-1}$ by $\E$ and
$^{\db}\P'_{k-1}$ by $\check \E$; let $\varphi$ and $\check{\varphi}$
denote the projections of $\E$ and $\check \E$, respectively, on $\Y$.

Consider the fiber product $\E\underset{\Y}\times \check \E$,
and let $\on{ev}$ be the natural evaluation
map from it to ${\mathbb A}^1$. The assertion of the lemma amounts to a description
of the functor $\on{D}(\E)\to \on{D}(\E\underset{\Y}\times \check \E)$ given by
\begin{equation}  \label{usual Fourier}
\F\mapsto (\varphi\times \on{id})^*(\on{Four}(\F))
\otimes \on{ev}^*(\on{A-Sch})[\on{dim.rel.}(\E,\Y)]
\end{equation}
in terms of an action of a certain groupoid on $\E$.

Namlely, for an integer $i$ consider another vector bundle over $\Y$,
denoted $\E^i$, whose fiber over 
$(\M',D_1,...,D_k)\in \Y$ is the vector space 
$\on{Hom}(\M',\Omega^{n-k}(i\cdot y)/\Omega^{n-k})$. Let $\check \E^i$
be the dual vector bundle, whose fiber, by Serre's duality can
be identified with $\on{Hom}(\Omega^{n-k-1},\M'/\M'(-i\cdot y))$.
We have a natural map $\E^i\to \E$. When working 
over an open substack of finite type in $\Bun_{n-k}$,
for a large enough integer $i$, the dual map $\check \E\to \check \E^i$
is a closed embedding.

Using the group-scheme structure, we can think of $\E^i$ as 
a groupoid acting on $\E$. Let
$a$ and $p$ denote the corresponding maps
$\E_i\underset{\Y}\times \E\to \E$. Thus, we can consider the
functor $\on{D}(\E)\to \on{D}(\check \E_i\underset{\Y}\times \E)$
given by
\begin{equation}  \label{groupoid Fourier}
\F\mapsto \on{Four}\left(a^*(\F)\right)[\on{dim.rel.}(\E^i,\Y)],
\end{equation}
where $\on{Four}$ is the relative Fourier transform functor
$\on{D}(\E_i\underset{\Y}\times\E)\to \on{D}(\check \E_i\underset{\Y}\times \E)$.

The assertion of the lemma follows from the fact the functor in
\eqref{groupoid Fourier} is isomorphic to the composition of the
functor of \eqref{usual Fourier}, followed by the direct image 
under the closed embedding $\E\underset{\Y}\times\check \E\to 
\E\underset{\Y}\times\check \E^i$.

\end{proof}

This lemma implies the proposition. Indeed for a given
$\F\in \on{D}^W(\Qb^y_k)$ to show that, over the preimage of $U$, 
$W^{y,i}_{k,k+1,ex}(\F)$ is supported on $\Qb^y_{k+1,ex}$, 
it is enough to do so over the preimage of each stratum 
$^{\db}\Qb^y_k\subset \Qb^y_k$.
The latter support property is insured by \lemref{compatible on strata}.

\end{proof}

\ssec{}   \label{construction of whittaker functors}

Since $\on{pr}_{i',i}{}_!(W^{y,i'}_{k,k+1,ex}(\F))\simeq 
W^{y,i}_{k,k+1,ex}(\F)$, the above proposition implies that
we obtain a well-defined functor
$W^y_{k,k+1,ex}:\on{D}^W(\Qb^y_k)\to \on{D}(\Qb^y_{k+1,ex})$.
\footnote{This definition of $W^y_{k,k+1,ex}$ has been inspired
by a certain construction of V.Drinfeld in the $n=2$ case, 
one incarnation of which is explained in \secref{Drinfeld's transform}.}

Moreover, by combining \lemref{compatible on strata} and
\lemref{restrictions determine}(3) we obtain that the image of
$W^y_{k,k+1,ex}$ lies in $\on{D}^W(\Qb^y_{k+1,ex})$.

\begin{prop}   \label{functors with y}
The direct image functor $\F\mapsto \pi_{k+1,ex,k}{}_!(\F)$ maps
$\on{D}^W(\Qb^y_{k+1,ex})$ to $\on{D}^W(\Qb^y_k)$ and is a 
quasi-inverse to $W^y_{k,k+1,ex}$. Moreover, for
$\F\in \on{D}^W(\Qb^y_{k+1,ex})$, 
$\pi_{k+1,ex,k}{}_!(\F)\to \pi_{k+1,ex,k}{}_*(\F)$
is an isomorphism.
\end{prop}

\begin{proof}

First, let us show that for $\F\in \on{D}^W(\Qb^y_k)$,
$\pi_{k+1,ex,k}{}_!(W^y_{k,k+1,ex}(\F))\simeq \F$.

Indeed, by working over a fixed stack $U$ of finite type
and a large enough integer $i$, we are reduced to showing that
$${'\check{\on{pr}}^i_k}{}_!
\biggl(\on{Four}\left({}'\on{act}^i_k{}^*(\F)[\dim(i,k)]\right)\biggr)\simeq
\F,$$
where ${'\check{\on{pr}}^i_k}$ is the projection $\check \E^i_k\to \Qb^y_k$.

However, by the general properties of the Fourier transform functor
we obtain that the left-hand side of the above expression is isomorphic
to the restriction of $'\on{act}^i_k{}^*(\F)$ to the unit section of
$\E^i_k\simeq{}'\H^y_{N_k}$, i.e., to $\F$ itself.

\medskip

Now let us show that $\pi_{k+1,ex,k}{}_!$ 
indeed maps $\on{D}^W(\Qb^y_{k+1,ex})$
to $\on{D}^W(\Qb^y_k)$. However, this follows immediately from the 
definitions: 

By unfolding the definition of $\on{D}^W(\Qb^y_k)$, we see that 
it is defined by means of an equivariance property with respect to 
the groupoid $''\H_{N^y_k}$ (cf. \secref{two groupoids}).
However, $''\H_{N^y_k}$ acts on $\Qb^y_{k+1,ex}$, being a part
of $\H_{N^y_k}$, i.e., we have a Cartesian diagram
$$
\CD
\Qb^y_{k+1,ex}  @<{\on{act}_{k,ex}}<< ''\H_{N^y_k}\underset
{\Qb^y_k}\times \Qb^y_{k+1,ex}   \\
@V{\pi_{k+1,ex,k}}VV               @V{\on{id}\times \pi_{k+1,ex,k}}VV   \\
\Qb^y_k  @<{\on{act}_k}<< ''\H_{N^y_k}.
\endCD
$$
Finally,
$$\chi_{k-1}|_{{}''\H_{N^y_k}\underset {\Qb^y_k}\times \Qb^y_{k+1,ex}}
=\chi_k|_{{}''\H_{N^y_k}\underset {\Qb^y_k}\times \Qb^y_{k+1,ex}}.$$

\medskip

Now, let us show that for $\F'\in \on{D}^W(\Qb^y_{k+1,ex})$, the object
$W^y_{k,k+1,ex}(\pi_{k+1,ex,k}{}_!(\F'))$
is canonically isomorphic to $\F'$.

Note that as in \lemref{groupoid lifts}, the groupoid 
$'\H^y_{N_k}$ ``lifts'' to $\check \E^i_k$, i.e., we have a Cartesian 
diagram
$$
\CD
\check \E^i_k   @<{\on{act}_{\E^i_k}}<< \E^i_k\underset
{\Qb^y_k}\times \check \E^i_k  \\
@V{{'\check{\on{pr}}^i_k}}VV    @V{\on{id}\times {'\check{\on{pr}}^i_k}}VV
\\
\Qb^y_k @<{'\on{act}^i_k}<<   \E^i_k.
\endCD
$$

Thus we may assume that we start with 
$\F'\in \on{D}(\check \E^i_k)$, which satisfies
\begin{equation} \label{complicated}
\on{act}_{\E^i_k}^*(\F')\simeq ('\on{pr}^i_k\times\on{id})^*(\F')\otimes 
\on{ev}^*(\on{A-Sch}),
\end{equation}
where $'\on{pr}^i_k\times\on{id}$ is the map 
$\E^i_k\underset {\Qb^y_k}\times \check \E^i_k\to \check \E^i_k$
and $\on{ev}$ is the evaluation map
$\E^i_k\underset {\Qb^y_k}\times \check \E^i_k\to \check \AA^1$.

However by looking at another Cartesian square:
$$
\CD
\check \E^i_k   @<{'\on{pr}^i_k\times\on{id}}<< \E^i_k\underset
{\Qb^y_k}\times \check \E^i_k  \\
@V{{'\check{\on{pr}}^i_k}}VV    @V{\on{id}\times {'\check{\on{pr}}^i_k}}VV \\
\Qb^y_k @<{'\on{pr}^i_k}<<   \E^i_k,
\endCD
$$
we obtain that \eqref{complicated} implies that
$$\on{Four}^{-1}(\F')\simeq 
{}'\on{act}^i_k{}^*\left({'\check{\on{pr}}^i_k}{}_!(\F')\right),$$
which is what we had to show. 

The last assertion that $\pi_k{}_!(\F')\to \pi_k{}_*(\F')$
also follows from the above diagram using the fact that the
$!-$ and $*-$ Fourier transforms coincide.

\end{proof}

\ssec{}

Finally, we are ready to construct the functor
$$W_{k,k+1,ex}:\on{D}^W(\Qb_k)\to \on{D}^W(\Qb_{k+1,ex}).$$

The construction of the corresponding functor on the level 
of abelian categories,
i.e., $W_{k,k+1,ex}:\on{P}^W(\Qb_k)\to \on{P}^W(\Qb_{k+1,ex})$
follows immediately from \propref{functors with y}, because
an object $\F\in \on{P}^W(\Qb_{k+1,ex})$ can be glued from its restrictions
$\F|_{\Qb^y_{k+1,ex}}$.

Moreover, \propref{functors with y} implies that
$\pi_{k+1,ex,k}{}_!$ maps $\on{D}^W(\Qb_{k+1,ex})$ to $\on{D}^W(\Qb_k)$
and the induced functor $\on{P}^W(\Qb_{k+1,ex})\to \on{P}^W(\Qb_k)$
is an equivalence inverse to $W_{k,k+1,ex}$.

Thus, it remains to show that 
$\pi_{k+1,ex,k}{}_!:\on{D}^W(\Qb_{k+1,ex})\to \on{D}^W(\Qb_k)$
is an equivalence.

\medskip

First, let us notice that for any substack $U\subset \Qb_k$
of finite type, as in \propref{support of sheaves}, one can find an integer
$i''(U)$ such that for $i\geq i''(U)$ the image of each $U\cap \Qb^y_k$ under the action
of the entire groupoid $\H_{\N^y_k}$ equals that of
$\H_{\N^y_k}^i$. 
Hence, any substack $U$ of finite type can be replaced by a 
bigger one, which is also of finite type, such that each
$U\cap \Qb^y_k$ is $\H_{\N^y_k}$-stable.
Let $U_{k+1,ex}$ denote its preimage in $\Qb_{k+1,ex}$.

\medskip

Then the categories $\on{D}^W(U)$ and $\on{D}^W(U_{k+1,ex})$ make
sense, and we have the functor 
$\pi_{k+1,ex,k}{}_!:\on{D}^W(U_{k+1,ex})\to \on{D}^W(U)$, and
it suffices to show that it is an equivalence.

We claim that this functor is fully-faithful. Since $U$ intersects
only finitely many strata $^{\db}\Qb_k$, it suffices to check
that for two objects $\F_1,\F_2\in \on{D}^W(U_{k+1,ex})$
$$\Hom_{\on{D}({}^{\db}\Qb_{k+1,ex})}(\F_1,\F_2)\to
\Hom_{\on{D}({}^{\db}\Qb_k)}
\left(\pi_{k+1,ex,k}{}_!(\F_1),\pi_{k+1,ex,k}{}_!(\F_2)\right)$$
is an isomorphism. But we know that from 
\lemref{Whittaker functor on strata}.

To finish the proof, we must show that 
$\pi_{k+1,ex,k}{}_!:\on{D}^W(U_{k+1,ex})\to \on{D}^W(U)$ is surjective
on objects. However, this we know, because every object of
$\on{D}^W(U_{k+1,ex})$ is obtained by gluing finitely many
perverse sheaves, and we know already that
$\pi_{k+1,ex,k}{}_!:\on{P}^W(\Qb_{k+1,ex})\to \on{P}^W(\Qb_k)$
is an equivalence.

\ssec{}

Thus, \thmref{Whittaker functors} is proved.

\medskip

We define the functor $W_{k,k+1}:\on{D}^W(\Qb_k)\to 
\on{D}^W(\Qb_{k+1})$ as the composition of $W_{k,k+1,ex}$
followed by the restriction $\on{D}^W(\Qb_{k+1,ex})\to 
\on{D}^W(\Qb_{k+1})$. By construction, $W_{k,k+1}$ is exact.

For two integers $1\leq k< k'\leq n$ we define 
$W_{k,k'}:\on{D}^W(\Qb_k)\to \on{D}^W(\Qb_{k'})$ as the composition
$W_{k'-1,k'}\circ...\circ W_{k,k+1}$. Finally, we set
$W:\on{D}^W(\Bun_n')\to \on{D}^W(\Qb)$ to be $W_{1,n}$. All these
functors are exact.

\ssec{}  \label{GGm equivariance}

Recall that in \secref{base S} we said that the categories
$\on{D}^W(S\times \Qb_k)$, $\on{D}^W(S\times \Qb_{k,ex})$
can be introduced for an arbitrary base scheme $S$.

Similarly, one has the functors
$W_{k,k+1,ex}:\on{D}^W(S\times \Qb_k)\to \on{D}^W(S\times \Qb_{k+1,ex})$,
which are equivalences of categories, and the corresponding functors 
$W_{k,k'}:\on{D}^W(S\times \Qb_k)\to \on{D}^W(S\times \Qb_{k'})$.
All these functors map perverse sheaves to perverse sheaves, and
commute with the Verdier duality.

Moreover, for a morphism of schemes $f:S_1\to S_2$, the $!$- and $*$-
direct and inverse image functors
functors $\on{D}^W(S_1\times \Qb_k)\rightleftarrows 
\on{D}^W(S_2\times \Qb_k)$ and
$\on{D}^W(S_1\times \Qb_{k+1,ex})\rightleftarrows 
\on{D}^W(S_2\times \Qb_{k+1,ex})$ commute in the natural sense
with $W_{k,k'}$.

This is evident for the $*$-inverse image and the $!$-direct image
via the description of the quasi-inverse functor
as $\pi_{k+1,ex,k}{}_!$, and for the $!$-inverse image and the 
$*$-direct image as $\pi_{k+1,ex,k}{}_*$.

\medskip

To conclude, let us remark that the group
$\GG_m$ acts on all the stacks $\Qb_k$ by simultaneously
scaling the maps $\kappa_i$. 
Thus, it makes sense to talk about the equivariant 
derived categories $\on{D}^{\GG_m}(\Qb_k)$.

We introduce the equivariant version of the Whittaker category
$\on{D}^{\GG_m,W}(\Qb_k)$ to be the full triangulated subcategory
of $\on{D}^{\GG_m}(\Qb_k)$ consisting of objects who perverse
cohomologies belong to $\on{D}^W(\Qb_k)$.

Thus, we have an equivariant version of
\thmref{Whittaker functors}, and, in particular, the equivalences
$W_{k,k+1,ex}:\on{D}^{\GG_m,W}(\Qb_k)\to \on{D}^{\GG_m,W}(\Qb_{k+1,ex})$,
and the Whittaker functors
$W_{k,k'}:\on{D}^{\GG_m,W}(\Qb_k)\to \on{D}^{\GG_m,W}(\Qb_{k'})$.

\ssec{} 

The rest of this section will not be used in the sequel. We would like 
to compare the Whittaker functor $W:\on{D}(\Bun_n')\to \on{D}^W(\Qb)$
defined above with another functor of related nature introduced by
G.~Laumon in \cite{La1}.

For an integer $k$, let $\on{Coh}_k$ denote the stack of
coherent sheaves of generic rank $k$; let $\on{Coh}'_k$
denote the stack of pairs: $(\M,\kappa)$, where
$\M\in \on{Coh}_k$, and $\kappa$ is an injective map
of sheaves $\Omega^{k-1}\to \M$. Let, in addition,
$\on{Coh}'_{k,ex}\supset \on{Coh}'_k$ denote the stack of pairs
$(\M,\kappa)$ as before, but where we omit the condition 
that $\kappa$ is injective.

We have a functor 
$W^{\on{Coh}}_{k,k-1}: \on{D}(\on{Coh}'_k)\to \on{D}(\on{Coh}'_{k-1})$.
Namely, note that $\on{Coh}'_k$ and $\on{Coh}'_{k-1,ex}$ form a pair
of mutually dual vector bundles over $\on{Coh}_{k-1}$.
We set $W^{\on{Coh}}_{k,k-1}$ to be the composition of the Fourier
transform functor $\on{D}(\on{Coh}'_k)\to \on{D}(\on{Coh}'_{k-1,ex})$
followed by the restriction $\on{D}(\on{Coh}'_{k-1,ex})\to 
\on{D}(\on{Coh}'_{k-1})$.

By composing, for any $n$ we obtain a functor 
$W^{\on{Coh}}_n:\on{D}(\on{Coh}'_n)\to \on{D}(\on{Coh}'_1)$.

\medskip

Recall now the stack $\wt{\Q}$ of \cite{FGV1}. We have a natural 
smooth projection $\phi^{\on{Coh}}:\wt{\Q}\to \on{D}(\on{Coh}'_1)$,
and a map $\on{ev}:\wt{\Q}\to \AA^1$.

We define the functor $\wt{W}^{\on{Coh}}_n:\on{D}(\on{Coh}'_n)\to
\on{D}(\wt{\Q})$ by
$$\F\mapsto \phi^{\on{Coh}}{}^*\left(W^{\on{Coh}}_n(\F)\right)[d]\otimes
\on{ev}^*(\on{A-Sch}),$$
where $[d]$ is the shift by $\on{dim.rel.}(\wt{\Q},\on{Coh}'_1)$.

\medskip

Recall also that we have a map $\nu:\wt{\Q}\to \Qb=\Q_n$.
Finally note that $\on{Coh}'_n$ contains $\Bun_n'=\Q_1$ as
an open substack.

\begin{prop} \label{Laumon's version}
Let $\F$ be an object of $\on{D}(\Bun_n')$, and let
$\F'$ be its any extension to an object of $\on{D}(\on{Coh}'_n)$.
Then 
$$\nu_!\left(\wt{W}^{\on{Coh}}_n(\F')\right)\simeq
\nu_*\left(\wt{W}^{\on{Coh}}_n(\F')\right)\simeq W(\F).$$
\end{prop}

Instead of giving the proof of this statement, 
we will sketch the argument when $n=2$. The proof in the 
general case follows the same lines.

\ssec{}   \label{Drinfeld's transform}

Consider the following set-up:

Let $\Y$ be a base, and $\E_1$, $\E_2$ be two vector bundles
viewed as group-schemes over $\Y$, and $\on{p}: \E_1\to \E_2$ a map.
Suppose that both $\E_1$ and $\E_2$ act on a scheme $\X$ over
$\Y$, i.e., we have the action maps 
$$\on{act}_i: \E_i\underset{\Y}\times \X\to \X,$$
with $\on{act}_1=\on{act}_2\circ \on{p}$.

Consider the functors $\on{F}_i:\on{D}(\X)\to 
\on{D}(\check\E_i\underset{\Y}\times \X)$, where
$\check\E_i$ is the dual vector bundle, given by
$$\F\mapsto \on{Four}\left(\on{act}_i^*(\F)[d_i]\right),$$
where $d_i=\on{dim.rel.}(\E_i,\Y)$.

Then for $\F\in \on{D}(\X)$, we have:
$$\on{F}_1(\F)\simeq (\check{\on{p}}\times\on{id})_!(\on{F}_2(\F)),$$
where $\check{\on{p}}\times\on{id}$ is the natural map
$\check\E_2\underset{\Y}\times \X\to \check\E_1\underset{\Y}\times \X$.

\medskip

We apply the above observation in the following circumstances:

We set $\Y=\on{Coh}_1$, $\X:=\on{Coh}'_2$. The vector bundle
$\E_2$ is isomorphic to $\X:=\on{Coh}'_2$ itself, i.e., its
fiber at $\L\in \on{Coh}_1$ is the stack of extensions
$$0\to \Omega\to \M\to \L\to 0.$$
The vector bundle $\E_1$ has its fiber over $\L$ as above 
the stack of extensions
$$0\to \Omega\to \M\to \on{det}(\L)\to 0,$$
where $\on{det}(\L)$ is the determinant of $\L$. The map
$\on{p}:\E_1\to \E_2$ comes from the canonical map of sheaves
$\L\to \on{det}(\L)$.

Note that the action of $\E_1$ preserves the open substack
$\Bun'_2\subset \on{Coh}_2'$, and 
$\check\E_1\underset{\on{Coh}_1}\times \Bun'_2\simeq \Qb$.
Moreover, the functor
$$\on{F}_1|_{\Bun'_2}:\on{D}(\Bun'_2)\to 
\check\E_1\underset{\on{Coh}_1}\times \Bun'_2$$
identifies with $W:\on{D}(\Bun'_2)\to \on{D}^W(\Qb)$.

To prove the assertion of the proposition, it suffices to notice
that we have a Cartesian square:
$$
\CD
\check\E_2\underset{\on{Coh}_1}\times \on{Coh}_2'   @<<<
\check\E_2\underset{\on{Coh}_1}\times \Bun'_2  @>{\sim}>> \wt{\Q}  \\
@V{\check{\on{p}}\times\on{id}}VV
@V{\check{\on{p}}\times\on{id}}VV   @V{\nu}VV  \\
\check\E_1\underset{\on{Coh}_1}\times \on{Coh}_2'   @<<<
\check\E_1\underset{\on{Coh}_1}\times \Bun'_2 @>{\sim}>> \Qb.
\endCD
$$

\section{Cuspidality}   \label{cuspidality}

\ssec{}

Let us first recall the notion of cuspidality on $\Bun_n$.

For $n=n_1+n_2$, let $\on{Fl}^n_{n_1,n_2}$ denote the stack
of extensions 
$$0\to \M^1\to \M \to \M^2\to 0,$$
where $\M^i\in \Bun_{n_i}$.

We have the natural projection $\p_{n_1,n_2}:\on{Fl}^n_{n_1,n_2}\to \Bun_n$,
which remembers the middle term of the above short exact sequence,
and the projection $\q_{n_1,n_2}:\on{Fl}^n_{n_1,n_2}\to 
\Bun_{n_1}\times\Bun_{n_2}$, which remembers $(\M^1,\M^2)$.

The projection $\q_{n_1,n_2}$ is in general non-representable,
but is a generalized vector bundle with the fiber over $(\M^1,\M^2)\in 
\Bun_{n_1}\times\Bun_{n_2}$ being the stack of extensions of $\M^2$ by
means of $\M^1$. Therefore, the direct image functors
$\q_{n_1,n_2}{}!:\on{D}(\on{Fl}^n_{n_1,n_2})\to \on{D}(\Bun_n)$
are well-defined.

The constant term functors 
$\on{CT}^n_{n_1,n_2}:\on{D}(\Bun_n)\to \on{D}(\Bun_{n_1}\times \Bun_{n_2})$
are defined by
$$\F\mapsto \q_{n_1,n_2}{}_!\left(\p_{n_1,n_2}^*(\F)\right).$$

\medskip

Recall that an object $\F\in \on{D}(\Bun_n)$ is called cuspidal if
$\on{CT}^n_{n_1,n_2}(\F)=0$ for all $1\leq n_1,n_2<n$.

\medskip

Since the projection $\q$ is not proper, the functor 
$\on{CT}^n_{n_1,n_2}$ does not commute with the Verdier duality.
Therefore, if $\F$ is cuspidal, it will not in general be true 
that $\DD(\F)$ is cuspidal.

\ssec{}

We will now introduce the notion of cuspidality on the stacks
$\Qb_k$. 

For $n_1$ as above and $k\leq n_1$, let $\Qb_{n_1,k}$ denote the stack
classifying the data of $\M^1\in \Bun_{n_1}$,
and a collection of maps
$\kappa_{n_1,i}:\Omega^{n-1+...+n-i}\to \Lambda^i(\M^1)$
for $1\leq i\leq k$, satisfying the Pl\"ucker relations.

\medskip

For $k\leq n_1$, let
$\on{Fl}^{\Qb_k}_{n_1,n_2}$ be the stack classifying the data of
a short exact sequence 
$$0\to \M^1\to \M \to \M^2\to 0 ,$$
as in the definition of $\on{Fl}^n_{n_1,n_2}$, and a collection of maps
$\kappa_i:\Omega^{n-1+...+n-i}\to \Lambda^i(\M^1)$ for $1\leq i\leq k$,
which satisfy the Pl\"ucker relations.

We have a natural map 
$\q_{n_1,n_2,k}:\on{Fl}^{\Qb_k}_{n_1,n_2}\to \Qb_{n_1,k}\times \Bun_{n_2}$,
which makes the following square Cartesian:
$$
\CD
\on{Fl}^{\Qb_k}_{n_1,n_2}   @>{\q_{n_1,n_2,k}}>> 
\Qb_{n_1,k}\times \Bun_{n_2}  \\
@VVV           @VVV   \\
\on{Fl}^n_{n_1,n_2}   @>{\q_{n_1,n_2}}>> \Bun_{n_1}\times \Bun_{n_2}.
\endCD
$$

In addition, we have a map 
$\p_{n_1,n_2,k}:\on{Fl}^{\Qb_k}_{n_1,n_2}\to \Qb_k$.

\medskip

For $k\leq n_1$ we define the constant term functors 
$$\on{CT}^{\Qb_k}_{n_1,n_2}:
\on{D}^W(\Qb_k)\to \on{D}(\Bun_{n_1}\times \Bun_{n_2})$$ by
$\F\mapsto \q_{n_1,n_2,k}{}_!\left(\p_{n_1,n_2,k}^*(\F)\right)$.

We call an object $\F\in \on{D}^W(\Qb_k)$ cuspidal if
$\on{CT}^{\Qb_k}_{n_1,n_2}(\F)=0$ for all $k\leq n_1<n$.

\medskip

In principle, one can introduce the constant term functors also for
$k>n_1$, and properly speaking, a complex $\F\in \on{D}(\Qb_k)$
should be called cuspidal if {\it all} the constant term 
functors vanish when applied to it, including those with $k>n_1$. 
However, for objects of the Whittaker category these other 
functors vanish automatically, so the two notions coincide.

\medskip

Let $\pi$ denote the natural projection $\Qb_1\simeq \Bun_n'\to \Bun_n$.
It is easy to see that for $\F\in\on{D}(\Bun_n)$,
$$\on{CT}^{\Qb_1}_{n_1,n_2}\left(\pi^*(\F)\right)\simeq
(\pi_{n_1}\times \on{id})^*\left(\on{CT}^n_{n_1,n_2}(\F)\right),$$
where $\pi_{n_1}\times \on{id}$ denotes the natural map
$\Qb_{n_1,1}\times \Bun_{n_2}\to \Bun_{n_1}\times\Bun_{n_2}$.

Therefore, if an object $\F\in \on{D}(\Bun_n)$ is cuspidal, then so
is $\pi^*(\F)$.

\ssec{}

The main result of this section is the following theorem:

\begin{thm}  \label{cuspidal orthogonal to degenerate}
Let $\F_1\in \on{D}(\Bun'_n)$ be cuspidal and
$\F_2\in \on{D}(\Bun'_n)$ be any object.
Then the map 
$\Hom_{\on{D}(\Bun'_n)}(\F_1,\F_2)\to 
\Hom_{\on{D}^W(\Qb)}(W(\F_1),W(\F_2))$
is an isomorphism.
\end{thm}

Of course, along with \thmref{cuspidal orthogonal to degenerate}
as it is stated, we have its $\GG_m$-equivariant version, and a version
involving a base $S$, cf. \secref{GGm equivariance}.

\thmref{cuspidal orthogonal to degenerate} follows by induction from the following assertion:

\begin{prop}  \label{Whittaker and constant term functors}
{\em (1)} The functor $W_{k,k+1}:\on{D}^W(\Qb_k)\to \on{D}^W(\Qb_{k+1})$
maps cuspidal objects to cuspidal. 

{\em (2)} If $\F\in \on{D}^W(\Qb_k)$ is cuspidal, then the $*$-restriction of
$W_{k,k+1,ex}(\F)$ to $\Qb_{k+1}-\Qb_{k+1,ex}$ is zero.
\end{prop}

Indeed, to prove \thmref{cuspidal orthogonal to degenerate},
it suffices to show that if we have two objects
$\F_1,\F_2\in \on{D}^W(\Qb_k)$ with $\F_1$ cuspidal, then
$$\Hom_{\on{D}^W(\Qb_k)}(\F_1,\F_2)\to \Hom_{\on{D}^W(\Qb_{k+1})}
(W_{k,k+1}(\F_1),W_{k,k+1}(\F_2))$$ is an isomorphism.

However, by \thmref{Whittaker functors}, we know that
$$\Hom_{\on{D}^W(\Qb_k)}(\F_1,\F_2)\to \Hom_{\on{D}^W(\Qb_{k+1,ex})}
(W_{k,k+1,ex}(\F_1),W_{k,k+1,ex}(\F_2))$$ is an isomorphism. And now,
the condition that 
$$W_{k,k+1,ex}(\F_1)|_{\Qb_{k+1}-\Qb_{k+1,ex}}=0$$ 
means that
\begin{align*}
&\Hom_{\on{D}^W(\Qb_{k+1,ex})}
(W_{k,k+1,ex}(\F_1),W_{k,k+1,ex}(\F_2)):=  \\
&\Hom_{\on{D}(\Qb_{k+1,ex})}
(W_{k,k+1,ex}(\F_1),W_{k,k+1,ex}(\F_2))\simeq \\
&\Hom_{\on{D}(\Qb_{k+1})}
(W_{k,k+1}(\F_1),W_{k,k+1}(\F_2))=: \\
&\Hom_{\on{D}^W(\Qb_{k+1})}
(W_{k,k+1}(\F_1),W_{k,k+1}(\F_2)).
\end{align*}

\ssec{Proof of \propref{Whittaker and constant term functors}(1)}

Let $n_1\geq k+1$. Note that in addition to the stack $\Qb_{n_1,k}$,
one can introduces its ``$ex$'' version $\Qb_{n_1,k+1,ex}$.
Moreover, proceeding just as in \secref{Whittaker categories}
and \secref{Whittaker functs}, we introduce the categories
$\on{D}^W(\Qb_{n_1,k})$, $\on{D}^W(\Qb_{n_1,k+1,ex})$, and the 
functors 
$W_{n_1,k,k+1,ex}:\on{D}^W(\Qb_{n_1,k})\to \on{D}^W(\Qb_{n_1,k+1,ex})$
and $W_{n_1,k,k+1}:\on{D}^W(\Qb_{n_1,k})\to \on{D}^W(\Qb_{n_1,k+1})$.

\medskip

In addition, we can introduce a stack $\on{Fl}^{\Qb_{k+1,ex}}_{n_1,n_2}$,
which fits into the diagram:
$$
\CD
\Qb_{k+1,ex}   @<{\p_{n_1,n_2,k+1,ex}}<<
\on{Fl}^{\Qb_{k+1,ex}}_{n_1,n_2}   @>{\q_{n_1,n_2,k+1,ex}}>> 
\Qb_{n_1,k+1,ex}\times \Bun_{n_2}  \\
@V{\pi_{k+1,ex,k}}VV           @VVV   @V{\pi_{n_1,k+1,ex,k}\times\on{id}}VV \\
\Qb_k @<{\p_{n_1,n_2,k}}<<  
\on{Fl}^{\Qb_k}_{n_1,n_2} @>{\q_{n_1,n_2,k}}>> 
\Qb_{n_1,k}\times \Bun_{n_2}   \\
@VVV       @VVV    @VVV  \\
\Bun_n   @<{\p_{n_1,n_2}}<<
\on{Fl}^n_{n_1,n_2}   @>{\q_{n_1,n_2}}>> \Bun_{n_1}\times \Bun_{n_2}.
\endCD
$$
In this diagram the right portion consists of Cartesian squares.

Using the stack $\on{Fl}^{\Qb_{k+1,ex}}_{n_1,n_2}$
we introduce the functor 
$$\on{CT}^{\Qb_{k+1,ex}}_{n_1,n_2}:
\on{D}^W(\Qb_{k+1,ex})\to \on{D}(\Qb_{n_1,k+1,ex}\times \Bun_{n_2}).$$

\begin{lem}  \label{CT and Whittaker}
The functor $\on{CT}^{\Qb_{k+1,ex}}_{n_1,n_2}$ maps
$\on{D}^W(\Qb_{k+1,ex})$ to 
$\on{D}^W(\Qb_{n_1,k+1,ex}\times \Bun_{n_2})$.
\end{lem}

\begin{proof}

We will use the description of $\on{D}^W(\Qb_{n_1,k+1,ex})$ similar to that
of \propref{description of category}. For a string of integers
$\db=d_1,...,d_k$, let $^{\db}\Qb_{n_1,k+1,ex}{}'\subset {}^{\db}\Qb_{n_1,k+1,ex}$
be the corresponding locally closed substacks of $\Qb_{n_1,k+1,ex}$. Let
also $^{\db}\P_{n_1,k}$ be the stack classifying the data of
$(D_1,...,D_k,\M^1{}',\Omega^{n-k-1}\to \M^1{}')$, as 
in the definition of $^{\db}\P_{k}$, with the difference that now
$\M^1{}'$ is a vector bundle of rank $n_1-k$. We have a smooth
map $\phi_{n_1,k}:{}^{\db}\Qb_{n_1,k+1,ex}{}'\to {}^{\db}\P_{n_1,k}$.

To prove the lemma it is sufficient to show that for
$\F\in \on{D}^W(\Qb_{k+1,ex})$, the restriction of
$\on{CT}^{\Qb_{k+1,ex}}_{n_1,n_2}(\F)$ to each
$^{\db}\Qb_{n_1,k+1,ex}{}'$ is isomorphic to the pull-back
of a complex on $^{\db}\P_{n_1,k}$, tensored by an appropriate
Artin-Schreier sheaf.

Consider the fiber product $$Z:=\on{Fl}^{n-k}_{n_1-k,n_2}\underset{\Bun_{n_1-k}}
\times {}^{\db}\P_{n_1,k}.$$ 
Let $^{\db}\on{Fl}^{\Qb_{k+1,ex}}_{n_1,n_2}$ be the preimage in 
$\on{Fl}^{\Qb_{k+1,ex}}_{n_1,n_2}$ of the substack
$^{\db}\Qb_{k+1,ex}:"'\subset \Qb_k$ under $\p_{n_1,n_2}$.
We have a commutative diagram
$$
\CD
^{\db}\Qb_{k+1,ex}{}' @<{\p_{n_1,n_2}}<< ^{\db}\on{Fl}^{\Qb_{k+1,ex}}_{n_1,n_2} @>{\q_{n_1,n_2}}>> 
^{\db}\Qb_{n_1,k+1,ex}{}'\times \Bun_{n_2} \\
@V{\phi_k}VV  @VVV  @V{\phi_{n_1,k}\times \on{id}}VV  \\
^{\db}\P_k @<<< Z @>>> ^{\db}\P_{n_1,k}\times \Bun_{n_2}.
\endCD
$$

The right portion of this diagram is not Cartesian. However, the map
$$^{\db}\on{Fl}^{\Qb_{k+1,ex}}_{n_1,n_2}\to
Z\underset{^{\db}\P_{n_1,k}\times \Bun_{n_2}}\times 
({}^{\db}\Qb_{n_1,k+1,ex}{}'\times \Bun_{n_2})$$ is smooth with contractible
fibers.

Hence, the assertion of the lemma follows from the projection formula.

\end{proof}

Taking into account the above lemma,
part (1) of \propref{Whittaker and constant term functors} would
follow once we are able to establish an isomorphism of functors:
\begin{equation}  \label{CT and W}
\on{CT}^{\Qb_{k+1,ex}}_{n_1,n_2}\circ W_{k,k+1,ex}\simeq
(W_{n_1,k,k+1,ex}\times\on{id})\circ \on{CT}^{\Qb_k}_{n_1,n_2},
\end{equation}
both of which map from $\on{D}^W(\Qb_k)$ to 
$\on{D}^W(\Qb_{n_1,k+1,ex}\times \Bun_{n_2})$.

\medskip

Let us observe that the functor $\on{CT}^{\Qb_k}_{n_1,n_2}$,
has a natural right adjoint,
which we will denote by $\on{Eis}^{\Qb_k}_{n_1,n_2}$,
that maps
$\F\in \on{D}(\Qb_{n_1,k+1}\times \Bun_{n_2})$ to
$\p_{n_1,n_2,k}{}_*\left(\q_{n_1,n_2,k}^!(\F)\right)$.
This functor also maps
$\on{D}^W(\Qb_{n_1,k}\times \Bun_{n_2})$ to $\on{D}^W(\Qb_k)$.

And similarly, we have a right adjoint of $\on{CT}^{\Qb_{k+1,ex}}_{n_1,n_2}$
$$\on{Eis}^{\Qb_{k+1,ex}}_{n_1,n_2}:
\on{D}^W(\Qb_{n_1,k+1,ex}\times \Bun_{n_2})\to \on{D}^W(\Qb_{k+1,ex}).$$

To prove \eqref{CT and W} it suffices to verify the isomorphism
on the level of the corresponding adjoint functors. In other words, we must show
that
$$\pi_{k+1,ex,k}{}_*\circ \on{Eis}^{\Qb_{k+1,ex}}_{n_1,n_2}\simeq
\on{Eis}^{\Qb_k}_{n_1,n_2}\circ (\pi_{n_1,k+1,ex,k}\times\on{id})_*.$$

However, the latter isomorphism follows from base change.

\ssec{Proof of \propref{Whittaker and constant term functors}(2)}
\label{proof of point 2}

Note that $\Qb_{k+1,ex}-\Qb_{k+1}\subset \Qb_{k+1}$ is naturally
isomorphic to $\Qb_k$. We would like to calculate 
$W_{k,k+1,ex}(\F)|_{\Qb_k}$ in terms of $\on{CT}^{\Qb_k}_{n_1,n_2}$
for $n_1=k$.

Recall that to a string of integers $\db=d_1,...,d_k$ we associated
a locally closed substack $^{\db}\Qb_k\subset \Qb_k$.

Note now that we have a natural map $\psi_k:{}^{\db}\Qb_k\to 
\Qb_{k,k}\times \Bun_{n-k}$. Namely, we can think of a point of
$^{\db}\Qb_k$ as a data of
$$0=\M_0\subset \M_1\subset...\subset\M_k\subset \M,$$
and identifications
$\M_i/\M_{i-1}\simeq \Omega^{n-i}(D_i-D_{i-1})$ for
$(D_1,...,D_k) \in X^{(\db)}$.

The corresponding point of $\Qb_{k,k}\times \Bun_{n-k}$
is $\M^1=\M_k$, with the data of $\kappa_{n_1,i}$ being given
by the old $\kappa_i$'s, and $\M^2:=\M/\M_k$.

We claim that up to a cohomological shift, for $\F\in \on{D}^W(\Qb_k)$,
\begin{equation}  \label{degenerate coefficient as a constant term}
W_{k,k+1,ex}(\F)|_{^{\db}\Qb_k}\simeq 
\psi_k^*(\on{CT}^{\Qb_k}_{k,n-k}(\F)).
\end{equation}

This follows immediately from the description of the functor
$^{\db}W_{k,k+1,ex}$ in \secref{explicit on strata}.

Thus, part (2) of \propref{Whittaker and constant term functors} follows,
because to show that $W_{k,k+1,ex}(\F)|_{\Qb_k}=0$ for $\F$ cuspidal, it
is enough to show that for all $\db$ as above $W_{k,k+1,ex}(\F)|_{^{\db}\Qb_k}=0$,
and the latter is given by \eqref{degenerate coefficient as a constant term}.

\medskip

Note that in the course of the proof we have shown that
$W_{k,k+1,ex}(\F)|_{\Qb_k}=0$ {\it if and only if} 
$\on{CT}^{\Qb_k}_{k,n-k}(\F)=0$.
This is because the stack $\Qb_{k,k}$ is also stratified by means
of $^{\db}\Qb_{k,k}$, and for every $\db$ the map
$$\psi: {}^{\db}\Qb_k\to {}^{\db}\Qb_{k,k}\times \Bun_{n-k}$$
is surjective.

\ssec{}   \label{degenerate sheaves on Bun'_n}

Thus, \thmref{cuspidal orthogonal to degenerate} is proved.
We will now give another categorical interpretation of it.

\medskip

Let $\on{D}^W_{cusp}(\Qb_k)$ denote the full subcategory 
consisting of cuspidal objects in $\on{D}^W(\Qb_k)$. This is
evidently a triangulated subcategory in $\on{D}^W(\Qb_k)$.

Now, let $\on{D}^W_{degen}(\Qb_k)\subset \on{D}^W(\Qb_k)$
denote the (full triangulated) subcategory of those objects $\F$ for
which $W_{k,n}(\F)=0$.
Let $\Dt^W(\Qb_k)$ denote the quotient triangulated category
$\on{D}^W(\Qb_k)/\on{D}^W_{degen}(\Qb_k)$.

Consider the composition
$$\on{D}^W_{cusp}(\Qb_k)\to \on{D}^W(\Qb_k)\to \Dt^W(\Qb_k).$$

\begin{thm}  \label{as in automorphic forms}
{\em (1)} The above functor
$\on{D}^W_{cusp}(\Qb_k)\to \Dt^W(\Qb_k)$
is an equivalence of categories. 

{\em (2)}
The functor
$\on{D}^W_{cusp}(\Qb_k) \overset{W_{k,n}}\longrightarrow \on{D}^W(\Qb)$
is an equivalence as well.  
\end{thm}

\ssec{Proof of \thmref{as in automorphic forms}}

Let $W_{k,k+1}^{-1}:\on{D}^W(\Qb_{k+1})\to \on{D}^W(\Qb_{k})$
be defined by sending 
$\F'\in \on{D}^W(\Qb_{k+1})$ to 
$\pi_{k+1,ex,k}{}_!(\F'')$, where $\F''$ is
the $!$-extension from $\Qb_{k+1}$ to $\Qb_{k+1,ex}$.
Since $W_{k,k+1,ex}$ is an equivalence,
we have an isomorphism of functors
$$W_{k,k+1}\circ W_{k,k+1}^{-1}\simeq \on{id}_{\on{D}^W(\Qb_{k+1})},$$
and an adjunction map
$W_{k,k+1}^{-1}\circ W_{k,k+1}\to \on{id}_{\on{D}^W(\Qb_k)}$.

\medskip

By construction, $W_{k,k+1}$ induces a functor
$\Dt^W(\Qb_k)\to \Dt^W(\Qb_{k+1})$. We claim this
functor is an equivalence for every $k$.

Indeed, it is easy to see that the functor
$$\on{D}^W(\Qb_{k+1})\overset{W_{k,k+1}^{-1}}\longrightarrow
\on{D}^W(\Qb_{k})\to \Dt^W(\Qb_k)$$
factors through $\Dt^W(\Qb_{k+1})$ and defines a quasi-inverse
for $W_{k,k+1}$.

Hence, $W_{k,n}:\Dt^W(\Qb_k)\to \Dt^W(\Qb_n)=\on{D}^W(\Qb)$
is an equivalence as well.

\medskip

Thus, it remains to prove the first assertion of the theorem.
For that it is enough to show that
$W_{k,k+1}$ induces an equivalence 
$\on{D}^W_{cusp}(\Qb_k)\to \on{D}^W_{cusp}(\Qb_{k+1})$
for every $k$. The fact that the image of $\on{D}^W_{cusp}(\Qb_k)$
under $W_{k,k+1}$ belongs to $\on{D}^W_{cusp}(\Qb_k)$ was proved in
\propref{Whittaker and constant term functors}.

We claim that $W_{k,k+1}^{-1}$ defines a quasi-inverse.
Indeed, for $\F\in\on{D}^W_{cusp}(\Qb_{k+1})$
to show that $W_{k,k+1}^{-1}(\F)\in \on{D}^W_{cusp}(\Qb_k)$ 
we must verify that
$\on{CT}^{\Qb_k}_{n_1,n_2}(W_{k,k+1}^{-1}(\F))=0$ for $n_1\geq k$.

Suppose first that $n_1\geq k+1$. Then, since $W_{n_1,k,k+1,ex}$
is an equivalence, what we need follows immediately from \eqref{CT and W}.
For $n_1=k$, the needed assertion follows from the last remark
of \secref{proof of point 2}.

\medskip

The fact that $W_{k,k+1}\circ W_{k,k+1}^{-1}\simeq \on{id}$ 
we know already. It remains, therefore to show that for
$\F\in \on{D}^W_{cusp}(\Qb_k)$,
$$W_{k,k+1}^{-1}(W_{k,k+1}(\F))\to \F$$
is an isomorphism.
Let $\F'$ be the cone of the above map. We know that
$\F'\in \on{D}^W_{cusp}(\Qb_k)$, and $W_{k,k+1}(\F')\simeq 0$.
Hence, $\F'\simeq 0$ by \thmref{cuspidal orthogonal to degenerate}.

\ssec{}

As a corollary of \thmref{as in automorphic forms} we obtain
that the category $\on{D}^W_{cusp}(\Qb_k)$, and hence, in
particular $\on{D}_{cusp}(\Bun_1)$, possesses a $t$-structure.
Indeed, it is equivalent to the category $\on{D}^W(\Qb)$,
for which the $t$-structure is manifest.

Note that this $t$-structure does not coincide with
the $t$-structure on the ambient category $\on{D}^W(\Qb_k)$.

\section{The Hecke functors}   \label{Hecke functors}

\ssec{}   \label{intr Hecke}

Recall the Hecke functor
$\on{H}:\on{D}(\Bun_n)\to\on{D}(X\times \Bun_n)$,
which was defined using the stack $\H=\on{Mod}^1_n$.
In this section we will introduce Hecke functors that map
from $\on{D}(\Qb_k)$ to $\on{D}(X\times \Qb_k)$. First
we will consider the case of $\Qb_1=\Bun_n'$.
 
\medskip

Set $\H^{\Bun_n'}:=\Bun_n'\underset{\Bun_n}\times \H$,
where we use the map $\hl:\H\to \Bun_n$ to define the
fiber product.

We have a commutative diagram
$$
\CD
X\times \Bun_n'  @<{s^{\Bun_n'}\times \hl{}^{\Bun_n'}}<< 
\H^{\Bun_n'} @>{\hr{}^{\Bun_n'}}>> 
\Bun_n'  \\
@V{\on{id}\times \pi}VV   @VVV   @V{\pi}VV   \\
X\times \Bun_n @<{s\times \hl}<< \H @>{\hr}>>  \Bun_n,
\endCD
$$
in which the left square is Cartesian. Indeed, the map
$\hr{}^{\Bun_n'}$ attaches to a point 
$(x,\M\hookrightarrow \M', \kappa:\Omega^{n-1}\to \M)\in 
\Bun_n'\underset{\Bun_n}\times \H$ the point
$(\M',\kappa':\Omega^{n-1}\to \M')$, where
$\kappa'$ is the composition $\Omega^{n-1}\overset{\kappa}\to \M\to \M'$.

\medskip

We define the functor $\on{H}^{\Bun_n'}:\on{D}(\Bun_n')\to 
\on{D}(X\times\Bun_n')$ by 
$$\F\mapsto (s^{\Bun_n'}\times \hl{}^{\Bun_n'})_!
\left(\hr{}^{\Bun_n'}{}^*(\F)\right)[n-1].$$

Note that the functors $\on{H}^{\Bun_n'}$ and $\on{H}$ are compatible
in the following way: for $\F\in \on{D}(\Bun_n)$ we have
\begin{equation}  \label{pull back and Hecke}
(\on{id}\times \pi)^*\left(\on{H}(\F)\right)\simeq
\on{H}^{\Bun_n'}\left(\pi^*(\F)\right)[1].
\end{equation}

Note also that since the map $\hr{}^{\Bun_n'}$ is not smooth,
the functor $\on{H}^{\Bun_n'}$ does not commute with the Verdier
duality. In particular, one could define its Verdier twin by
$\F\mapsto (s^{\Bun_n'}\times \hl{}^{\Bun_n'})_*
\left(\hr{}^{\Bun_n'}{}^!(\F)\right)[1-n]$.

\ssec{}

For $1\leq k\leq n$ we introduce the appropriate Hecke functors
in a similar fashion. Namely,
we set $\H^{\Qb_k}:=\Qb_k\underset{\Bun_n}\times \H$, which fits
into a commutative diagram
$$
\CD
X\times \Qb_k  @<{s^{\Qb_k}\times \hl{}^{\Qb_k}}<< 
\H^{\Qb_k} @>{\hr{}^{\Qb_k}}>> 
\Bun_n'  \\
@V{\on{id}\times \pi}VV   @VVV   @V{\pi}VV   \\
X\times \Bun_n @<{s\times \hl}<< \H @>{\hr}>>  \Bun_n,
\endCD
$$
in which the left square is Cartesian.

The functor
$\on{H}^{\Qb_k}:\on{D}(\Qb_k)\to \on{D}(X\times \Qb_k)$ is defined
by means of
$$\F\mapsto (s^{\Qb_k}\times \hl{}^{\Qb_k})_!
\left(\hr{}^{\Qb_k}{}^*(\F)\right)[n-1].$$

Set $_x\H^{\Qb_k}$ to be the preimage of $x\in X$ in $\H^{\Qb_k}$.

For a point $x\in X$ we will denote by $_x\on{H}^{\Qb_k}$ the functor
$\on{D}(\Qb_k)\to \on{D}(\Qb_k)$ obtained as a composition
of $\on{H}^{\Qb_k}$ followed by the $*$-restriction to
$x\times \Qb_k\subset X\times \Qb_k$.

In other words, $_x\on{H}^{\Qb_k}$ can be defined using the substack
$_x\H^{\Qb_k}$ as a correspondence.

\medskip

In a similar way we define the stack $\H^{\Qb_{k+1,ex}}$
and the corresponding functor
$\on{H}^{\Qb_{k+1,ex}}:\on{D}(\Qb_{k+1,ex})\to 
\on{D}(X\times \Qb_{k+1,ex})$.

\ssec{}

\begin{prop}  \label{Hecke preserves Whittaker}
The functor $\on{H}^{\Qb_{k+1,ex}}$ maps 
$\on{D}^W(\Qb_{k+1,ex})$ to $\on{D}^W(X\times \Qb_{k+1,ex})$.
\end{prop}

Of course, as a corollary of this proposition we obtain that
$\on{H}^{\Qb_k}$ maps $\on{D}^W(\Qb_k)$ $\on{D}^W(X\times\Qb_k)$.

\begin{proof}

To simplify the notation we will show that for any $x\in X$,
the functor 
$_x\on{H}^{\Qb_{k+1,ex}}:\on{D}(\Qb_{k+1,ex})\to 
\on{D}(\Qb_{k+1,ex})$ preserves the subcategory $\on{D}^W(\Qb_{k+1,ex})$.

\medskip

Let $y\in X$ be a point different from $x$. It is easy to see that
we have a well-defined functor
$_x\on{H}^{\Qb^y_{k+1,ex}}:\on{D}(\Qb^y_{k+1,ex})\to 
\on{D}(\Qb^y_{k+1,ex})$, constructed using the stack that
we will denote by $_x\H^{\Qb^y_{k+1,ex}}$. 
We will first show that this functor preserves $\on{D}^W(\Qb^y_{k+1,ex})$.

However, this is almost obvious from the definitions:

Recall the groupoid $\H_{N^y_k}\underset{\Qb^y_k}\times \Qb^y_{k+1,ex}$ 
acting on $\Qb^y_{k+1,ex}$. We claim that it lifts to the stack
$_x\H^{\Qb^y_{k+1,ex}}$, i.e., we have a groupoid 
$_x\H^{\Qb^y_{k+1,ex}}_{N^y_k}$ which fits into two commutative diagrams
$$
\CD
_x\H^{\Qb^y_{k+1,ex}}  @<{\on{act}}<< _x\H^{\Qb^y_{k+1,ex}}_{N^y_k}
@>{\on{pr}}>>  _x\H^{\Qb^y_{k+1,ex}}  \\
@V{\hl{}^{\Qb^y_{k+1,ex}}}VV   @V{\hl}VV    
@V{\hl{}^{\Qb^y_{k+1,ex}}}VV   \\
\Qb^y_{k+1,ex}  @<{\on{act}_{k,ex}}<< 
\H_{N^y_k}\underset{\Qb^y_k}\times \Qb^y_{k+1,ex}
@>{\on{pr}_{k,ex}}>>  \Qb^y_{k+1,ex},
\endCD
$$
and
$$
\CD
_x\H^{\Qb^y_{k+1,ex}}  @<{\on{act}}<< _x\H^{\Qb^y_{k+1,ex}}_{N^y_k}
@>{\on{pr}}>>  _x\H^{\Qb^y_{k+1,ex}}  \\
@V{\hr{}^{\Qb^y_{k+1,ex}}}VV   @V{\hr}VV    
@V{\hr{}^{\Qb^y_{k+1,ex}}}VV   \\
\Qb^y_{k+1,ex}  @<{\on{act}_{k,ex}}<< 
\H_{N^y_k}\underset{\Qb^y_k}\times \Qb^y_{k+1,ex}
@>{\on{pr}_{k,ex}}>>  \Qb^y_{k+1,ex},
\endCD
$$
in both of which both squares are Cartesian.
Moreover, the compositions
$$_x\H^{\Qb^y_{k+1,ex}}_{N^y_k}
\overset{\hl}\longrightarrow 
\H_{N^y_k}\underset{\Qb^y_k}\times \Qb^y_{k+1,ex}\overset{\chi_k}\to\AA^1$$
and
$$_x\H^{\Qb^y_{k+1,ex}}_{N^y_k}
\overset{\hr}\longrightarrow 
\H_{N^y_k}\underset{\Qb^y_k}\times \Qb^y_{k+1,ex}\overset{\chi_k}\to\AA^1$$
coincide.

Therefore, if an object $\F\in \on{D}(\Qb^y_{k+1,ex})$ satisfies
the equivariance condition \eqref{N_k equivariance}, then so does
$\hl{}^{\Qb^y_k}_! \left(\hr{}^{\Qb^y_k}{}^*(\F)\right)[n-1]$.

\medskip

Let now $\F$ be an arbitrary object of $\on{D}^W(\Qb_{k+1,ex})$.
To show that $_x\on{H}^{\Qb_{k+1,ex}}(\F)$ also belongs to
$\on{D}^W(\Qb_{k+1,ex})$, from \lemref{general equivariance}
it follows that it is sufficient to show that any irreducible
sub-quotient of any perverse cohomology sheaf of 
$_x\on{H}^{\Qb_{k+1,ex}}(\F)$ belongs to $\on{D}^W(\Qb_{k+1,ex})$.

Let $\K$ be such a sub-quotient. Then there exists $y\in X$, such
that the restriction of $\K$ to $\Qb^y_{k+1,ex}$ is non-zero.
Hence, again by \lemref{general equivariance} and 
\corref{restrictions well-behaved}, it suffices to show that
$\K|_{\Qb^y_{k+1,ex}}$ belongs to $\on{D}^W(\Qb^y_{k+1,ex})$.
But above we have shown that the entire 
$_x\on{H}^{\Qb_{k+1,ex}}(\F)|_{\Qb^y_{k+1,ex}}$ 
belongs to $\on{D}^W(\Qb^y_{k+1,ex})$, and hence also $\K|_{\Qb^y_{k+1,ex}}$,
which is its sub-quotient.

\end{proof}

\ssec{}

Our next goal is to show that the Hecke functors and Whittaker 
functors commute with each other.

\begin{prop}  \label{Hecke and Whittaker commute} 
We have a natural isomorphism of functors
$$\on{H}^{\Qb_{k+1,ex}}\circ W_{k,k+1,ex}\simeq
(\on{id}\times W_{k,k+1,ex})\circ \on{H}^{\Qb_k}:
\on{D}^W(\Qb_k)\to \on{D}^W(X\times \Qb_{k+1,ex}).$$
\end{prop}

Of course, the proposition implies that the functors
$\on{H}^{\Qb_{k+1}}\circ W_{k,k+1}$ and
$(\on{id}\times W_{k,k+1})\circ \on{H}^{\Qb_k}$ from
$\on{D}^W(\Qb_k)$ to $\on{D}^W(X\times \Qb_{k+1})$ are 
isomorphic.

\begin{proof}

As in the proof of the previous proposition, in order
to simplify the notation, we will consider the functors
$_x\on{H}^{\Qb_{k+1,ex}}$ and $_x\on{H}^{\Qb_k}$ instead of
$\on{H}^{\Qb_{k+1,ex}}$ and $\on{H}^{\Qb_k}$.

In fact, from the proof of \propref{Hecke preserves Whittaker}
given above one can directly deduce that for $y\neq x$, 
$_x\on{H}^{\Qb^y_{k+1,ex}}\circ W^y_{k,k+1,ex}\simeq
W^y_{k,k+1,ex}\circ {}_x\on{H}^{\Qb^y_k}$, using the definition
of $W^y_{k,k+1,ex}$ via the Fourier transform functor as in 
\secref{construction of whittaker functors}. We will proceed 
differently. Namely, we will prove that for
$\F\in \on{D}^W(\Qb_{k+1,ex})$,
\begin{equation}  \label{expression for Hecke and Whittaker}
\pi_{k+1,ex,k}{}_!\left({}_x\on{H}^{\Qb_{k+1,ex}}(\F)\right)\simeq
{}_x\on{H}^{\Qb_k}\left(\pi_{k+1,ex,k}{}_!(\F)\right),
\end{equation}
which is equivalent to the statement of 
\propref{Hecke and Whittaker commute}, since $\pi_{k+1,ex,k}{}_!$
induces an equivalence of categories.

\medskip

For a point $x\in X$,
let $\Qb_{k+1,ex,x}$ denote the stack that classifies the data
$(\M,\kappa_1,...,\kappa_k,\\ \kappa_{k+1})$ as before, with the
difference that now the last map 
$\kappa_{k+1}:\Omega^{n-1+...+n-(k+1)}\to \Lambda^{k+1}(\M)$
is allowed to have a simple pole at $x$.
We have a natural closed embedding $\Qb_{k+1,ex}\hookrightarrow
\Qb_{k+1,ex,x}$.

Let $_x\on{H}^{\Qb_{k+1,ex,x}}$ denote the Cartesian product
$$_x\on{H}^{\Qb_{k+1,ex,x}}:={}_x\on{H}^{\Qb_k}\underset{\Qb_k}\times
\Qb_{k+1,ex},$$
where we have used the map 
$\hr{}^{\Qb_k}:{}_x\on{H}^{\Qb_k}\to \Qb_k$ to define the product.

We have a commutative diagram:
$$
\CD
\Qb_{k+1,ex,x}  @<{\hl{}^{\Qb_{k+1,ex,x}}}<<  
_x\on{H}^{\Qb_{k+1,ex,x}}  @>{\hr{}^{\Qb_{k+1,ex,x}}}>> \Qb_{k+1,ex}   \\
@V{\pi_{k+1,ex,k,x}}VV   @VVV   @V{\pi_{k+1,ex,k}}VV   \\
\Qb_k  @<{\hl{}^{\Qb_k}}<<  _x\on{H}^{\Qb_k}  @>{\hr{}^{\Qb_k}}>> \Qb_k, 
\endCD
$$
in which the right square is Cartesian.

By base change, for $\F\in \on{D}^W(\Qb_{k+1,ex})$,
the right-hand side of \eqref{expression for Hecke and Whittaker}
equals
\begin{equation} \label{another expression for Hecke and Whittaker}
(\pi_{k+1,ex,k,x})_!
\biggl(\hl{}^{\Qb_{k+1,ex,x}}_!\left(\hr{}^{\Qb_{k+1,ex,x}}{}^*(\F)\right)\biggr).
\end{equation}

\begin{lem}  \label{support after Hecke}
For $\F\in \on{D}^W(\Qb_{k+1,ex})$, 
the object 
$$\hl{}^{\Qb_{k+1,ex,x}}_!\left(\hr{}^{\Qb_{k+1,ex,x}}{}^*(\F)\right)
\in \on{D}(\Qb_{k+1,ex,x})$$ is supported on $\Qb_{k+1,ex}$.
\end{lem}

\begin{proof}

For $y\neq x$ let $\Qb^y_{k+1,ex,x}$ denote the open substack of
$\Qb_{k+1,ex,x}$ equal to the preimage of $\Qb^y_k$ under
$\pi_{k+1,ex,k,x}$. It would be sufficient to show that for
any such $y$, the restriction of 
$\hl{}^{\Qb_{k+1,ex,x}}_!\left(\hr{}^{\Qb_{k+1,ex,x}}{}^*(\F)\right)$
(as in the lemma) to $\Qb^y_{k+1,ex,x}$
is supported on $\Qb^y_{k+1,ex}$.

As in \secref{def of cat with y} we can introduce the category
$\on{D}^W(\Qb^y_{k+1,ex,x})$, and, as in \propref{Hecke preserves Whittaker}, we 
show that the Hecke functor
$\F\mapsto \hl{}^{\Qb_{k+1,ex,x}}_!\left(\hr{}^{\Qb_{k+1,ex,x}}{}^*(\F)\right)$
maps $\on{D}^W(\Qb^y_{k+1,ex})$ to $\on{D}^W(\Qb^y_{k+1,ex,x})$.

However, we claim that every object of the category $\on{D}^W(\Qb^y_{k+1,ex,x})$
is supported on $\Qb^y_{k+1,ex}$. This is done by introducing a
stratification on $\Qb^y_{k+1,ex,x}$ analogous to the stratification by
$^{\db}\Qb^y_{k+1,ex}$ on $\Qb^y_{k+1,ex}$ and using an analog of 
\propref{description of category}(1).

\end{proof}

To finish the proof of \propref{Hecke and Whittaker commute}, let us
observe that there is another diagram:
$$
\CD
\Qb_k  @<{\pi_{k+1,ex,k}}<< 
\Qb_{k+1,ex} @<{\hl{}^{\Qb_{k+1,ex}}}<<  _x\on{H}^{\Qb_{k+1,ex}} 
@>{\hr{}^{\Qb_{k+1,ex}}}>> \Qb_{k+1,ex}  \\
@V{\on{id}}VV   @VVV    @VVV   @V{\on{id}}VV \\
\Qb_k  @<{\pi_{k+1,ex,k,x}}<<
\Qb_{k+1,ex,x}  @<{\hl{}^{\Qb_{k+1,ex,x}}}<<  _x\on{H}^{\Qb_{k+1,ex,x}}  
@>{\hr{}^{\Qb_{k+1,ex,x}}}>> \Qb_{k+1,ex},  
\endCD
$$
in which the middle square is Cartesian.

Therefore, using \lemref{support after Hecke}, the expression in
\eqref{another expression for Hecke and Whittaker} can be rewritten as
$$(\pi_{k+1,ex,k})_!
\biggl(\hl{}^{\Qb_{k+1,ex}}_!\left(\hr{}^{\Qb_{k+1,ex}}{}^*(\F)\right)\biggr),
$$
which equals the expression in the left-hand side of 
\eqref{expression for Hecke and Whittaker}.

\end{proof}

\ssec{}  \label{sect Hecke right exact}

The following theorem is one of the main technical results
of this paper:

\begin{thm}  \label{Hecke is right-exact}
The functor $\on{H}^{\Qb_n}:\on{D}(\Qb_n)\to \on{D}(X\times \Qb_n)$
is right-exact.
\end{thm}

The rest of this section is devoted to the proof of this theorem.

\medskip

Let us restrict our attention to the connected component of
$\Qb_n$ corresponding to vector bundles $\M$ of a fixed degree.
We set $d=\on{deg}(\Lambda^{n}(\M))-\on{deg}(\Omega^{n-1+n-2+...+1+0})$.
According to the conventions of \cite{FGV1}, the corresponding connected component 
of $\Bun_n$ is denoted by $\Bun_n^d$, and we keep a similar notation for
$\Qb^d_n$.

The data of $\kappa_n$ in the definition of $\Qb_n$ defines a map
$\tau_d:\Qb^d_n\to X^{(d)}$.
Observe that we have a commutative diagram:
$$
\CD
\H^{\Qb_n}  @>{\hr{}^{\Qb_n}}>>  \Qb^{d+1}_n  \\
@V{s^{\Qb_n}\times \hl{}^{\Qb_n}}VV    @V{\tau_{d+1}}VV  \\
X\times \Qb^d_n @>>>  X^{(d+1)},
\endCD
$$
where the bottom horizontal arrow is the composition
$$X\times \Qb^d_n\overset{\on{id}\times \tau_d}\longrightarrow
X\times X^{(d)}\to X^{(d+1)}.$$

From the above diagram we obtain the following

\begin{lem}  \label{preimage is contained in support}
For a given point $(\M',\kappa'_1,...,\kappa'_n)\in \Qb^{d+1}_n$, 
its preimage in $\H^{\Qb_n}$ is contained in
$$\underset{x\in \on{supp}(D'_n)}\bigcup\, {}_x\H^{\Qb_n},$$
where $D'_n\in X^{(d+1)}$ is the image of the above point under 
$\tau_{d+1}$.
\end{lem}

The proof of \thmref{Hecke is right-exact} will be obtained 
from the following general result:

Let
$$\Y  \overset{f}\leftarrow Z \overset{f'}\to \Y'$$
be a diagram of stacks with the morphism $f$ being representable.
Suppose that $Z$ can be decomposed into locally closed substacks
$Z=\cup \, Z_\alpha$ (the decomposition being locally finite)
such that if we denote by $m_\alpha$ (resp., $m'_\alpha$) the maximum
of the dimensions of fibers of $f:Z_\alpha\to \Y$ (resp.,
$f':Z_\alpha\to \Y'$), we have:
$$m_\alpha+m'_\alpha\leq m$$
for some integer $m$.

\begin{lem}
Under the above circumstances, the functor
$\on{D}(\Y')\to \on{D}(\Y)$ given by
$$\F\mapsto f_!\left(f'{}^*(\F)\right)$$
sends objects of $\on{D}(\Y')^{\leq 0}$ to $\on{D}(\Y)^{\leq m}$.
\end{lem}

The proof of the lemma follows from the definition of the perverse
$t$-structure.

We apply this lemma for $\Y=X\times \Qb_n$, $\Y'=\Qb_n$,
$Z=\H^{\Qb_n}$, $f=s^{\Qb_n}\times \hl{}^{\Qb_n}$, $f'=\hr{}^{\Qb_n}$,
and $m=n-1$. Thus, our task is to find
a suitable stratification of $\H^{\Qb_n}$.

\ssec{}

For two strings of non-negative integers
$\db^1=d^1_1,...,d^1_n$, $\db^2=d^2_1,...,d^2_n$ with
$d^2_n=d^1_n+1$, and $d^1_i\leq d^2_i\leq d^1_i+1$,
let $^{\db^1,\db^2}\H^{\Qb_n}$ denote the following
locally closed substack of $\H^{\Qb_n}$:

Recall that $\H^{\Qb_n}$ classifies the data of
$$(x\in X,\,\,\M\in \Bun_n, \,\,
\kappa_i:\Omega^{n-1+...+n-i}\to \Lambda^i(\M),\,\,\M'\in \Bun_n,\,\,
\beta:\M\hookrightarrow \M'),$$
with $\M'/\M$ being a skyscraper at $x$.

We say that a point as above belongs to $^{\db^1,\db^2}\H^{\Qb_n}$
if

\noindent {(a)}
Each map $\kappa_i:\Omega^{n-1+...+n-i}\to \Lambda^i(\M)$ has a
zero of order $d^1_i$ at $x$.

\noindent {(b)}
Each composed map
$\kappa'_i:\Omega^{n-1+...+n-i}\to \Lambda^i(\M')$ has a 
zero of order $d^2_i$ at $x$.

\medskip

As in the case of $\Qb_k=\cup\, {}^{\db}\Qb_k$, it is easy 
to show that the substacks $^{\db^1,\db^2}\H^{\Qb_n}$ define
a locally finite decomposition of $\H^{\Qb_n}$ into locally closed
substacks.

We now need to verify the estimate on the dimensions of fibers
$^{\db^1,\db^2}\H^{\Qb_n}$ under the maps $(s^{\Qb_n}\times \hl{}^{\Qb_n})$ 
and $\hr{}^{\Qb_n}$.

\medskip

Let $^{\db^1,\db^2}_x\H^{\Qb_n}$ denote the intersection
$^{\db^1,\db^2}\H^{\Qb_n}\cap {}_x\H^{\Qb_n}$. 
In view of \lemref{preimage is contained in support}, it suffices
to check that for any fixed $x\in X$, the sum of the dimensions
of fibers of
$$\hl{}^{\Qb_n}:^{\db^1,\db^2}_x\H^{\Qb_n}\to \Qb_n \text{ and }
\hr{}^{\Qb_n}:^{\db^1,\db^2}_x\H^{\Qb_n}\to \Qb_n$$
does not exceed $n-1$.

\medskip

For fixed $\db^1,\db^2$, let $k$ be the first integer for
which $d^2_k=d^1_k+1$. We claim that the dimension of the fibers 
of $^{\db^1,\db^2}_x\H^{\Qb_n}$ under $\hl{}^{\Qb_n}$ are exactly
$k-1$, and those for $\hr{}^{\Qb_n}$ are of dimension $n-k$.

Indeed, let first $(\M,\kappa_1,...,\kappa_n)$ be a point of
$\Qb^d_n$ such that each $\kappa_i$ has a zero of order $d^1_i$
at $x$. Then on the formal disk around $x$ we have a filtration
$$0=\M_0\subset \M_1\subset...\subset \M_n=\M$$
with $\M_i/\M_{i-1}\simeq \Omega^{n-1+...n-i}
\left((d^1_i-d^1_{i-1})(x)\right)$.

The variety of all possible upper modifications $\M'$ of $\M$ at the
given $x$ is the projective space $\PP(\M_x)$. Now, the condition that
the point that $\M'$ defines in $_x\H^{\Qb_n}$ belongs to
$^{\db^1,\db^2}_x\H^{\Qb_n}$ with the above condition on $(\db^1,\db^2)$ 
means that the corresponding line
$\ell\subset \M_x$ belongs to $(\M_k)_x\subset \M_x$, and does not
belong to $(\M_{k-1})_x$.

The dimension of the variety of these lines is exactly $k-1$.

\medskip

Similarly, if we start with a point 
$(\M',\kappa'_1,...,\kappa'_n)\in \Qb^{d+1}_n$ with each 
$\kappa'_i$ having a zero of order $d^2_i$ at $x$, we obtain a flag
$$0=\M'_0\subset \M'_1\subset...\subset \M'_n=\M'$$
defined on the formal disk around $x$, and
$\M'_i/\M'_{i-1}\simeq \Omega^{n-1+...n-i}
\left((d^2_i-d^2_{i-1})(x)\right)$.

The variety of all possible lower modifications $\M$ of $\M'$
constitutes the projective space of hyperplanes in $\M'_x$.
The condition that $\M$ defines a point of 
$^{\db^1,\db^2}_x\H^{\Qb_n}$ means that the corresponding
hyperplane contains $(\M'_{k-1})_x$, and does not contain
$(\M'_k)_x$, and the variety of these hyperplanes has 
dimension $n-k$.

\ssec{}

As usual, everything said in this section carries over to the
relative situation, i.e., for a base $S$ we have the Hecke functors
$\on{H}^{\Qb_k}:\on{D}^W(S\times \Qb_k) \to \on{D}^W(S\times X\times \Qb_k)$.
Moreover, for $k=n$ this functor is right-exact.

Note, however, that for a map $g:S_1\to S_2$, the functors
$\on{H}^{\Qb_k}$ commute only with the $!$-push forward and
the $*$-pull back. 

\section{Construction of quotients}   \label{construction of quotients}

In this section we will complete the construction of the quotient categories.
Recall the category $\Dt(\Bun_n')$ introduced in \secref{degenerate sheaves on Bun'_n}.
A naive idea would be to define $\Dt(\Bun_n)$ as a quotient of $\on{D}(\Bun_n)$ by
the kernel of the composition
$$\on{D}(\Bun_n)\overset{\pi^*}\to \on{D}(\Bun'_n)\to \Dt(\Bun_n),$$
i.e., to ``kill'' those sheaves $\F$ on $\Bun_n$, for which $\pi^*(\F)\in \on{D}(\Bun'_n)$
is degenerate.
However, this definition does not work, because since the map $\pi:\Bun'_n\to \Bun_n$
is not smooth, the functor $\pi^*$ is not exact, and the resulting kernel would
not in general be compatible with the t-structure.
To remedy this, we will ``kill'' even more objects in $\on{D}(\Bun_n)$.

\ssec{}   \label{beginning of constr}

Let $\U\subset \Bun_n$ be the open substack corresponding
to $\M\in \Bun_n$ for which $Ext^1(\Omega^{n-1},\M)=0$.
It is well-known that each $\U\cap\Bun_n^d$ is of finite type.
Obviously, the map $\pi:\Bun_n'\to \Bun_n$ is smooth over $\U$.
Set $\V=\Bun_n-\U$, $\U^d=\U\cap \Bun^d_n$, and $\V^d=\Bun^d_n-\U^d$.

Recall (cf. \cite{FGV1}, Sect. 3.2) that a vector bundle $\M$ is called
very unstable if $\M$ can represented as a direct sum $\M=\M^1\oplus \M^2$,
with $\M^i\neq 0$, and $Ext^1(\M_1,\M_2)=0$.

It is well-known (cf. \cite{FGV1}, Lemma 6.11) that if $\F$ is a 
cuspidal object of $\on{D}(\Bun_n)$, then its $*$-stalk at every
very unstable point $\M\in \Bun_n$ vanishes. The following is also
well-known (cf. \cite{FGV1}, Lemma 3.3):

\begin{lem}  \label{support of cuspidal}
There exists an integer $d_0$, depending only on the genus of $X$,
such that for $d\geq d_0$ every point of $\M\in \V^d$ is
very unstable.
\end{lem}

\ssec{}

Set $\V'\subset \Bun_n'$, $\U'\subset \Bun_n'$ to be the preimages of
$\V$ and $\U$, respectively, in $\Bun'_n$. We denote by $\jmath:\U\to \Bun_n$,
$\jmath':\U'\to \Bun'_n$ the corresponding open embeddings.

The category $\on{D}(\V')$ is a full triangulated subcategory
of $\on{D}(\Bun'_n)$. It is compatible with the t-structure
on $\on{D}(\Bun'_n)$, cf. \secref{subcategory compatible with t-structure}.

Recall now the subcategory
$\on{D}_{degen}(\Bun_n')\subset \on{D}(\Bun'_n)$ of 
\secref{degenerate sheaves on Bun'_n}, which by definition consists 
of objects annihilated by the functor
$W:\on{D}(\Bun'_n)\to \on{D}^W(\Qb)$. Since the functor $W$ is
exact, $\on{D}_{degen}(\Bun_n')$ is also compatible with the t-structure
on $\on{D}(\Bun'_n)$, cf. \lemref{kernel is compatible}.

Let $\on{D}(\V'+\on{degen})\subset \on{D}(\Bun'_n)$ be the 
triangulated category generated by $\on{D}(\V')$ and
$\on{D}_{degen}(\Bun_n')$, i.e.,
$\on{D}(\V'+\on{degen})$ is the minimal
full triangulated subcategory of $\on{D}(\Bun'_n)$, which contains
both $\on{D}(\V')$ and $\on{D}_{degen}(\Bun_n')$.

We have:
\begin{lem}
Let $\C$ be a triangulated subcategory endowed with a t-structure
and let $\C',\C''\subset \C$ be two full triangulated subcategories,
both compatible with the t-structure on $\C$. Let $\C'+\C''\subset \C$ be the 
triangulated subcategory generated by $\C'$ and $\C''$.
Then $\C'+\C''$ is also compatible with the t-structure on $\C$.
\end{lem}

\begin{proof}

By definition, $(\C'+\C'')\cap \on{P}(\C)$ is the full abelian
subcategory of $\on{P}(\C)$, consisting of objects, 
which admit a finite filtration
with successive quotients being objects of either $\on{P}(\C')$ or 
$\on{P}(\C'')$. Clearly, $(\C'+\C'')\cap \on{P}(\C)$ is a Serre
subcategory of $\on{P}(\C)$.

Thus, we have to show that if $\S$ is an object of $\C'+\C''$,
then so is $\tau^{\leq 0}(\S)$. Suppose that $\S$ can be obtained
by an iterated $i$-fold procedure of taking cones, starting
from objects of either $\C'$ or $\C''$. By induction on $i$,
we may assume that $\S$ fits into an exact triangle
$$\S_1\to \S\to \S_2$$
with $\S_1,\S_2\in \C'+\C''$ and $\tau^{\leq 0}(\S_1),\tau^{\leq 0}(\S_2)$
being also in $\C'+\C''$. Let $\S_3$ be the image of $h^0(\S_2)$
in $h^1(\S_1)$; it belongs to $(\C'+\C'')\cap \on{P}(\C)$, by the above.
Let $\S_4$ be the cone of $\tau^{\leq 0}(\S_2)\to \S_3$.
Then $\tau^{\leq 0}(\S)$ fits into the exact triangle
$$\tau^{\leq 0}(\S_1)\to\tau^{\leq 0}(\S)\to \S_4.$$

\end{proof}

By applying this lemma to $\on{D}(\V'+\on{degen})$, we obtain
from \propref{t-structure on quotient} that the quotient triangulated category
$$\wt{\Dt}(\Bun_n'):=\on{D}(\Bun'_n)/\on{D}(\V'+\on{degen})$$
carries a t-structure.

\medskip

For an arbitrary base scheme $S$, the category 
$\wt{\Dt}(S\times \Bun_n')$ is defined in an absolutely similar way, 
as a quotient of $\on{D}(S\times \Bun_n')$ by a subcategory
denoted $\on{D}(S,\V'+\on{degen})$.
Thus quotient is stable under the standard functors, i.e.,
for a map $S_1\to S_2$ the four functors 
$\on{D}(S_1\times \Bun'_n)\rightleftarrows 
\on{D}(S_2\times \Bun'_n)$ give rise to well-defined functors
on the quotients
$\wt{\Dt}(S_1\times \Bun_n') \rightleftarrows\wt{\Dt}(S_2\times \Bun_n')$.

Moreover, the Verdier duality functor on $\on{D}(S\times \Bun'_n)$ descends
to a well-defined self-functor on $\wt{\Dt}(S\times \Bun_n')$.
Finally, the ``tensor product along $S$'' functor
$$\on{D}(S)\times \on{D}(S\times \Bun'_n)\to \on{D}(S\times \Bun'_n)$$
is also well-defined on the quotient.

\ssec{}

We define the functor
$\ol{\pi^*_S}:\on{D}(S\times \Bun_n)\to \wt{\Dt}(S\times \Bun_n')$
as follows.

For $\F\in \on{D}(S\times \Bun^d_n)$ we set $\ol{\pi^*_S}(\F)$
to be the image of $(\on{id}\times\pi)^*(\F)[\on{dim}(d)]$ under
$\on{D}(S\times \Bun'_n)\to \wt{\Dt}(S\times \Bun_n')$,
where $\on{dim}(d)=\on{dim.rel.}(\U',\U^d)$.

Note that $\on{dim}(d+1)=\on{dim}(d)+1$, by the Riemann-Roch
theorem.

\begin{prop}  \label{pull-back and projection}
The functor $\ol{\pi^*_S}$ is exact.
Moreover, it commutes with the Verdier duality, the tensor
product along $S$, and for a map $S_1\to S_2$ it is 
compatible with the four functors 
$\on{D}(S_1\times \Bun_n)\rightleftarrows 
\on{D}(S_2\times \Bun_n)$ and
$\wt{\Dt}(S_1\times \Bun_n') \rightleftarrows\wt{\Dt}(S_2\times \Bun_n')$.
\end{prop}

\begin{proof}

The functor $\F\mapsto (\on{id}\times \pi)^*(\F)[\on{dim}(d)]$
from $\on{D}(S\times \Bun^d_n)$ to $\on{D}(S\times \Bun_n')$
is not exact, because the map $\pi$ is not smooth.

However, for a perverse sheaf $\F\in \on{P}(S\times \Bun^d_n)$
all the non-zero cohomology sheaves of 
$(\on{id}\times\pi)^*(\F)[\on{dim}(d)]$
are supported on $\V'$. Hence they vanish after the projection
to $\wt{\Dt}(S\times \Bun_n')$.
This establishes the exactness of $\ol{\pi^*_S}$.

The other assertions of the proposition follow in a similar way.
For example, to show that $\ol{\pi^*_S}$ commutes with the Verdier
duality functor, it suffices to observe that
$$\jmath'{}^*\circ\left(\on{id}\times \pi)^*(\F)[\on{dim}(d)]\right)\simeq
\jmath'{}^*\circ \DD\circ \left(\on{id}\times \pi)^*(\DD\F)[\on{dim}(d)]\right),$$
and for any $\F'\in \on{D}(S\times \Bun_n')$ the map
$\jmath'_!\circ \jmath'{}^*(\F')\to \F'$ becomes an isomorphism in
$\wt{\Dt}(S\times \Bun_n')$.

\end{proof}

\ssec{}

Since the functor $\ol{\pi^*_S}$ is exact, the subcategory
$$\on{D}_{degen}(S\times \Bun_n):=
ker(\ol{\pi^*_S})\subset \on{D}(S\times \Bun_n)$$ 
is compatible with the t-structure.

We define the category $\Dt(S\times \Bun_n)$ 
as the quotient $\on{D}(S\times \Bun_n)/\on{D}_{degen}(S\times \Bun_n)$.
By \propref{t-structure on quotient}, $\Dt(S\times \Bun_n)$ inherits a t-structure from
$\on{D}(S\times \Bun_n)$.


\medskip

By \propref{pull-back and projection}, the standard 6 functors
that act on $\on{D}(S\times \Bun_n)$ are well-defined on the
quotient $\Dt(S\times \Bun_n)$.

Thus, it remains to show that $\Dt(S\times \Bun_n)$ satisfies
Properties 1 and 2 of \secref{properties}.

\ssec{Verification of Property 1}

We must show that the Hecke functor
$$\on{H}_S:\on{D}(S\times \Bun_n)\to \on{D}(S\times X\times \Bun_n)$$
descends to the quotient $\Dt(S\times \Bun_n)$, and the corresponding
functor $\wt{\on{H}}_S$ is exact.

\medskip

To prove the fact that it is well-defined, we must show that
$\on{H}_S$ maps $ker(\ol{\pi^*_S})$ to 
$ker(\ol{\pi^*_{S\times X}})$.
Using \eqref{pull back and Hecke}, cf. \secref{intr Hecke},
this reduces to showing that the subcategory 
$\on{D}(S,\V'+\on{degen})\subset \on{D}(S\times \Bun_n')$
is preserved by 
$\on{H}_S^{\Bun_n'}:
\on{D}(S\times \Bun_n')\to \on{D}(S\times X\times \Bun_n')$.

For that, it suffices to show that $\on{H}_S^{\Bun_n'}$ maps
$\on{D}_{degen}(S\times \Bun_n')$ to
$\on{D}_{degen}(S\times X\times \Bun_n')$ and
$\on{D}(S\times \V')$ to $\on{D}(S\times X\times \V')$.

The former follows immediately from \propref{Hecke and Whittaker commute}.
To prove the latter, it suffices to observe that in the diagram
$$\Bun_n'  \overset{\hl{}^{\Bun_n'}}\longleftarrow \H^{\Bun_n'}
\overset{\hr{}^{\Bun_n'}}\longrightarrow \Bun_n',$$
the subset $(\hr{}^{\Bun_n'})^{-1}(\V')$ is contained in 
$(\hl{}^{\Bun_n'})^{-1}(\V')$.

\medskip

Now we will prove the exactness of $\wt{\on{H}}_S$ on
$\Dt(S\times \Bun_n)$. Since the functor 
$\on{H}_S:\on{D}(S\times \Bun_n)\to \on{D}(S\times X\times \Bun_n)$
commutes with the Verdier duality, it suffices to show that
$\wt{\on{H}}_S$ is right-exact on $\Dt(S\times \Bun_n)$.

We have the following general assertion:

\medskip

Let
$$
\CD
\C_1  @>F>> \C_2  \\
@V{G_1}VV   @V{G_2}VV  \\
\C'_1 @>{F'}>> \C'_2
\endCD
$$
be a commutative diagram of functors
between triangulated categories endowed
with t-structures. Suppose that the functors $F$, and $F'$ are exact,
and the functor $G_2$ is right-exact (resp., left-exact, exact).

\begin{lem}  \label{general exactness}
Under the above circumstances, $G_1$ gives rise
to a well-defined functor 
$$\wt{G_1}:\C_1/ker(F)\to \C'_1/ker(F'),$$ 
and the latter functor is right-exact (resp., left-exact, exact).
\end{lem}

\begin{proof}

The fact that the functor $\wt{G_1}:\C_1/ker(F)\to \C'_1/ker(F')$
is well-defined is immediate. Let us assume that $G_2$ is
right-exact. To prove that $\wt{G_1}$ is then also right-exact,
we must show that for
$\S\in \C_1^{\leq 0}$, the projection to $\C'_1/ker(F')$ of
$\tau^{>0}(G_1(\S))$ vanishes. This amounts to showing that
$F'\left(\tau^{>0}(G_1(\S))\right)=0$. Since $F$ is exact,
$$F'\left(\tau^{>0}(G_1(\S))\right)\simeq 
\tau^{>0}\left(F'\circ G_1(\S)\right),$$
which, in turn, is isomorphic to $\tau^{>0}\left(G_2\circ F(\S)\right)$.
Since $F$ is exact, $F(\S)\in \C_2^{\leq 0}$, and since $G_2$
is right-exact, $G_2\circ F(\S)\in \C_2'{}^{\leq 0}$, which is
what we had to show.

\end{proof}

We apply this lemma first to
$\C_1=\on{D}(S\times \Bun'_n)$,
$\C'_1=\on{D}(S\times X\times \Bun'_n)$,
$\C_2=\on{D}(S\times \Qb)$, 
$\C'_2=\on{D}(S\times X\times\Qb)$,
$F,F'=W$, $G_1=\on{H}_S^{\Bun_n'}$, and
$G_2=\on{H}_S^{\Qb}$.

\medskip

From \thmref{Hecke is right-exact} we know that 
$\on{H}^{\Qb}_S$ is right exact, which by \lemref{general exactness}
implies that $\wt{\on{H}}_S^{\Bun_n'}:\Dt(S\times \Bun'_n)\to 
\Dt(S\times X\times \Bun'_n)$ is right exact. Hence, the corresponding
functor $\wt{\wt{\on{H}}}_S^{\Bun_n'}:
\wt{\Dt}(S\times \Bun'_n)\to \wt{\Dt}(S\times X\times \Bun'_n)$
is also right-exact.

\medskip

We apply \lemref{general exactness} the second time
to $\C_1=\on{D}(S\times \Bun_n)$, $\C'_1=\on{D}(S\times X\times \Bun_n)$,
$\C_2=\wt{\Dt}(S\times \Bun'_n)$, $\C'_2=\wt{\Dt}(S\times X\times \Bun'_n)$,
$F=\ol{\pi^*_S}$, $F'=\ol{\pi^*_{S\times X}}$, $G_1=\on{H}_S$, and
$G_2=\wt{\wt{\on{H}}}_S^{\Bun_n'}$.

We conclude that $\wt{\on{H}}_S$ is exact as a functor
$\Dt(S\times \Bun_n)\to \Dt(S\times X\times \Bun_n)$.

\ssec{Verification of property 2}  \label{verification of prop 2}

We must show that if $\F_1$ is a cuspidal object of
$\on{D}(\Bun^d_n)$ with $d\geq d_0$ (cf. \lemref{support of cuspidal})
and $\F_2\in \on{D}_{degen}(\Bun_n)$, then 
$$\Hom_{\on{D}(\Bun_n)}(\F_1,\F_2)=0.$$

First, from \lemref{support of cuspidal}, we obtain that
$\F_1\simeq \jmath_!(\jmath^*(\F_1))$. Therefore,
$$\Hom_{\on{D}(\Bun_n)}(\F_1,\F_2)\simeq
\Hom_{\on{D}(\U)}(\F_1|_{\U},\F_2|_{\U}).$$

Consider the natural map
\begin{equation} \label{expression for Hom}
\Hom_{\on{D}(\U)}(\F_1|_{\U},\F_2|_{\U})\to 
\Hom_{\on{D}^{\GG_m}(\U')}\left(\pi^*(\F_1)|_{\U'},\pi^*(\F_2)|_{\U'}\right).
\end{equation}

We claim that this map is an injection. Indeed, the quotient stack
$\U'/\GG_m$ is fibration into projective spaces over $\U$, and 
the required injectivity follows from the 
fact the direct image of the constant 
sheaf from $\U'/\GG_m$ to $\U$ contains the constant sheaf on $\U$
as a direct summand.

Thus, it will be sufficient to show that the right-hand side of
\eqref{expression for Hom} vanishes. Note that since
$\pi^*(\F_1)|_{\V'}=0$, we can rewrite \eqref{expression for Hom} as
$$\Hom_{\on{D}^{\GG_m}(\Bun_n')}
\left(\pi^*(\F_1),\pi^*(\F_2)\right).$$

We will show that for any 
$\F'_2\in \on{D}(\V'+\on{degen})$,
$$\Hom_{\on{D}^{\GG_m}(\Bun_n')}(\pi^*(\F_1),\F'_2)=0.$$

By definition of $\on{D}(\V'+\on{degen})$, we must analyze two cases:

\noindent Case 1: $\F'_2\in \on{D}(\V')$. 
In this case the above $\Hom$ vanishes, because $\pi^*(\F_1)|_{\V'}=0$,
as was noticed before.

\medskip

\noindent Case 2: $\F'_2\in \on{D}_{degen}(\Bun'_n)$.

We know that $\F'_1:=\pi^*(\F_1)$ is cuspidal, and from
\thmref{cuspidal orthogonal to degenerate} (or rather from its 
$\GG_m$-equivariant version) we obtain that
$$\Hom_{\on{D}^{\GG_m}(\Bun_n')}(\F'_1,\F'_2)\simeq
\Hom_{\on{D}^{\GG_m}(\Qb_n)}(W(\F'_1),W(\F'_2))=0,$$
since it was assumed that $W(\F'_2)=0$.

\section*{Appendix}   \label{App}

\ssec*{A.1.}

We will present now a different way of how \conjref{vanishing conjecture}
can be deduced from \thmref{exactness of Av}. This argument  
is due to A.~Braverman. 

\medskip

By induction, we are assuming \conjref{vanishing conjecture} 
for all integers $n'<n$. 
It is enough to show that $\on{Av}_E^d(\F)$ vanishes
for a perverse sheaf $\F\in \on{D}(\Bun_n)$, where $d$ is as in 
\conjref{vanishing conjecture}. We know that $\on{Av}_E^d(\F)$
is a perverse sheaf (by \thmref{exactness of Av}) and that it
is cuspidal, by \lemref{Av is cuspidal}.

Recall the functor $\on{Av}^{-d}_{E^*}$, which is left and
right adjoint to $\on{Av}_E^d$. Since
$$\on{Hom}_{\on{D}(\Bun_n)}(\on{Av}_E^d(\F),\on{Av}_E^d(\F))\simeq
\on{Hom}_{\on{D}(\Bun_n)}(\on{Av}^{-d}_{E^*}(\on{Av}_E^d(\F)),\F),$$
we obtain that it is enough to show that the functor
$\on{Av}^{-d}_{E^*}$ annihilates every cuspidal perverse sheaf.

The stack $\Bun_n$ admits a natural automorphism, which sends a
bundle to its dual. This automorphism transforms the functor 
$\on{Av}^{-d}_{E^*}$ to $\on{Av}_{E^*}^d$. Since $E^*$ is irreducible
if and only if $E$ is, we deduce that it is enough to show that 
$\on{Av}_E^d(\F)=0$, where $\F$ is both perverse and cuspidal.

By \lemref{van Euler char} and \thmref{exactness of Av}, the above vanishing
is equivalent to a weaker statement. Namely, it is sufficient
to show that for a cuspidal perverse sheaf $\F$, the 
Euler-Poincar\'e characteristics of the stalks of $\on{Av}_E^d(\F)$ is zero.
Finally, by \lemref{independence} we conclude that it is enough to show that 
the Euler-Poincar\'e characteristics of the stalks of $\on{Av}_{E_0}^d(\F)$
vanish, where $E_0$ is {\it the trivial local system} of rank equal to
that of $E$, and $\F\in \on{D}(\Bun_n)$ is cuspidal and perverse.

We will prove a stronger statement. Namely, we will show that the object
$\on{Av}_{E_0}^d(\F)$ vanishes, where $E_0$ is trivial local system of rank $m$, and 
$d>(2g-2)\cdot n\cdot m$ for {\it every cuspidal object} $\F\in \on{D}(\Bun_n)$.

\ssec*{A.2}

First, let us express the functor $\on{Av}_{E_0}^d$ in terms of the corresponding
averaging functor for the trivial $1$-dimensional local system. We have:

\medskip

\noindent{\bf Proposition A.3.} {\it
Let a local system $E$ be the direct sum $E=E_1\oplus E_2$. Then we have canonically:
$$\on{Av}_{E}^d(\F)\simeq \oplus \on{Av}_{E_1}^{d_1}\circ \on{Av}_{E_2}^{d_2}(\F),$$
where the direct sum is taken over all pairs $(d_1,d_2)$ with $d_i\geq 0$, $d_1+d_2=0$.}

\begin{proof}

For two non-negative integers $d_1,d_2$ 
consider the stack $\Mod^{d_1}_n\underset{\Bun_n}\times \Mod^{d_2}_n$,
where the fiber product is formed using the maps $\hr:\Mod^{d_1}_n\to \Bun_n$
and $\hl:\Mod^{d_2}_n\to \Bun_n$. In other words, this stack classifies successive
extensions $\M\subset \M'\subset \M''$, where $\M'/\M$ is of length $d_1$ and 
$\M''/\M'$ is of length $d_2$. There is a natural projection 
$\r_{d_1,d_2}:\Mod^{d_1}_n\underset{\Bun_n}\times \Mod^{d_2}_n\to \Mod^d_n$,
where $d=d_1+d_2$. 

We have:
$$\s^*(\L^d_E)\simeq \underset{(d_1,d_2)}\oplus \, 
\r_{d_1,d_2}{}_!(\s^*(\L^{d_1}_{E_1})\boxtimes \s^*(\L^{d_2}_{E_2})).$$
Indeed, the isomorphism is evident over the open substack
$\overset{\circ}\Mod{}^d_n$, and it extends to the entire $\Mod^d_n$, since the
maps $\r_{d_1,d_2}$ are small.

By definition, this implies the required property of the functor $\on{Av}_{E}^d$.

\end{proof}

The same proof shows that the functors $\on{Av}_{E_1}^{d_1}$ and 
$\on{Av}_{E_2}^{d_2}$ mutually commute. 

Let $\on{Av}^d$, with the subscript omitted, denote the averaging functor
with respect to the trivial $1$-dimensional local system. Note that for the trivial 
$1$-dimensional local system, Laumon's sheaf on $\on{Coh}^d_0$ is the constant sheaf. 
Therefore, the functor $\on{Av}^d$ is just
\begin{equation}  \label{simple averaging}
\F\mapsto \hl_!\circ \hr{}^!(\F)[nd].
\end{equation}

From Proposition A.3 we obtain that
$$\on{Av}_{E}^d(\F)\simeq \underset{\db}\oplus\, \on{Av}^{d_1}\circ...\circ  \on{Av}^{d_m}(\F),$$
where the direct sum is taken over the set of $m$-tuples of non-negative integers
$\db=(d_1,...,d_m)$ with $d_1+...+d_m=d$. If $d>(2g-2)\cdot m\cdot n$, then for
every such $\db$ at least one $d_i$ satisfies $d_i>(2g-2)\cdot n$. Hence, we are reduced
to showing the following:

\medskip

\noindent{\bf Theroem A.4.}  
{\it If $\F\in \on{D}(\Bun_n)$ is cuspidal, then $\on{Av}^d(\F)=0$
for $d>(2g-2)\cdot n$.}

\medskip

This theorem is a geometric analog of the classical statement that the L-function of
a cuspidal automorphic representation of $GL_n$ over a function field is a polynomial.
The proof will be a geometrization of the Jacquet-Godement proof of the above
classical fact, in the spirit of how the functional equation is established for
geometric Eisenstein series in \cite{BG}, Sect. 7.3. 

\ssec*{A.5.}

The starting point is the following observation, due to V.~Drinfeld and proved
in \cite{BG}, Sect. 7.3.:

Let $\Y$ be a stack and $\E_1,\E_2$ two vector bundles on it, and 
$p:\E_1\to \E_2$ a map between them as coherent sheaves. 
Let $K_{p}$ be the kernel of $p$, considered as a group-scheme over $\Y$
and $\varphi$ be its projection onto $\Y$. Consider the object $\K_p$ of
$\on{D}(\Y)$ equal to $\varphi_!({\Ql}_{K_{p}})[\on{dim.rel.}(\E_1,\Y)]$, where ${\Ql}_{K_{p}}$
denotes the constant sheaf on $K_p$.
Let $\check p:\check E_2\to \check E_1$ denote the dual map, and consider also
the object $\K_{\check p}:=\check\varphi_!({\Ql}_{K_{\check p}})[\on{dim.rel.}(\check\E_2,\Y)]$.
We have:

\medskip

\noindent{\bf Lemma A.6.}
{\it There is a canonical isomorphism $\K_p\simeq \K_{\check p}$.}

\medskip

We will apply this lemma in the following situation. Let
$\F$ be a cuspidal object of $\on{D}(\Bun_n)$ supported on a
connected component $\Bun_n^{d'}$. As in \secref{beginning of constr},
we can assume that $\F$ is the extension by zero from an open
substack of finite type $U'\subset \Bun_n^{d'}$. Let $U$ be a
scheme of finite type, which maps smoothly to $\Bun_n^{d'-d}$;
moreover, we can assume that $U'$ was chosen large enough so that
the image of $\hr:U\underset{\Bun_n}\times \Mod^d_n\to \Bun_n^{d'}$ is
contained in $U'$. We are going to show that $\on{Av}^d(\F)|_{U}$
vanishes.

\medskip

We set the base $\Y$ to be $U\times U'$. To define $\E_1$ and $\E_2$
we pick an arbitrary point $y\in X$ and let $i$ be a large enough
integer so that $\on{Ext}^1(\M,\M'(i\cdot y))=0$, whenever 
$(\M,\M')\in \Bun_n\times \Bun_n$ is in the image of $U\times U'$.
We set $\E_1$ (resp., $\E_2$) to be the vector bundle, whose fiber
at a point of $U\times U'$ mapping to a point $(\M,\M')$ as above
is $\on{Hom}(\M,\M'(i\cdot y))$ (resp., $\on{Hom}(\M,\M'(i\cdot y)/\M')$).

The group-scheme $K_p$ has as its fiber over $(\M,\M')$
the vector space $\on{Hom}(\M,\M')$.
By Serre's duality, the fiber of $K_{\check p}$ is 
$\on{Hom}(\M',\M\otimes \Omega)$. Let $\overset{\circ}K_p$ 
(resp., $\overset{\circ}K_{\check p}$) be the open subscheme
corresponding to the condition that the map of sheaves $\M\to \M'$
(resp., $\M'\to \M\otimes \Omega$) is injective. Note that if
$d=\on{deg}(\M')-\on{deg}(\M)>(2g-2)\cdot n$, then 
$\overset{\circ}K_{\check p}$ is empty.
Let $\overset{\circ}\K_p$ be $\varphi_!({\Ql}_{\overset{\circ}K_{p}})[\on{dim.rel.}(\E_1,\Y)]$
(resp.,  $\overset{\circ}\K_{\check p}=
\varphi_!({\Ql}_{\overset{\circ}K_{\check p}})[\on{dim.rel.}(\check \E_2,\Y)]$).
Finally, let $\K^c_p$ (resp., $\K^c_{\check p}$) denote the cone of the natural
arrow $\overset{\circ}\K_p\to \K_p$ (resp., $\overset{\circ}\K_{\check p}\to 
\K_{\check p}$).

\medskip

Let us denote by $q,q'$ the projections from $U\times U'$ to $U$ and $U'$,
respectively. Consider the two functors $\on{D}(U')\to \on{D}^-(U)$ defined
by 
$$\F\mapsto q_!(q'{}^*(\F)\otimes \K_p) \text{ and }
q_!(q'{}^*(\F)\otimes \overset{\circ}\K_p).$$
Here $\on{D}^-(U)$ denotes the bounded from above derived category of sheaves on $U$,
which appears due to the fact that the map $q$ is not representable.
Note, however, that because of \eqref{simple averaging},
$$q_!(q'{}^*(\F)\otimes \overset{\circ}\K_p)\simeq \on{Av}^d(\F)|_{U}.$$

Taking into account Lemma A.6, we have reduced 
Theorem A.4 to the fact that the functors
$$\F\mapsto q_!(q'{}^*(\F)\otimes \K^c_p) \text{ and }
\F\mapsto q_!(q'{}^*(\F)\otimes \K^c_{\check p})$$
annihilate cuspidal objects. We will prove it in the case
of $\K^c_p$, as the other assertion is completely analogous.

\ssec*{A.7.}

Let $K^c_p$ denote the complement to $\overset{\circ}K_p$ in $K_p$.
By definition, it can be decomposed into the union of $n$ locally
closed substacks, where the $k$-th substack, classifies the
data of a pair of points $(u,u')\in U\times U'$ and a map
between the corresponding sheaves $\M\to \M'$, which is of
generic rank $k$, with $k$ running from $0$ to $n-1$.
Each such substack admits a further decomposition
into locally closed substacks according to the length of the
torsion of the quotient $\M'/\M$. 

It is enough to show that the correspondence $\on{D}(U')\to \on{D}^-(U)$
defined by the constant sheaf on each of these locally closed substacks
annihilates $\F\in \on{D}(U')$, provided that $\F$ is cuspidal.

\medskip

Let us consider separately the cases when
$k=0$ and when $k>0$. In the former case, the corresponding
(closed) substack of $K^c_p$ is the zero-section, i.e.,
the product $U\times U'$. Thus, we must show that
$H_c(U',\F)=0$, when $\F$ is cuspidal. In other words, we must
show that $\on{Hom}_{\on{D(U)}}(\F,{\Ql}_U)=0$. However, this follows
from \secref{verification of prop 2}: with no restriction
of generality we may assume that $d'\geq d_0$, and the object
${\Ql}_{\Bun^{d'}_n}$ clearly belongs to $\on{D}_{degen}(\Bun_n)$.

\medskip

Now let us suppose that $k>0$, and consider the stack
$$Z:=\on{Fl}^n_{n-k,k}\underset{\Bun_k}\times \Mod^a_k\underset{\Bun_k}\times
\on{Fl}^n_{k,n-k},$$
where we have used the map $(\hl\times \hr):\Mod^a_k\to \Bun_k\times \Bun_k$
to define the fiber product. By definition, a point of $Z$ is the data of
$$0\to \M_{n-k}\to \M\to \M_k\to 0;\,\, \M_k\hookrightarrow \M'_k;\,\,
0\to \M'_k\to \M'\to \M'_{n-k},$$
where $\M_{n-k},\M'_{n-k}$ are vector bundles of rank $n-k$,
$\M_k,\M'_k$ are vector bundles of rank $k$, and the quotient $\M'_k/\M_k$
is of length $a$.
 
\medskip

The stack $Z$ maps to $\Bun_n\times \Bun_n$ by remembering the data of
$(\M,\M')$ and the fiber product
$_UZ_{U'}:=U\underset{\Bun_n}\times Z\underset{\Bun_n}\times U'$ is the required 
locally closed substack of $K^c_p$. By taking the constant sheaf on
$_UZ_{U'}$ we obtain a functor $\on{D}(U')\to \on{D}^-(U)$, and we have to
show that this functor annihilates every cuspidal object $\F\in \on{D}(U')$.
However, this follows by base change from the following diagram:
$$
\CD
_UZ_{U'} @>>> \on{Fl}^n_{k,n-k}\underset{\Bun_n}\times U' @>>> U' \\
@VVV  @VVV \\
U\underset{\Bun_n}\times\on{Fl}^n_{n-k,k}\underset{\Bun_k}\times \Mod^a_k @>>> \Bun_k  \\
@VVV  \\
U\underset{\Bun_n}\times \on{Fl}^n_{n-k,k} \\
@VVV \\
U.
\endCD
$$
Indeed, the functor $\on{D}(U')\to \on{D}^-(\Bun_k)$ corresponding to the
upper-right corner of the above diagram annihilates cuspidal objects,
by the definition of the constant term functor $\on{CT}^n_{k,n-k}$, because
the vertical arrow $\on{Fl}^n_{k,n-k}\underset{\Bun_n}\times U'\to \Bun_k$,
appearing in the diagram, factors as
$$\on{Fl}^n_{k,n-k}\underset{\Bun_n}\times U'\overset{\q_{k,n-k}}\longrightarrow \Bun_k\times
\Bun_{n-k}\to \Bun_k.$$

\end{document}